\documentclass[a4paper,10pt,dvipsnames]{article}
\usepackage[utf8]{inputenc}
\usepackage[english]{babel}
\usepackage{amsthm}
\usepackage{amsfonts}
\usepackage{amssymb}	
\usepackage{amsmath}
\allowdisplaybreaks
\usepackage{chemarr}
\usepackage{slashed}
\usepackage{enumerate}
\usepackage{graphicx}
\usepackage{systeme}
\usepackage{dsfont}
\usepackage{relsize}
\usepackage[margin=0.90in]{geometry}
\usepackage{float}
\usepackage{upgreek}
\usepackage{titlesec}
\usepackage{tikz}
\usepackage[labelformat=simple]{subcaption}

\usepackage{graphicx}
\usepackage{bmpsize}
\usepackage{epstopdf}
\usepackage{bbm}
\usepackage{etoolbox}
\usepackage{xcolor}
\usepackage{titling}

\usepackage{blindtext,graphicx}
\usepackage[absolute]{textpos}
\usepackage[hang,flushmargin]{footmisc}

\usepackage{xfrac}
\usepackage{nicefrac}
\usepackage{empheq}
\usepackage{stmaryrd}
\usepackage{cancel}
\usepackage[linguistics]{forest}
\usepackage{url}
\usepackage[font=small]{caption}
\usepackage{array}
\usepackage{dsfont}
\usepackage{tikz}
\usetikzlibrary{trees}
\usetikzlibrary{calc}
\usepackage{pdflscape}
\usepackage{rotating}

\usepackage{todonotes}
\usepackage{xcolor,hyperref}
\definecolor{darkblue}{rgb}{0,0,0.6}
\hypersetup{
    colorlinks=true,       
    linkcolor=darkblue,          
    citecolor=darkblue,        
    filecolor=darkblue,      
    urlcolor=darkblue           
}

\newcommand{\AS}[1]{{\color{blue}#1}}
\usepackage[giveninits=true,sortcites=true,date=year,maxbibnames=99,doi=false,isbn=false,url=false,eprint=false]{biblatex}
\renewbibmacro{in:}{}
\DeclareFieldFormat{pages}{#1}
\AtEveryBibitem{%
\clearlist{language}%
}
\usepackage{csquotes}

\usepackage{setspace}
\makeatletter
\newcommand{\thickhline}{%
\noalign {\ifnum 0=`}\fi \hrule height 1pt\section{The relaxed micromorphic model in plane-strain}

Under the plane-strain hypothesis only the in-plane components of the kinematic fields are different from zero and they only depend on $(x_1,x_2)$.
The structure of the kinematic fields ($\widetilde{u}$,$\widetilde{P}$) are \cite{iecsan2001plane}
\begin{align}
\widetilde{u}
=
\left(
\begin{array}{ccc}
u_{1} \\ 
u_{2} \\
0
\end{array}
\right)
\, ,
&&
\widetilde{u}^{\sharp}
=
\left(
\begin{array}{ccc}
u_{1} \\ 
u_{2}
\end{array}
\right)
\, ,
&&
\widetilde{P}
=
\left(
\begin{array}{ccc}
P_{11} & P_{12} & 0 \\ 
P_{21} & P_{22} & 0 \\ 
0 & 0 & 0
\end{array}
\right)
\, ,
&&
\widetilde{P}^{\sharp}
=
\left(
\begin{array}{ccc}
P_{11} & P_{12}\\
P_{21} & P_{22}
\end{array}
\right)
\, ,
\label{eq:gen_kine}
\end{align}
while the structures of the gradient of the displacement $\text{D}\widetilde{u}$, of the Curl of the micro distortion tensor $\text{Curl}\,\widetilde{P}$, and of the double Curl of the micro distortion tensor $\text{Curl}\,\text{Curl}\,\widetilde{P}$ are
\begin{gather}
\text{D}\widetilde{u}
=
\left(
\begin{array}{ccc}
u_{1,1} & u_{1,2} & 0 \\
u_{2,1} & u_{2,2} & 0 \\
0 & 0 & 0 \\
\end{array}
\right)
\, ,
\qquad
\text{D}\widetilde{u}^{\sharp}
=
\left(
\begin{array}{ccc}
u_{1,1} & u_{1,2} \\
u_{2,1} & u_{2,2}
\end{array}
\right)
\, ,
\\*[5pt]
\text{Curl} \,\widetilde{P}
= 
\left(
\begin{array}{cc|c}
0 & 0 & P_{12,1}-P_{11,2} \\
0 & 0 & P_{22,1}-P_{21,2} \\
\hline
0 & 0 & 0 \\
\end{array}
\right)
= 
\left(
\begin{array}{cc|c}
0 & 0 & \\
0 & 0 & \multicolumn{1}{c}{\smash{\raisebox{.5\normalbaselineskip}{$\text{Curl}_{\text{2D}}\,\widetilde{P}$}}} \\
\hline
0 & 0 & 0 \\
\end{array}
\right)
\, ,
\qquad
\text{Curl}_{\text{2D}} \,\widetilde{P}
= 
\left(
\begin{array}{ccc}
P_{12,1}-P_{11,2} \\
P_{22,1}-P_{21,2}
\end{array}
\right)
\in \mathbb{R}^2
\, ,
\\*[5pt]
\text{Curl} \,\text{Curl} \,\widetilde{P}
= 
\left(
\begin{array}{cc|c}
P_{12,12}-P_{11,22} &
P_{11,12}-P_{12,11} & 0 \\
P_{22,12}-P_{21,22} &
P_{21,12}-P_{22,11} & 0 \\
\hline
0 & 0 & 0 \\
\end{array}
\right)
= 
\left(
\begin{array}{cc|c}
&   & 0 \\
\multicolumn{2}{c|}{\smash{\raisebox{.5\normalbaselineskip}{$\text{Curl} \, \text{Curl}_{\text{2D}}\,\widetilde{P}$}}} & 0 \\
\hline
0 & 0 & 0 \\
\end{array}
\right)
\, ,
\\*[5pt]
\text{Curl} \, \text{Curl}_{\text{2D}}\,\widetilde{P}=
\left(
\begin{array}{cc}
P_{12,12}-P_{11,22} &
P_{11,12}-P_{12,11} \\
P_{22,12}-P_{21,22} &
P_{21,12}-P_{22,11}
\end{array}
\right)
\, .
\label{eq:gen_kine_deriv}
\end{gather}
Because of the nature of the Curl operator, it is noted that $\text{Curl} \,P$ just has out of plane components that depend on the in plane components of $P$, while $\text{Curl} \,\text{Curl} \,P$ fully preserves the in-plane structure.
Moreover, since 
\begin{gather}
\text{tr}(\text{Curl} \,\widetilde{P})=0
\\*[5pt]
\left\lVert \text{dev sym} \, \text{Curl} \, \widetilde{P}\right\rVert^2
=\left\lVert \text{sym} \, \text{Curl} \, \widetilde{P}\right\rVert^2
=\left\lVert \text{skew} \,  \text{Curl} \, \widetilde{P}\right\rVert^2
=\frac{1}{2}\left\lVert\text{Curl} \, \widetilde{P}\right\rVert^2
=\frac{1}{2}\left\lVert\text{Curl}_{\text{2D}} \, \widetilde{P}\right\rVert^2
\, ,
\notag
\end{gather}
the model will just depend on one cumulative parameter $\widetilde{a}\coloneqq\frac{a_1+a_2}{2}$, and the strain energy density in  (\ref{eq:energy_RM_3}) reduces to
\begin{align}
W (\text{D}\widetilde{u}, \widetilde{P},\text{Curl}\,\widetilde{P})
= &
\, \mu_{\rm e} \left\lVert \text{sym} (\text{D}\widetilde{u} - \widetilde{P} ) \right\rVert^{2}
+ \mu_{\rm c} \left\lVert \text{skew} (\text{D}\widetilde{u} - \widetilde{P} ) \right\rVert^{2}
+ \frac{\lambda_{\rm e}}{2} \text{tr}^2 (\text{D}\widetilde{u} - \widetilde{P} )
\label{eq:energy_RM_2D}
\\*
&
+ \mu_{\rm micro} \left\lVert \text{sym}\,\widetilde{P} \right\rVert^{2}
+ \frac{\lambda_{\rm micro}}{2} \text{tr}^2 (\widetilde{P})
+ \frac{L_{\rm c}^2 }{2} \, \widetilde{a} \,
\left\lVert \text{Curl} \, \widetilde{P}\right\rVert^2 
\, ,
\notag
\\*
= &
\, \mu_{\rm e} \left\lVert \text{sym} (\text{D}\widetilde{u}^{\sharp} - \widetilde{P}^{\sharp} ) \right\rVert^{2}
+ \mu_{\rm c} \left\lVert \text{skew} (\text{D}\widetilde{u}^{\sharp} - \widetilde{P}^{\sharp} ) \right\rVert^{2}
+ \frac{\lambda_{\rm e}}{2} \text{tr}^2 (\text{D}\widetilde{u}^{\sharp} - \widetilde{P}^{\sharp} )
\notag
\\*
&
+ \mu_{\rm micro} \left\lVert \text{sym}\,\widetilde{P}^{\sharp} \right\rVert^{2}
+ \frac{\lambda_{\rm micro}}{2} \text{tr}^2 (\widetilde{P}^{\sharp})
+ \frac{L_{\rm c}^2 }{2} \, \widetilde{a} \,
\left\lVert \text{Curl}_{\text{2D}} \, \widetilde{P}^{\sharp}\right\rVert^2 
\, .
\notag
\end{align}
The relations for the macro Lamé parameters ($\mu_{\rm macro},\lambda_{\rm macro}$) and the macroscopic bulk modulus for plane strain ($\widetilde{\kappa}_{\rm macro}$) are then
\begin{align}
\mu_{\rm macro} &\coloneqq \dfrac{\mu_{\rm e} \, \mu_{\rm micro}}{\mu_{\rm e} + \mu_{\rm micro}} \, ,
\qquad\qquad\qquad
\widetilde{\kappa}_{\rm macro} \coloneqq \dfrac{\widetilde{\kappa}_{\rm e} \, \widetilde{\kappa}_{\rm micro}}{\widetilde{\kappa}_{\rm e} + \widetilde{\kappa}_{\rm micro}} \, ,
\label{eq:Cmacro_2D}
\\*[5pt]
\widetilde{\lambda}_{\rm macro} &\coloneqq
\dfrac{(\mu_{\rm e} + \lambda_{\rm e})(\mu_{\rm micro} + \lambda_{\rm micro})}{(\mu_{\rm e} + \lambda_{\rm e}) + (\mu_{\rm micro} + \lambda_{\rm micro})} -\dfrac{\mu_{\rm e} \, \mu_{\rm micro}}{\mu_{\rm e} + \mu_{\rm micro}}
\, ,
\notag
\end{align}
where  $\widetilde{\kappa}_{\rm macro}= \mu_{\rm macro}+\widetilde{\lambda}_{\rm macro} $.
In order to have $\widetilde{\lambda}_{\rm macro}=\lambda_{\rm micro}=0$, the only possible condition is again $\lambda_{\rm micro}=\lambda_{\rm e}=0$.

Taking the first variation of the strain energy $I=\displaystyle\int_{\Omega} W \, \text{d} x$ under the plain strain hypothesis with respect to ($\widetilde{u}^{\sharp}$,$\widetilde{P}^{\sharp}$) leads to
\begin{align}
\delta I^{\widetilde{u}^{\sharp}}
&=
\int_{\Omega}
\hspace{-0.1cm}
\Big(
2\mu_{\rm e} \,
\langle \text{sym} (\text{D}\widetilde{u}^{\sharp} - \widetilde{P}^{\sharp} ),\text{D}\delta \widetilde{u}^{\sharp} \rangle
+ 2\mu_{\rm c} \,
\langle \text{skew} (\text{D}\widetilde{u}^{\sharp} - \widetilde{P}^{\sharp} ),\text{D}\delta \widetilde{u}^{\sharp} \rangle
+ \lambda_{\rm e} \langle \text{tr} (\text{D}\widetilde{u}^{\sharp} - \widetilde{P}^{\sharp} ) \mathbbm{1}_2, \text{D}\delta \widetilde{u}^{\sharp} \rangle
\Big)
\text{d}x
\, ,
\label{eq:first_varia_energy_RM_2D_u}
\\*[5pt]
\delta I^{\widetilde{P}^{\sharp}}
&=
\int_{\Omega}
\hspace{-0.1cm}
\Big(
-2\mu_{\rm e} \,
\langle \text{sym} (\text{D}\widetilde{u}^{\sharp} - \widetilde{P}^{\sharp} ),\delta \widetilde{P}^{\sharp} \rangle
-2\mu_{\rm c} \,
\langle \text{skew} (\text{D}\widetilde{u}^{\sharp} - \widetilde{P}^{\sharp} ),\delta \widetilde{P}^{\sharp} \rangle
-\lambda_{\rm e} \langle \text{tr} (\text{D}\widetilde{u}^{\sharp} - \widetilde{P}^{\sharp} ) \mathbbm{1}_2, \delta \widetilde{P}^{\sharp} \rangle
\label{eq:first_varia_energy_RM_2D_P}
\\
&
\phantom{=\int_{\Omega}\Big(}
+
2\mu_{\rm micro} \langle \text{sym}\,\widetilde{P}^{\sharp},\delta \widetilde{P}^{\sharp} \rangle
+ \lambda_{\rm micro} \langle \text{tr} (\widetilde{P}^{\sharp} ) \mathbbm{1}_2, \delta \widetilde{P}^{\sharp} \rangle
+ L_{\rm c}^2 \, \widetilde{a} \,
\langle \text{Curl}_{\text{2D}} \, \widetilde{P}^{\sharp}, \text{Curl}_{\text{2D}} \, \delta \widetilde{P}^{\sharp} \rangle
\Big)
\text{d}x
\, .
\notag
\end{align}
The equilibrium equation are now obtained by requiring 
\begin{equation}
\delta I^{\widetilde{u}^{\sharp}}= \langle \widetilde{f} , \delta \widetilde{u}^{\sharp} \rangle \, ,
\quad
\forall
\,
\delta \widetilde{u}^{\sharp}
\qquad\qquad
\text{and}
\qquad\qquad
\delta I^{\widetilde{P}^{\sharp}}= \langle \widetilde{M} , \delta \widetilde{P}^{\sharp} \rangle \, ,
\quad
\forall
\,
\delta \widetilde{P}^{\sharp}  \, .
\label{eq:first_variation_2D}
\end{equation}
We define the following quantities
\begin{align}
\widetilde{\sigma}&\coloneqq
2\mu_{\rm e} \, \text{sym} (\text{D}\widetilde{u}^{\sharp} - \widetilde{P}^{\sharp} )
+ 2\mu_{\rm c} \, \text{skew} (\text{D}\widetilde{u}^{\sharp} - \widetilde{P}^{\sharp} )
+ \lambda_{\rm e} \, \text{tr} (\text{D}\widetilde{u}^{\sharp} - \widetilde{P}^{\sharp} ) \mathbbm{1}_2
\, ,
\notag
\\*
\widetilde{\sigma}_{\rm micro}&\coloneqq
2\mu_{\rm micro} \, \text{sym} \widetilde{P}^{\sharp}
+ \lambda_{\rm micro} \, \text{tr} (\widetilde{P}^{\sharp}) \mathbbm{1}_2
\in \mathbb{R}^{2\times2}\, ,
\\*
\widetilde{m}&\coloneqq
\text{Curl}_{\text{2D}} \widetilde{P}^{\sharp}
\in \mathbb{R}^{2}\, ,
\notag
\end{align}
where we used the tilde $\widetilde{\sigma}$ and $\widetilde{\sigma}_{\rm micro}$ to emphasize that here we are only considering the in plane components.
We can rewrite the first variation $\delta I^{\widetilde{u}}$ as
\begin{align}
\delta I^{\widetilde{u}^{\sharp}}
&=
\int_{\Omega}
\hspace{-0.1cm}
\langle
\widetilde{\sigma}
,
\text{D}\delta \widetilde{u}^{\sharp}
\rangle
\,
\text{d}x
=
\int_{\Omega}
\hspace{-0.1cm}
\text{div}
(
\widetilde{\sigma}^{\rm T}
\, \delta \widetilde{u}^{\sharp}
)
-
\langle 
\text{Div} \, 
\widetilde{\sigma}
,
\delta \widetilde{u}^{\sharp}
\rangle
\,
\text{d}x
=
\int_{\partial \Omega}
\hspace{-0.1cm}
\langle
\widetilde{\sigma}^{\rm T}
\, \delta \widetilde{u}^{\sharp}
,
n
\rangle
\,
\text{d}s
-
\int_{\Omega}
\hspace{-0.1cm}
\langle 
\text{Div} \, 
\widetilde{\sigma}
,
\delta \widetilde{u}^{\sharp}
\rangle
\,
\text{d}x
\label{eq:first_varia_energy_RM_2D_u_2}
\\*[5pt]
&=
\int_{\partial \Omega}
\hspace{-0.1cm}
\langle
\widetilde{\sigma}
\, n
,
\delta \widetilde{u}^{\sharp}
\rangle
\,
\text{d}s
-
\int_{\Omega}
\hspace{-0.1cm}
\langle 
\text{Div} \, 
\widetilde{\sigma}
,
\delta \widetilde{u}^{\sharp}
\rangle
\,
\text{d}x \, ,
\notag
\end{align}
and that, because of the equation (\ref{eq:first_variation_2D}), and highlighting that $\widetilde{u}$ is orthogonal with respect to the out-of-plane displacement, implies
\begin{align}
\text{Div} \, 
\widetilde{\sigma}
= \widetilde{f}
\quad
\text{in }
\,
\Omega \, ,
\qquad\qquad
\widetilde{\sigma} \, n = 0
\quad
\text{on }
\,
\partial \Omega \, .
\label{eq:equa_var_u}
\end{align}
where the out-of-plane components of $\text{Div}\, \widetilde{\sigma}$ and $\widetilde{\sigma} \, n$ must not be considered.
We can now rewrite the first variation $\delta I^{\widetilde{P}}$ as
\begin{align}
\delta I^{\widetilde{P}^{\sharp}}
&=
\int_{\Omega}
\hspace{-0.1cm}
-
\langle \widetilde{\sigma},\delta \widetilde{P}^{\sharp} \rangle
+
\langle \widetilde{\sigma}_{\rm micro}, \delta \widetilde{P}^{\sharp} \rangle
+
\langle L_{\rm c}^2\widetilde{a}\,\text{Curl}_{\text{2D}}\,\widetilde{P}^{\sharp}, \text{Curl}_{\text{2D}} \, \delta \widetilde{P}^{\sharp} \rangle
\,
\text{d}x
\notag
\\*
&=
\int_{\Omega}
\hspace{-0.1cm}
\langle
-
\widetilde{\sigma}
+
\widetilde{\sigma}_{\rm micro}
,
\delta \widetilde{P}^{\sharp}
\rangle
+
\langle L_{\rm c}^2\widetilde{a}\,\text{Curl}_{\text{2D}}\,\widetilde{P}^{\sharp}, \text{Curl}_{\text{2D}} \, \delta \widetilde{P} \rangle
\,
\text{d}x
\notag
\\*
&=
\int_{\Omega}
\hspace{-0.1cm}
\langle
-
\widetilde{\sigma}
+
\widetilde{\sigma}_{\rm micro}
,
\delta \widetilde{P}^{\sharp}
\rangle
+
\langle L_{\rm c}^2\widetilde{a}\,\text{Curl} \, \text{Curl}_{\text{2D}} \, \widetilde{P}, \delta \widetilde{P}^{\sharp} \rangle
-
\text{div} [ \displaystyle\sum_{i=1}^{3} (L_{\rm c}^2\widetilde{a}\,\text{Curl}_{\text{2D}}\,\widetilde{P})_i \times (\widetilde{P}^{\sharp})_i ] 
\,
\text{d}x
\label{eq:first_varia_energy_RM_2D_P_2}
\\*
&=
\int_{\Omega}
\hspace{-0.1cm}
\langle
-
\widetilde{\sigma}
+
\widetilde{\sigma}_{\rm micro}
+
L_{\rm c}^2\widetilde{a}\,\text{Curl} \, \text{Curl}_{\text{2D}}\,\widetilde{P}^{\sharp}
,
\delta \widetilde{P}^{\sharp}
\rangle
-
\int_{\partial \Omega}
\hspace{-0.1cm}
\langle \displaystyle\sum_{i=1}^{3} (L_{\rm c}^2\widetilde{a}\,\text{Curl}_{\text{2D}}\,\widetilde{P}^{\sharp})_i \times (\widetilde{\delta P}^{\sharp})_i, n \rangle 
\,
\text{d}s
\notag
\\*
&=
\int_{\Omega}
\hspace{-0.1cm}
\langle
-
\widetilde{\sigma}
+
\widetilde{\sigma}_{\rm micro}
+
L_{\rm c}^2\widetilde{a}\,\text{Curl} \, \text{Curl}_{\text{2D}}\,\widetilde{P}^{\sharp}
,
\delta \widetilde{P}^{\sharp}
\rangle
-
\int_{\partial \Omega}
\hspace{-0.1cm}
\langle \displaystyle (L_{\rm c}^2\widetilde{a}\,\text{Curl}_{\text{2D}}\,\widetilde{P}^{\sharp}) \times n , \delta \widetilde{P}^{\sharp} \rangle 
\,
\text{d}s
\, ,
\notag
\end{align}
and that, because of (\ref{eq:first_variation_2D}), and highlighting that $\widetilde{P}$ is orthogonal with respect the out-of-plane micro distortion tensor $P$ (their scalar product is zero), implies
\begin{align}
\widetilde{\sigma}
-
\widetilde{\sigma}_{\rm micro}
-
L_{\rm c}^2\widetilde{a}\,\text{Curl} \, \text{Curl}_{\text{2D}} \, \widetilde{P}^{\sharp} = \widetilde{M}
\quad
\text{in }
\,
\Omega \, ,
\qquad\qquad
\displaystyle (L_{\rm c}^2\widetilde{a}\,\text{Curl}_{\text{2D}}\,\widetilde{P}^{\sharp}) \times n = 0
\quad
\text{on }
\,
\partial \Omega \, .
\label{eq:equa_var_P}
\end{align}
where the out-of-plane components of (\ref{eq:equa_var_P})$_1$ and (\ref{eq:equa_var_P})$_2$ must not be considered.
We can now collect all the homogeneous equilibrium equations obtained and the homogeneous Neumann boundary conditions
\begin{align}
\left.
\begin{array}{rrr}
\text{Div}\, \widetilde{\sigma} = \widetilde{f}
\\*
\widetilde{\sigma} - \widetilde{\sigma}_{\rm micro} - L_{\rm c}^2\widetilde{a}\,\text{Curl} \, \text{Curl}_{\text{2D}} \, \widetilde{P} = \widetilde{M}
\end{array}
\right\}
\quad
\text{in }
\,
\Omega
\, ,
\qquad\qquad
\left.
\begin{array}{rrr}
\widetilde{\sigma} \, n = 0
\\*
(L_{\rm c}^2\widetilde{a}\,\text{Curl}_{\text{2D}}\,\widetilde{P}) \times n = 0
\end{array}
\right\}
\quad
\text{on }
\,
\partial \Omega
\, .
\label{eq:equi_RM_2D}
\end{align}
Since $\text{Div} \, (L_{\rm c}^2\widetilde{a}\,\text{Curl} \, \text{Curl}_{\text{2D}}\,\widetilde{P}) =0$, combining the two equation in (\ref{eq:equi_RM_2D})$_1$ gives rise to another equilibrium equation that depends only on $\text{sym}\,P$
\begin{align}
\left.
\begin{array}{rrr}
\text{Div}\, \widetilde{\sigma} = \widetilde{f}
\\*
\widetilde{\sigma} - \widetilde{\sigma}_{\rm micro} - L_{\rm c}^2\widetilde{a}\,\text{Curl} \, \text{Curl}_{\text{2D}} \, \widetilde{P} = \widetilde{M}
\\*
\{
\text{Div}\, \widetilde{\sigma}_{\rm micro} = \widetilde{f} + \text{Div}\, \widetilde{M}
\}
\end{array}
\right\}
\quad
\text{in }
\,
\Omega
\, ,
\qquad\qquad
\left.
\begin{array}{rrr}
\widetilde{\sigma} \, n = 0
\\*
(L_{\rm c}^2\widetilde{a}\,\text{Curl}_{\text{2D}}\,\widetilde{P}) \times n = 0
\end{array}
\right\}
\quad
\text{on }
\,
\partial \Omega
\, .
\label{eq:equi_RM_2D_2}
\end{align}
The extra equation $\text{Div}\, \widetilde{\sigma}_{\rm micro} = \widetilde{f} + \text{Div}\, \widetilde{M}$ is not independent with respect the other two, and any smooth solution of (\ref{eq:equi_RM_2D})$_1$ automatically satisfy it.
This equation can nevertheless substitute $\text{Div}\, \widetilde{\sigma} = \widetilde{f}$, but, although it depends solely on $\text{sym} P$, it is an undetermined system of equations since we just have two equations for three unknown functions ($P_{11},P_{22},P_{12}$). 
The equilibrium equations~(\ref{eq:equi_RM_2D})$_1$ in components are
\begin{align}
\label{pstrain1}
(\lambda_{\rm e}+2 \mu_{\rm e})
\left(u_{1,11}-P_{11,1}\right)
+\lambda_{\rm e} \left(u_{2,12}-P_{22,1}\right)
\quad
&
\notag
\\*
+\mu_{\rm c} \left(u_{1,22}-u_{2,12}-P_{12,2}+P_{21,2}\right)
+\mu_{\rm e} \left(u_{1,22}+u_{2,12}-P_{12,2}-P_{21,2}\right)
&=f_1 \, ,
\\*[5pt]
\label{pstrain2}
(\lambda_{\rm e}+2 \mu_{\rm e}) \left(u_{2,22}-P_{22,2}\right)
+\lambda_{\rm e} \left(u_{1,12}-P_{11,2}\right)
\quad
&
\notag
\\*
+\mu_{\rm c} \left(u_{2,11}-u_{1,12}-P_{21,1}+P_{12,1}\right)
+\mu_{\rm e} \left(u_{2,11}+u_{1,12}-P_{21,1}-P_{12,1}\right)
&=f_2 \, ,
\\*[10pt]
\label{pstrain3}
L_{\rm c}^2\, \widetilde{a} \left(P_{11,22}-P_{12,12}\right)
\quad
&
\notag
\\*
-P_{11} (\lambda_{\rm e}+\lambda_{\rm m}+2 (\mu_{\rm e}+\mu_{\rm m}))
-(\lambda_{\rm e}+\lambda_{\rm m}) P_{22}
+(\lambda_{\rm e}+2 \mu_{\rm e}) u_{1,1}+
\lambda_{\rm e} u_{2,2}
&=M_{11} \, ,
\\*[5pt]
\label{pstrain4}
-L_{\rm c}^2\, \widetilde{a} \left(P_{11,12}-P_{12,11}\right)
\quad
&
\notag
\\*
-(\mu_{\rm c}+\mu_{\rm e}+\mu_{\rm m}) P_{12}
+(\mu_{\rm c}-\mu_{\rm e}-\mu_{\rm m}) P_{21}
+(\mu_{\rm c}+\mu_{\rm e}) u_{1,2}
+(\mu_{\rm e}-\mu_{\rm c}) u_{2,1}
&=M_{12} \, ,
\\*[5pt]
\label{pstrain5}
L_{\rm c}^2\, \widetilde{a} \left(P_{21,22}-P_{22,12}\right)
\quad
&
\notag
\\*
+(\mu_{\rm c}-\mu_{\rm e}-\mu_{\rm m}) P_{12}
-(\mu_{\rm c}+\mu_{\rm e}+\mu_{\rm m}) P_{21}
+(\mu_{\rm e}-\mu_{\rm c}) u_{1,2}
+(\mu_{\rm c}+\mu_{\rm e}) u_{2,1}
&=M_{21} \, ,
\\*[5pt]
\label{pstrain6}
-L_{\rm c}^2\, \widetilde{a} \left(P_{21,12}-P_{22,11}\right)
\quad
&
\notag
\\*
-P_{22} (\lambda_{\rm e}+\lambda_{\rm m}+2 (\mu_{\rm e}+\mu_{\rm m}))
-(\lambda_{\rm e}+\lambda_{\rm m}) P_{11}
+(\lambda_{\rm e}+2 \mu_{\rm e}) u_{2,2}
+\lambda_{\rm e} u_{1,1}
&=M_{22} \, .
\end{align}
which, according to \eqref{constEq}, are accompanied by the following constitutive plane strain equations
\begin{align}
\label{const1}
\sigma_{11}=(\lambda_{\rm e}+2\mu_{\rm e})u_{1,1}+\lambda_{\rm e} u_{2,2}-(\lambda_{\rm e}+2\mu_{\rm e})P_{11}-\lambda{\rm e} P_{22}\, ,
\\[5pt]
\label{const2}
\sigma_{22}=(\lambda_{\rm e}+2\mu_{\rm e})u_{2,2}+\lambda_{\rm e} u_{1,1}-(\lambda_{\rm e}+2\mu_{\rm e})P_{22}-\lambda{\rm e} P_{11}\, ,
\\*[5pt]
\label{const3}
\sigma_{12}=(\mu_{\rm e}+\mu_{\rm c})u_{1,2}+(\mu_{\rm e}-\mu_{\rm c})u_{2,1}-(\mu_{\rm e}+\mu_{\rm c})P_{12}-(\mu_{\rm e}-\mu_{\rm c})P_{21}\, ,
\\*[5pt]
\label{const4}
\sigma_{21}=(\mu_{\rm e}+\mu_{\rm c})u_{2,1}+(\mu_{\rm e}-\mu_{\rm c})u_{1,2}-(\mu_{\rm e}+\mu_{\rm c})P_{21}-(\mu_{\rm e}-\mu_{\rm c})P_{12}\, ,
\\*[5pt]
\label{const5}
m_{13}=-\frac{1}{2} L_{\rm c}^2 (a_1+a_2)(P_{11,2}-P_{12,1})\, ,
\\*[5pt]
\label{const6}
m_{23}=-\frac{1}{2} L_{\rm c}^2 (a_1+a_2)(P_{21,2}-P_{22,1}) \, .
\end{align}
\subsection{A true two-scale model: the relaxed micromorphic model limit for $L_{\rm c}\to 0$ and $L_{\rm c}\to \infty$ in plane strain}
The relaxed micromorphic model reduces to a classical Cauchy model for both $L_{\rm c}\to 0$ and $L_{\rm c}\to \infty$ but with two different stiffnesses, $\mathbb{C}_{\rm macro}$ and $\mathbb{C}_{\rm micro}$, respectively.
The expressions of such stiffnesses are presented in the next two sections for the plane strain problem.
\subsubsection{Limit for $L_{\rm c}\to 0$: macroscopic stiffness $\mathbb{C}_{\rm macro}$}
For the limit $L_{\rm c} \to 0$, the equilibrium equations (\ref{eq:equi_RM_2D}) reduce to
\begin{align}
\text{Div}
\big[
2\mu_{\rm e}\,\text{sym}  (\text{D}\widetilde{u}^{\sharp} - \widetilde{P}^{\sharp} )
+ 2\mu_{\rm c}\,\text{skew} (\text{D}\widetilde{u}^{\sharp} - \widetilde{P}^{\sharp} )
+ \lambda_{\rm e} \text{tr} (\text{D}\widetilde{u}^{\sharp} - \widetilde{P}^{\sharp} ) \mathbbm{1}_2
\big]
&=
\widetilde{f}
\, ,
\label{eq:equi_RM_2D_Lc_0}
\\*
2\mu_{\rm e}\,\text{sym}  (\text{D}\widetilde{u}^{\sharp} - \widetilde{P}^{\sharp} )
+ 2\mu_{\rm c}\,\text{skew} (\text{D}\widetilde{u}^{\sharp} - \widetilde{P}^{\sharp} )
+ \lambda_{\rm e} \text{tr} (\text{D}\widetilde{u}^{\sharp} - \widetilde{P}^{\sharp} ) \mathbbm{1}_2
- 2 \mu_{\rm micro}\,\text{sym}\,\widetilde{P}^{\sharp}
- \lambda_{\rm micro} \text{tr} ( \widetilde{P}^{\sharp} ) \mathbbm{1}_2
&=
\widetilde{M}
\, .
\notag
\end{align}
The equation (\ref{eq:equi_RM_2D_Lc_0})$_2$ is now algebraic in $\widetilde{P}^{\sharp}$.
Thanks to the orthogonality of the ``sym/skew'' decomposition, the equation (\ref{eq:equi_RM_2D_Lc_0})$_2$ requires that
\begin{align}
2\mu_{\rm c} \, \text{skew} (\text{D}\widetilde{u}^{\sharp} - \widetilde{P}^{\sharp} ) 
&=
\text{sym} \, \widetilde{M}
\, ,
\\
2\mu_{\rm e}\,\text{sym}  (\text{D}\widetilde{u}^{\sharp} - \widetilde{P}^{\sharp} )
+ \lambda_{\rm e} \text{tr} (\text{D}\widetilde{u}^{\sharp} - \widetilde{P}^{\sharp} ) \mathbbm{1}_2
- 2 \mu_{\rm micro}\,\text{sym}\,\widetilde{P}^{\sharp}
- \lambda_{\rm micro} \text{tr} ( \widetilde{P}^{\sharp} ) \mathbbm{1}_2
&=
\text{skew} \, \widetilde{M}
\, .
\notag
\end{align}
Since the ``sym'' operator is not ortogonal to the ``tr'' operator, we further decompose ``sym'' into ``dev sym'' and ``tr sym'' so that
\begin{align}
2\mu_{\rm c} \, \text{skew} (\text{D}\widetilde{u}^{\sharp} - \widetilde{P}^{\sharp} ) 
&=
\text{skew} \, \widetilde{M}
\, ,
\quad
\notag
\\
2\mu_{\rm e}\,\text{dev sym}  (\text{D}\widetilde{u}^{\sharp} - \widetilde{P}^{\sharp} )
+ \mu_{\rm e}\,\text{tr}  (\text{D}\widetilde{u}^{\sharp} - \widetilde{P}^{\sharp} ) \mathbbm{1}_2
+ \lambda_{\rm e} \text{tr} (\text{D}\widetilde{u}^{\sharp} - \widetilde{P}^{\sharp} ) \mathbbm{1}_2
\quad
&
\label{eq:equi_RM_2D_Lc_0_dev_sym}
\\*
- 2 \mu_{\rm micro}\,\text{dev sym}\,\widetilde{P}^{\sharp}
- \mu_{\rm micro}\,\text{tr}\,(\widetilde{P}^{\sharp}) \mathbbm{1}_2
- \lambda_{\rm micro} \text{tr} ( \widetilde{P}^{\sharp} ) \mathbbm{1}_2
&=
\text{sym} \, \widetilde{M}
\, .
\notag
\end{align}
note that ``tr sym'' is the same as ``tr''.
We also recall here the definition of the volumetric part, the deviatoric part, and the skew-symmetric parts in plane strain case
\begin{align}
\text{2D volumetric part}
&\coloneqq
\frac{1}{2} \text{tr}(\widetilde{P}^{\sharp})\mathbbm{1}_2 \, , \quad \mathbbm{1}_2=\begin{pmatrix} 1 & 0 \\ 0 & 1 \end{pmatrix}
\\*
\text{2D deviatoric symmetric part}
&\coloneqq
\frac{\widetilde{P}^{\sharp}+\widetilde{P}^{\sharp^{\rm T}}}{2}-\frac{1}{2} \text{tr}(\widetilde{P}^{\sharp})\mathbbm{1}_2 \, ,
\\*
\label{skewP}
\text{2D skew symmetric part}
&\coloneqq
\frac{\widetilde{P}^{\sharp}-\widetilde{P}^{\sharp^{\rm T}}}{2} \, .
\end{align}
With further manipulations and thanks to the orthogonality of the operator ``skew'', ``dev sym'', and ``tr'', the system (\ref{eq:equi_RM_2D_Lc_0_dev_sym}) requires that
\begin{align}
2\mu_{\rm c} \, \text{skew} (\text{D}\widetilde{u}^{\sharp} - \widetilde{P}^{\sharp} ) 
&=
\text{skew} \, \widetilde{M}
\, ,
\notag
\\
\mu_{\rm e}\,\text{dev sym}  (\text{D}\widetilde{u}^{\sharp} - \widetilde{P}^{\sharp} )
- \mu_{\rm micro}\,\text{dev sym}\,\widetilde{P}^{\sharp}
&=
\text{dev sym} \, \widetilde{M}
\label{eq:equi_RM_2D_Lc_0_dev_sym_2}
\, ,
\\*
\left(\mu_{\rm e}+ \lambda_{\rm e}\right)\,\text{tr}  (\text{D}\widetilde{u}^{\sharp} - \widetilde{P}^{\sharp} ) \mathbbm{1}_2
- \left(\mu_{\rm micro} + \lambda_{\rm micro}\right)\,\text{tr}\,(\widetilde{P}^{\sharp}) \mathbbm{1}_2
&=
\frac{1}{2}\text{tr} (\widetilde{M}) \mathbbm{1}_2
\, .
\notag
\end{align}
From equation (\ref{eq:equi_RM_2D_Lc_0_dev_sym_2}) we can evaluate the expressions for $\text{skew} \, \widetilde{P}^{\sharp}$, $\text{dev sym} \, \widetilde{P}^{\sharp}$, and $\text{tr}(\widetilde{P}^{\sharp})$ as
\begin{align}
\text{skew} \, \text{D}\widetilde{u}^{\sharp}
-
\frac{1}{2 \mu_{\rm c}}\text{skew}\,\widetilde{M}
&=
\text{skew} \, \widetilde{P}^{\sharp}
\, ,
\notag
\\*
\frac{\mu_{\rm e}}{\mu_{\rm e}+\mu_{\rm micro}}
\text{dev sym}  \, \text{D}\widetilde{u}^{\sharp}
-
\frac{1}{2 (\mu_{\rm e}+\mu_{\rm micro})}\text{sym}\,\widetilde{M}
&=
\text{dev sym}  \, \widetilde{P}^{\sharp}
\, ,
\label{eq:equi_RM_2D_Lc_0_dev_sym_3}
\\*
\frac{\widetilde{\kappa}_{\rm e}}{\widetilde{\kappa}_{\rm e}+\widetilde{\kappa}_{\rm micro}}
\text{tr}  \, \text{D}\widetilde{u}^{\sharp}
-
\frac{1}{2 (\widetilde{\kappa}_{\rm e}+\widetilde{\kappa}_{\rm micro})}\text{tr}\,\widetilde{M}
&=
\text{tr}(\widetilde{P}^{\sharp})
\, .
\notag
\end{align}
where  $\widetilde{\kappa}_{\rm e}=\mu_{\rm e}+\lambda_{\rm e}$ and $\widetilde{\kappa}_{\rm micro}=\mu_{\rm micro}+\lambda_{\rm micro}$ are the plane strain bulk moduli.

Substituting back the relations (\ref{eq:equi_RM_2D_Lc_0_dev_sym_3}) in the equation~(\ref{eq:equi_RM_2D_Lc_0})$_1$ while also applying the ``dev sym'', and ``tr'' decomposition, we have
\begin{align}
\text{Div}
\big[
2\mu_{\rm e} \, \text{dev sym} 
\left(\text{D}\widetilde{u}^{\sharp} - \left(\frac{\mu_{\rm e}}{\mu_{\rm e}+\mu_{\rm micro}}\text{D}\widetilde{u}^{\sharp}\right) \right)
+ \widetilde{\kappa}_{\rm e} \, \text{tr}
\left(\text{D}\widetilde{u}^{\sharp} - \left(\frac{\widetilde{\kappa}_{\rm e}}{\widetilde{\kappa}_{\rm e}+\widetilde{\kappa}_{\rm micro}} \, \text{D}\widetilde{u}^{\sharp}\right) \right) \mathbbm{1}_2
\big]
&=
\widetilde{f}^{*}
\, ,
\notag
\\*[5pt]
\Longleftrightarrow
\qquad
\text{Div}
\big[
2\dfrac{\mu_{\rm e} \, \mu_{\rm micro}}{\mu_{\rm e} + \mu_{\rm micro}} \, \text{dev sym} \,\text{D}\widetilde{u}^{\sharp}
+ \dfrac{\widetilde{\kappa}_{\rm e} \, \widetilde{\kappa}_{\rm micro}}{\widetilde{\kappa}_{\rm e} + \widetilde{\kappa}_{\rm micro}} \, \text{tr}
\left(\text{D}\widetilde{u}^{\sharp} \right) \mathbbm{1}_2
\big]
&=
\widetilde{f}^{*}
\, ,
\label{eq:equi_RM_2D_Lc_0_2}
\\*[5pt]
\Longleftrightarrow
\qquad
\text{Div}
\big[
2\mu_{\rm macro} \, \text{dev sym}\,\text{D}\widetilde{u}^{\sharp}
+ \widetilde{\kappa}_{\rm macro} \, \text{tr}
(\text{D}\widetilde{u}^{\sharp}) \mathbbm{1}_2
\big]
&=
\widetilde{f}^{*}
\, .
\notag
\end{align}
where $\widetilde{f}^*$ is defined as
\begin{align}
\widetilde{f}^{*} \coloneqq
\widetilde{f} -
\text{Div}\big[
\frac{\mu_{\rm macro}}{\mu_{\rm micro}} \, \text{dev sym} \, \widetilde{M}
+
\text{skew} \, \widetilde{M}
+
\frac{1}{2}\frac{\widetilde{\kappa}_{\rm macro}}{\widetilde{\kappa}_{\rm micro}} \, \text{tr} (\widetilde{M}) \mathbbm{1}_2
\big]
\, .
\end{align}
It is noted that $\widetilde{f}^{*}$ depends on $\text{skew} \, \widetilde{M}$ without any multiplicative elastic coefficient because of the choice of an isotropic constitutive law (an isotropic second order skew-symmetric tensor depends on one coefficient). This limit with a concentrated double body force may be instrumental in order to identify the \textit{micro} parameters. The equation (\ref{eq:equi_RM_2D_Lc_0_2})$_3$ is the equilibrium equation for a classical linear elastic isotropic Cauchy continuum with stiffness $\mu_{\rm macro}$ and $\kappa_{\rm macro}$.

The relations for the macro Lamé parameters ($\mu_{\rm macro},\lambda_{\rm macro}$) and the macroscopic bulk modulus for plane strain ($\widetilde{\kappa}_{\rm macro}$) are then
\begin{align}
\mu_{\rm macro} &\coloneqq \dfrac{\mu_{\rm e} \, \mu_{\rm micro}}{\mu_{\rm e} + \mu_{\rm micro}} \, ,
\qquad\qquad\qquad
\widetilde{\kappa}_{\rm macro} \coloneqq \dfrac{\widetilde{\kappa}_{\rm e} \, \widetilde{\kappa}_{\rm micro}}{\widetilde{\kappa}_{\rm e} + \widetilde{\kappa}_{\rm micro}} \, ,
\\*[5pt]
\widetilde{\lambda}_{\rm macro} &\coloneqq
\dfrac{(\mu_{\rm e} + \lambda_{\rm e})(\mu_{\rm micro} + \lambda_{\rm micro})}{(\mu_{\rm e} + \lambda_{\rm e}) + (\mu_{\rm micro} + \lambda_{\rm micro})} -\dfrac{\mu_{\rm e} \, \mu_{\rm micro}}{\mu_{\rm e} + \mu_{\rm micro}}
\, ,
\notag
\end{align}
where  $\widetilde{\kappa}_{\rm macro}= \mu_{\rm macro}+\widetilde{\lambda}_{\rm macro} $.
In order to have $\widetilde{\lambda}_{\rm macro}=\lambda_{\rm micro}=0$, the only possible condition is again $\lambda_{\rm micro}=\lambda_{\rm e}=0$.
\subsubsection{Limit for $L_{\rm c}\to \infty$: microscopic stiffness $\mathbb{C}_{\rm micro}$}
The minimization of an energy functional that incorporate $L_{\rm c}^2 \, \| \text{Curl} \widetilde{P}^{\sharp} \|^2$, for the limit $L_{\rm c} \to \infty$, requires $\text{Curl} \widetilde{P}^{\sharp}=0$, and this implies that the micro-distortion tensor $P$ has to reduce to a gradient field $\widetilde{P}^{\sharp} \to \text{D}\widetilde{v}^{\sharp}$ on a simply connected domain and
\begin{align}
\text{Curl} \, \text{D}\widetilde{v}^{\sharp} = 0 \qquad \forall \, \widetilde{v}^{\sharp} \in [\mathit{C}^\infty(\Omega)]^3 \, ,
\end{align}
thus asserting finite energies of the relaxed micromorphic model for arbitrarily large characteristic length values $L_{\rm c}$.
The corresponding strain energy density in terms of the reduced kinematics $\{ \widetilde{u}, \widetilde{v}^{\sharp} \} : \Omega\to \mathbb{R}^{3}$ now reads
\begin{align}
W \left(\text{D}\widetilde{u}, \text{D}\widetilde{v}^{\sharp}\right) = &
\, \mu_{\rm e} \left\lVert \text{sym} (\text{D}\widetilde{u}^{\sharp} - \text{D}\widetilde{v}^{\sharp} ) \right\rVert^{2}
+ \mu_{\rm c} \left\lVert \text{skew} (\text{D}\widetilde{u}^{\sharp} - \text{D}\widetilde{v}^{\sharp} ) \right\rVert^{2}
+ \frac{\lambda_{\rm e}}{2} \text{tr}^2 (\text{D}\widetilde{u}^{\sharp} - \text{D}\widetilde{v}^{\sharp} )
\label{eq:energy_RM_2D_Lc_infty}
\\*
&
+ \mu_{\rm micro} \left\lVert \text{sym}\,\text{D}\widetilde{v}^{\sharp} \right\rVert^{2}
+ \frac{\lambda_{\rm micro}}{2} \text{tr}^2 \left(\text{D}\widetilde{v}^{\sharp} \right)
\, .
\notag
\end{align}
The first variation of the strain energy $I=\displaystyle\int_{\Omega} W \, \text{d}x$ with respect to the two vector fields $\widetilde{u}^{\sharp}$ and $\widetilde{v}^{\sharp}$ leads to
\begin{align}
\delta I^{\widetilde{u}}
&=
\int_{\Omega}
\hspace{-0.1cm}
\Big(
2\mu_{\rm e} \,
\langle \text{sym} (\text{D}\widetilde{u}^{\sharp} - \text{D}\widetilde{v}^{\sharp} ),\text{D}\delta \widetilde{u}^{\sharp} \rangle
+ 2\mu_{\rm c} \,
\langle \text{skew} (\text{D}\widetilde{u}^{\sharp} - \text{D}\widetilde{v}^{\sharp} ),\text{D}\delta \widetilde{u}^{\sharp} \rangle
\label{eq:first_varia_energy_RM_2D_u_Lc_infy}
\\*[5pt]
&
\phantom{=\int_{\Omega}\Big(}
+ \lambda_{\rm e} \langle \text{tr} (\text{D}\widetilde{u}^{\sharp} - \text{D}\widetilde{v}^{\sharp} ) \mathbbm{1}_2, \text{D}\delta \widetilde{u}^{\sharp} \rangle
\Big)
\text{d}x
\, ,
\notag
\\*[5pt]
\delta I^{\widetilde{v}^{\sharp}}
&=
\int_{\Omega}
\hspace{-0.1cm}
\Big(
-2\mu_{\rm e} \,
\langle \text{sym} (\text{D}\widetilde{u}^{\sharp} - \text{D}\widetilde{v}^{\sharp} ),\text{D}\delta \widetilde{v}^{\sharp} \rangle
-2\mu_{\rm c} \,
\langle \text{skew} (\text{D}\widetilde{u}^{\sharp} - \text{D}\widetilde{v}^{\sharp} ),\text{D}\delta \widetilde{v}^{\sharp} \rangle
-\lambda_{\rm e} \langle \text{tr} (\text{D}\widetilde{u}^{\sharp} - \text{D}\widetilde{v}^{\sharp} ) \mathbbm{1}_2, \text{D}\delta \widetilde{v}^{\sharp} \rangle
\label{eq:first_varia_energy_RM_2D_v_Lc_infy}
\\*
&
\phantom{=\int_{\Omega}\Big(}
+
2\mu_{\rm micro} \langle \text{sym}\,\text{D}\widetilde{v}^{\sharp},\text{D}\delta \widetilde{v}^{\sharp} \rangle
+ \lambda_{\rm micro} \langle \text{tr} (\text{D}\widetilde{v}^{\sharp} ) \mathbbm{1}_2, \text{D}\delta \widetilde{v}^{\sharp} \rangle
\Big)
\text{d}x
\, .
\notag
\end{align}
The equilibrium equations are now obtained by requiring 
\begin{equation}
\delta I^{\widetilde{u}^{\sharp}}= \langle \widetilde{f} , \delta \widetilde{u}^{\sharp} \rangle \, ,
\quad
\forall
\,
\delta \widetilde{u}^{\sharp}
\qquad\qquad
\text{and}
\qquad\qquad
\delta I^{\widetilde{v}^{\sharp}}=
\langle \widetilde{M} , \text{D}\delta \widetilde{v}^{\sharp} \rangle \, ,
\quad
\forall
\,
\delta \widetilde{v}^{\sharp}  \, .
\label{eq:first_variation_2D_Lc_infy}
\end{equation}
where the contributions on the right sides are the virtual work of the external forces $\widetilde{f}$ (classical body force) and $\widetilde{M}$ (non-symmetric second order double body force tensor), and the equilibrium equations read
\begin{align}
\text{Div}
\big[
2\mu_{\rm e}\,\text{sym}  (\text{D}\widetilde{u}^{\sharp} - \text{D}\widetilde{v}^{\sharp} )
+ 2\mu_{\rm c}\,\text{skew} (\text{D}\widetilde{u}^{\sharp} - \text{D}\widetilde{v}^{\sharp} )
+ \lambda_{\rm e} \, \text{tr} (\text{D}\widetilde{u}^{\sharp} - \text{D}\widetilde{v}^{\sharp} ) \mathbbm{1}_2
\big]
&=
\widetilde{f}
\, ,
\label{eq:equi_RM_2D_Lc_infy_2}
\\*
-\text{Div}
\big[
2\mu_{\rm e}\,\text{sym}  (\text{D}\widetilde{u}^{\sharp} - \text{D}\widetilde{v}^{\sharp} )
+ 2\mu_{\rm c}\,\text{skew} (\text{D}\widetilde{u}^{\sharp} - \text{D}\widetilde{v}^{\sharp} )
+ \lambda_{\rm e} \,\text{tr} (\text{D}\widetilde{u}^{\sharp} - \text{D}\widetilde{v}^{\sharp} ) \mathbbm{1}_2
\big]
\quad
\notag
\\*
+
\text{Div}
\big[
2\mu_{\rm micro}\,\text{sym}  \text{D}\widetilde{v}^{\sharp}
+ \lambda_{\rm micro} \,\text{tr} ( \text{D}\widetilde{v}^{\sharp} ) \mathbbm{1}_2
\big]
&=
\text{Div} \, \widetilde{M}
\, ,
\notag
\end{align}
where the constraint $\widetilde{M} \, n =0$ is required on the boundary, with $n$ the normal to the boundary.
The term on the left-hand side of equation (\ref{eq:equi_RM_2D_Lc_infy_2})$_2$ can be substituted with the right-hand side of (\ref{eq:equi_RM_2D_Lc_infy_2})$_1$ and, while keeping the equation (\ref{eq:equi_RM_2D_Lc_infy_2})$_1$, we can re-write the system of equations (\ref{eq:equi_RM_2D_Lc_infy_2}) as
\begin{align}
\text{Div}
\big[
2\mu_{\rm e}\,\text{sym}  (\text{D}\widetilde{u}^{\sharp} - \text{D}\widetilde{v}^{\sharp} )
+ 2\mu_{\rm c}\,\text{skew} (\text{D}\widetilde{u}^{\sharp} - \text{D}\widetilde{v}^{\sharp} )
+ \lambda_{\rm e} \,\text{tr} (\text{D}\widetilde{u}^{\sharp} - \text{D}\widetilde{v}^{\sharp} ) \mathbbm{1}_2
\big]
&=
\widetilde{f}
\, ,
\label{eq:equi_RM_2D_Lc_infy_3}
\\*
\text{Div}
\big[
2\mu_{\rm micro}\,\text{sym} \, \text{D}\widetilde{v}^{\sharp}
+ \lambda_{\rm micro} \text{tr} \, ( \text{D}\widetilde{v}^{\sharp} ) \mathbbm{1}_2
\big]
&=
\widetilde{f} + \text{Div} \, \widetilde{M}
\, .
\notag
\end{align}
The only case in which $\widetilde{v}^{\sharp}=\widetilde{u}^{\sharp}$ is an admissible solution is if the classical body forces $\widetilde{f}$ are zero.
In this case  (\ref{eq:equi_RM_2D_Lc_infy_3}) reduces to
\begin{align}
\text{Div} \, \sigma_{\rm micro}
=
\text{Div}
\big[
2\mu_{\rm micro}\,\text{sym} \, \text{D}\widetilde{u}^{\sharp}
+ \lambda_{\rm micro} \, \text{tr} ( \text{D}\widetilde{u}^{\sharp} ) \mathbbm{1}_2
\big]
=
\text{Div} \, \widetilde{M}
\, ,
\label{eq:equi_RM_2D_Lc_infy_4}
\end{align}
which is an equilibrium equation of the classical elasticity type with a micro stiffness given by $\mu_{\rm micro}$ and $\lambda_{\rm micro}$ and a body force vector equal to $\text{Div} \, \widetilde{M}$.
\subsection{A particular case of the relaxed micromorphic model with a zero Cosserat couple modulus}
\label{sec:zero_coss}
If we additionally set $\mu_{\rm c}=0$, the force stress tensor $\sigma$ becomes symmetric and the model further reduces to
\begin{align}
\text{Div}
\big[
\overbrace{
2\mu_{\rm e}\,\text{sym}  (\text{D}\widetilde{u}^{\sharp} - \widetilde{P}^{\sharp} )
+ \lambda_{\rm e} \text{tr} (\text{D}\widetilde{u}^{\sharp} - \widetilde{P}^{\sharp} ) \mathbbm{1}_2
}^{
\mathlarger{\sigma}\coloneqq
}
\big]
&=
\widetilde{f}
\, ,
\notag
\\*
\sigma
- 2 \mu_{\rm micro}\,\text{sym}\,\widetilde{P}^{\sharp}
- \lambda_{\rm micro} \text{tr} ( \widetilde{P}^{\sharp} ) \mathbbm{1}_2
- L_{\rm c}^2 \, \widetilde{a}\,
\text{Curl} \, \text{Curl}_{\text{2D}} \, \widetilde{P}^{\sharp}
&=
\widetilde{M}
\, ,
\label{eq:equi_RM_2D_CC_mc0}
\end{align}
\begin{align}
\widetilde{M}
&=
\left(
\begin{array}{ccccccc}
M_{11} & M_{12} & 0\\
M_{21} & M_{22} & 0\\
0                  & 0                  & 0
\end{array}
\right)
\, ,
\qquad\qquad\qquad
\widetilde{f}=
\left(
\begin{array}{ccccccc}
f_{1} \\
f_{2} \\
0                  
\end{array}
\right)
\, .
\end{align}
\subsection{A simpler case of the relaxed micromorphic model with zero micro and macro Poisson's ratio}
\label{sec:zero_poisson}
If we set $\lambda_{\rm micro}=\lambda_{\rm e}=0$, which implies $\lambda_{\rm macro}=0$, the equilibrium equations (\ref{eq:equi_RM_2D}) further reduce to 
\begin{align}
\text{Div}
\big[
2\mu_{\rm e}\,\text{sym}  (\text{D}\widetilde{u}^{\sharp} - \widetilde{P}^{\sharp} )
+ 2\mu_{\rm c}\,\text{skew} (\text{D}\widetilde{u}^{\sharp} - \widetilde{P}^{\sharp} )
\big]
&=
\widetilde{f}
\, ,
\label{eq:equi_RM_2D_CC_l0}
\\*
2\mu_{\rm e}\,\text{sym}  (\text{D}\widetilde{u}^{\sharp} - \widetilde{P}^{\sharp} )
+ 2\mu_{\rm c}\,\text{skew} (\text{D}\widetilde{u}^{\sharp} - \widetilde{P}^{\sharp} )
- 2 \mu_{\rm micro}\,\text{sym}\,\widetilde{P}^{\sharp}
- L_{\rm c}^2 \, \widetilde{a} \, \text{Curl} \, \text{Curl}_{\text{2D}} \, \widetilde{P}^{\sharp}
&=
\widetilde{M}
\, .
\notag
\end{align}
Componentwise we have
\begin{align}\begin{split}
\mu_{\rm c} \left(u_{1,22}-u_{2,12}+P_{21,2}-P_{12,2}\right)
+\mu_{\rm e} \left(u_{1,22}+2 u_{1,11}+u_{2,12}-2 P_{11,1}-P_{12,2}-P_{21,2}\right)
&=f_1 \, ,
\\*[5pt]
\mu_{\rm c} \left(P_{12,1}-P_{21,1}-u_{1,12}+u_{2,11}\right)
+\mu_{\rm e} \left(u_{1,12}+2 u_{2,22}+u_{2,11}-P_{12,1}-P_{21,1}-2 P_{22,2}\right)
&=f_2 \, ,
\\*[10pt]
\widetilde{a} L_{\rm c}^2 \left(P_{11,22}-P_{12,12}\right)
+2\mu_{\rm e} \left(u_{1,1}-P_{11}\right)
-2 \mu_{\rm m} P_{11}
&=M_{11} \, ,
\\*[5pt]
\widetilde{a} L_{\rm c}^2 \left(P_{12,11}-P_{11,12}\right)
\quad
&
\\*
+\mu_{\rm c} \left(u_{1,2}-u_{2,1}-P_{12}+P_{21}\right)
+\mu_{\rm e} \left(u_{1,2}+u_{2,1}-P_{12}-P_{21}\right)
-\mu_{\rm m} (P_{12}+P_{21})
&=M_{12} \, ,
\\*[5pt]
\widetilde{a} L_{\rm c}^2 \left(P_{21,22}-P_{22,12}\right)
\quad
&
\\*
+\mu_{\rm c} \left(u_{2,1}-u_{1,2}+P_{12}-P_{21}\right)
+\mu_{\rm e} \left(u_{1,2}+u_{2,1}-P_{12}-P_{21}\right)
-\mu_{\rm m} (P_{12}+P_{21}) 
&=M_{21} \, ,
\\*[5pt]
\widetilde{a} L_{\rm c}^2 \left(P_{22,11}-P_{21,12}\right)
+2\mu_{\rm e} \left(u_{2,2}-P_{22}\right)
-2 \mu_{\rm m} P_{22}
&=M_{22} \, ,\end{split}
\end{align}
where we used the abbreviation $\mu_{\rm m}=\mu_{\rm micro}$.
The conditions for existence and uniqueness for the model in  (\ref{eq:equi_RM_2D_CC_l0}) are
\begin{align}
\mu_{\rm e}>0 \, ,
\qquad\qquad
\mu_{\rm micro}>0 \, ,
\qquad\qquad
L_{\rm c}^2\widetilde{a}>0 \, ,
\qquad\qquad
\mu_{\rm c}\geq0 \, .
\end{align}
For $\mu_{\rm c}\equiv0$, in order to guarantee existence and uniqueness, one needs tangential boundary conditions for $\widetilde{P}$, while for $\mu_{\rm c}>0$, one does not need boundary conditions for $\widetilde{P}$ in order to guarantee existence and uniqueness.
\subsection{A simpler case of the relaxed micromorphic model with one curvature parameter, a zero Cosserat couple modulus, and a zero micro and macro Poisson's ratio}
If in addition to the simplifications of Sec.~\ref{sec:zero_poisson} we also set $\mu_{\rm c}=0$, the equilibrium equations (\ref{eq:equi_RM_2D_CC_l0}) further reduce to 
\begin{align}
\text{Div}
\big[
2\mu_{\rm e}\,\text{sym}  (\text{D}\widetilde{u}^{\sharp} - \widetilde{P}^{\sharp} )
\big]
&=
\widetilde{f}
\, ,
\label{eq:equi_RM_2D_CC_mc0_l0}
\\*
2\mu_{\rm e}\,\text{sym}  (\text{D}\widetilde{u}^{\sharp} - \widetilde{P}^{\sharp} )
- 2 \mu_{\rm micro}\,\text{sym}\,\widetilde{P}^{\sharp}
- L_{\rm c}^2 \, \widetilde{a} \, \text{Curl} \, \text{Curl}_{\text{2D}} \, \widetilde{P}^{\sharp}
&=
\widetilde{M}
\, .
\notag
\end{align}

\futurelet \reserved@a \@xhline
}
\newcolumntype{"}{@{\hskip\tabcolsep\vrule width 1pt\hskip\tabcolsep}}
\makeatother
\theoremstyle{definition}

\makeatletter
\newcommand{\pushright}[1]{\ifmeasuring@#1\else\omit\hfill$\displaystyle#1$\fi\ignorespaces}
\newcommand{\pushleft}[1]{\ifmeasuring@#1\else\omit$\displaystyle#1$\hfill\fi\ignorespaces}
\makeatother

\titleformat{\paragraph}
{\normalfont\normalsize\bfseries}{\theparagraph}{1em}{}
\titlespacing*{\paragraph}
{0pt}{3.25ex plus 1ex minus .2ex}{1.5ex plus .2ex}

\newcounter{alphasect}
\def\alphainsection{0}

\let\oldsection=\section
\def\section{%
  \ifnum\alphainsection=1%
    \addtocounter{alphasect}{1}
  \fi%
\oldsection}%
\renewcommand\thesection{%
  \ifnum\alphainsection=1%
    \Alph{alphasect}%
  \else%
    \arabic{section}%
  \fi%
}%

\newenvironment{alphasection}{%
  \ifnum\alphainsection=1%
    \errhelp={Let other blocks end at the beginning of the next block.}
    \errmessage{Nested Alpha section not allowed}
  \fi%
  \setcounter{alphasect}{0}
  \def\alphainsection{1}
}{%
  \setcounter{alphasect}{0}
  \def\alphainsection{0}
}%

\title{Green's functions for the isotropic planar relaxed micromorphic model - concentrated force and concentrated couple}
\author{
Panos Gourgiotis\,\thanks{Corresponding author: Panos Gourgiotis, Mechanics Division, SAMPS, National Technical University of Athens, Zographou, GR-15773, Greece, email: pgourgiotis@mail.ntua.gr},
\quad
Gianluca Rizzi\,\thanks{Gianlucy Rizzi, Faculty of Architecture and Civil Engineering, TU Dortmund, August-Schmidt-Str. 8, 44227 Dortmund, Germany, email: gianluca.rizzi@tu-dortmund.de},
\quad
Peter Lewintan\,\thanks{Peter Lewintan, Faculty of Mathematics, University of Duisburg-Essen, Thea-Leymann-Straße 9, 45127 Essen, Germany, email: peter.lewintan@uni-due.de},
\quad
Davide Bernardini\,\thanks{Davide Bernardini, Department of Structural and Geotechnical Engineering, Sapienza University of Rome, Rome, Italy, email: davide.bernardini@uniroma1.it},
\\[5pt]
Adam Sky\,\thanks{Adam Sky, Institute of Computational Engineering and Sciences, Department of Engineering, Faculty of Science, Technology and Medicine, University of Luxembourg, 6, Avenue de la Fonte, L-4362 Esch-sur-Alzette, Luxembourg, email: adam.sky@uni.lu},\quad
Angela Madeo\,\thanks{Angela Madeo, Faculty of Architecture and Civil Engineering, TU Dortmund, August-Schmidt-Str. 8, 44227 Dortmund, Germany, email: angela.madeo@tu-dortmund.de},
\quad and \quad
Patrizio Neff\,\thanks{Patrizio Neff, Head of Chair for Nonlinear Analysis and Modelling, Fakultät für Mathematik, Universität Duisburg-Essen, Thea-Leymann-Straße 9, 45127 Essen, Germany, email: patrizio.neff@uni-due.de}
}

\thanksmarkseries{arabic}
\date{\today}
\begin{document}
\maketitle
\begin{abstract}
\noindent
We derive the Green's functions (concentrated force and couple in an infinite space) for the isotropic planar relaxed micromorphic model. Since the relaxed micromorphic model particularises into the micro-stretch, Cosserat (micropolar), couple-stress, and linear elasticity model for certain choices of material parameters, we recover the fundamental solutions in all these cases.
\end{abstract}
\textbf{Keywords}: generalized continua, fundamental solutions, Kelvin problem, micromorphic continuum, micro-stretch, Cosserat, micropolar, couple stress, relaxed micromorphic, gauge invariant dislocation model.
{\scriptsize
\tableofcontents
}
%
%
%
%
%
\section{Introduction}
\label{sec:intro}
\begin{figure}[htb!]\center	\unitlength=1mm	\begin{picture}(140,47)	\put(20,2){\def\svgwidth{11 cm}{\small\input{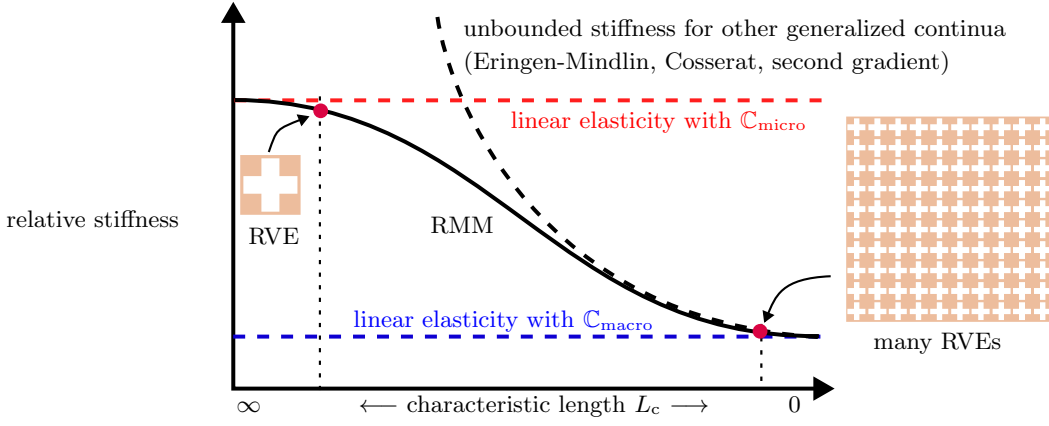}}}	\end{picture}	\caption{The stiffness of the relaxed micromorphic model (RMM) is bounded from above and below. Other generalized continua exhibit unbounded  stiffness for small sizes. For large values of the characteristic length $L_{\rm c}$, linear elasticity with a micro elasticity tensor is recovered (one RVE) while linear elasticity with a macro elasticity tensor is obtained for small values of  the characteristic length (many RVEs).   }	\label{Figure:intro}\end{figure}

The relaxed micromorphic is a new generalised continuum model that allows to describe size-effects and band-gaps appearing in metamaterials \cite{rizzi2021exploring,voss2023modeling,rizzi2022boundary,rizzi2022metamaterial,ramirez2023multi,demore2022unfolding,rizzi2022towards,madeo2015wave} (in its dynamic setting).
The relaxed micromorphic model belongs to the family of micromorphic models \cite{Mindlin1964,Eringen1999} where the kinematics are given by the classical displacement $u:\Omega\to \mathbb{R}^3$ and the non-symmetric micro-distortion $P:\Omega\to \mathbb{R}^{3\times 3}$. The solution is then determined from the variational two-field problem
\begin{align}
I(u,P)=
&\int_{\Omega}
\frac{1}{2}
\bigg(
\langle
\mathbb{C}_{\rm e}\,\text{sym} (\text{D}u - P), \text{sym} (\text{D}u - P)
\rangle
+
\langle
\mathbb{C}_{\rm c}\,\text{skew} (\text{D}u - P), \text{skew} (\text{D}u - P)
\rangle
\label{eq:energy_rmm}
\\*
&\phantom{\int_{\Omega} \frac{1}{2} \bigg(}
+
\langle
\mathbb{C}_{\rm micro}\,\text{sym} \, P, \text{sym} \, P
\rangle
+
\mu_{\rm macro}  L_{\rm c}^2
\langle
\mathbb{L}\,\text{Curl} \, P, \text{Curl} \, P
\rangle
\bigg)
\text{d}x
\quad
\longrightarrow
\quad
\text{min}\,(u,P)
\, .
\notag
\end{align}
Here $\mathbb{C}_{\rm e},\mathbb{C}_{\rm micro},\mathbb{L}$ are positive-definite fourth-order tensors, and $L_{\rm c}$ is a characteristic length and $\mu_{\rm macro}=\mu_{\rm M}$ is the macroscopic shear modulus. Furthermore,
$\mathbb{C}_{\rm c}$ is a positive semi-definite fourth order tensor and we note the homogenization relations \cite{neff2014unifying,neff2020identification}
\begin{gather}
\mathbb{C}_{\rm e}
=
\mathbb{C}_{\rm micro}
\bigg(\mathbb{C}_{\rm micro}-\mathbb{C}_{\rm macro}
\bigg)^{-1}
\mathbb{C}_{\rm macro}
\quad
\Longleftrightarrow
\quad
\mathbb{C}_{\rm macro}
=
\mathbb{C}_{\rm micro}
\bigg(\mathbb{C}_{\rm micro}+\mathbb{C}_{\rm e}
\bigg)^{-1}
\mathbb{C}_{\rm e}
\, ,
\label{eq:Ce_Cmac}
\\
\mathbb{C}_{\rm micro}
=
\mathbb{C}_{\rm e}
\bigg(\mathbb{C}_{\rm e}-\mathbb{C}_{\rm macro}
\bigg)^{-1}
\mathbb{C}_{\rm macro}
\, ,
\notag
\end{gather}
connecting the macroscopic stiffness $\mathbb{C}_{\rm macro}$ uniquely known from classical homogenization for a periodic metamaterial to the stiffness tensors $\mathbb{C}_{\rm micro}$ and $\mathbb{C}_{\rm e}$ of the relaxed micromorphic model.
This new model leverages some of the main shortcomings of the classical Eringen-Mindlin micromorphic model (unbounded stiffness, multitude of parameters).
This is achieved by reducing the complexity of the strain energy function in two ways:
first (i) by excluding some generalities in the local part of the energy,
and second and foremost (ii) by reducing the dependency of the curvature energy acting on a full gradient of the micro-distortion in the classical Mindlin-Eringen model to only a dependency on its Curl.
The consequences of this choice are remarkable: the additional balance equation remains of the second order (Curl is a second order tensor) and the model still includes the better known micro-stretch and Cosserat (micropolar) models (which can be alternatively written in dislocation format with a Curl in the curvature part \cite{ghiba2023cosserat}). Compared to the classical Eringen-Mindlin micromorphic model, note the absence of mixed coupling terms between the elastic strain $\text{sym}\,(\text{D}u- P)$ and the microstrain $\text{sym}\, P$, i.e, terms like $\langle
\widehat{\mathbb{C}}\,\text{sym} (\text{D}u - P), \text{sym} (\text{D}u - P)
\rangle$. This is the reason for which the crucial homogenization formula \eqref{eq:Ce_Cmac} for $L_{\rm c}\to 0$ can be obtained. Unlike for the linear Cosserat (micropolar) model, the relaxed micromorphic model remains operative and well posed \cite{neff2014unifying,d2022consistent,ghiba2015relaxed} also for zero Cosserat couple modulus $\mu_{\rm c}\equiv0$ ($\mathbb{C}_{\rm c}\equiv0$), in which case the force stress tensor remains symmetric.
The well-posedness is established using novel generalized Korn's inequalities for incompatible tensor fields \cite{lewintan2022lp,lewintan2021nevcas,lewintan2021korn,lewintan2021p,lewintan2021p,neff2015poincare,neff2012maxwell}, whereby sharp criteria for the validity of such coercivity estimates were given in the recent works \cite{GmeinederSpector,GLN1, GLN2}. In addition, the relaxed micromorphic model now operates as a true two-scale model between two clearly defined scales: the macroscopic scale with stiffness tensor $\mathbb{C}_{\rm macro}$ appearing for the characteristic length $L_{\rm c}\to 0$ (arbitrary large sample) and the microscopic scale with stiffness tensor $\mathbb{C}_{\rm micro}$ appearing for $L_{\rm c}\to \infty$.
Again, see Fig. \ref{Figure:intro}, the limit $L_{\rm c}\to \infty$ diverges as such in the classical micromorphic, second gradient, Cosserat model, along with others.

The above mentioned advantages have led to a multitude of investigations in short-time from the numerical side \cite{SkyHsymCurl,sky2021hybrid,sky2022primal,sky2022h,schroder2022lagrange,sarhil2022size}, from the modelling side \cite{rizzi2021exploring,voss2023modeling,rizzi2022boundary,rizzi2022metamaterial,ramirez2023multi,demore2022unfolding,rizzi2022towards,madeo2015wave}, analytical solutions \cite{rizzi2022torsion,rizzi2021uniaxial,rizzi2021cylindrical,rizzi2021shear}, regularity of solutions \cite{knees2023local,knees2023global}, and many others.

In this paper we continue our investigations from the theoretical side by determining the Green's functions for the case of a concentrated force and a concentrated couple in an infinite relaxed micromorphic medium.  Closed form solutions are derived using a Fourier transform analysis and results from generalized functions. It is shown that several well known generalized continuum fundamental solutions can be obtained as singular limiting cases of the relaxed micromorphic solution. In particular, from the relaxed micromorphic solutions we can readily derive the couple-stress, Cosserat-micropolar, micro-stretch, and classical elasticity fundamental solutions (\cite{mindlin1962effects,huilgol1967concentrated,sandru1966some,hattori2023isogeometric,khan1972singular,weitsman1967note,
mindlin1965stress,Dyszlewicz2004,mindlin1963influence,cowin1969singular,lakes2016physical,LiaHuang96,iecsan2001plane,timoshenko1970}), showing thus how versatile the relaxed micromorphic theory is. On the other hand, the full Eringen-Mindlin micromorphic model is at present too complicated for analytical or even numerical solutions to be sought. Here we take again advantage of the relaxed micromorphic model which drastically simplifies the situation in the isotropic planar case (only one curvature parameter remains operative). In the appendix we exhibit the two-scale elasticity nature relaxed micromorphic model. Moreover, we show how other generalised continua (micro-stretch, Cosserat-micropolar, couple stress) appear as limits of the relaxed micromorphic model.
%
%
%
%
%
%
\subsection{Notation}
For vectors $a,b\in\mathbb{R}^n$, we define the scalar product $\langle a,b \rangle \coloneqq \sum_{i=1}^n a_i\,b_i \in \mathbb{R}$, the (squared) euclidean norm  $\|a\|^2\coloneqq\langle a,a \rangle$ and  the dyadic product  $a\otimes b \coloneqq \left(a_i\,b_j\right)_{ij}\in \mathbb{R}^{n\times n}$. In the same way, for tensors  $P,Q\in\mathbb{R}^{n\times n}$, we define the scalar product $\langle P,Q \rangle \coloneqq\sum_{i,j=1}^n P_{ij}\,Q_{ij} \in \mathbb{R}$ and the (squared) Frobenius-norm $\|P\|^2\coloneqq\langle P,P \rangle$.
Moreover, $P^T$ denotes the transposition of the matrix $P$, which decomposes orthogonally into the skew-symmetric part $\text{skew} \, P \coloneqq \frac{1}{2} (P-P^T )$ and the symmetric part $\text{sym} \, P \coloneqq \frac{1}{2} (P+P^T)$.
The identity matrix is denoted by $\mathbbm{1}$, so that the trace of a matrix $P$ is given by \ $\text{tr} P \coloneqq \langle P,\mathbbm{1} \rangle$, while the deviatoric component of a matrix is given by $\text{dev} \, P \coloneqq P - \frac{\text{tr} \left( P\right)}{3} \, \mathbbm{1}$.
Given this, the  orthogonal decomposition possible for a matrix is $P = \text{dev} \,\text{sym} \, P + \text{skew} \, P + \frac{\text{tr} \left( P\right)}{3} \, \mathbbm{1}$.
The Lie-Algebra of skew-symmetric matrices is denoted by $\mathfrak{so}(3)\coloneqq \{A\in\mathbb{R}^{3\times 3}\mid A^T = -A\}$.
The derivative $\text{D}u$ and the curl of a vector field $u$ are defined as
\begin{equation}
\text{D}u
=
\left(
\begin{array}{ccc}
u_{1,1} & u_{1,2} & u_{1,3} \\
u_{2,1} & u_{2,2} & u_{2,3} \\
u_{3,1} & u_{3,2} & u_{3,3}
\end{array}
\right)\, ,
\qquad
\text{curl} \, u = \nabla \times u
=
\left(
\begin{array}{ccc}
u_{3,2} - u_{2,3}  \\
u_{1,3} - u_{3,1}  \\
u_{2,1} - u_{1,2}
\end{array}
\right)
\, .
\end{equation}
We also introduce the $\text{Curl}$ and the $\text{Div}$ operators for $P\in\mathbb{R}^{3\times 3}$ as
\begin{equation}
\text{Curl} \, P
=\!
\left(
\begin{array}{c}
(\text{curl}\left( P_{11} , \right. P_{12} , \left. P_{13} \right))^T \\
(\text{curl}\left( P_{21} , \right. P_{22} , \left. P_{23} \right))^T \\
(\text{curl}\left( P_{31} , \right. P_{32} , \left. P_{33} \right))^T
\end{array}
\right) \!,
\qquad
\text{Div}  \, P
=\!
\left(
\begin{array}{c}
\text{div}\left( P_{11} , \right. P_{12} , \left. P_{13} \right)^{T} \\
\text{div}\left( P_{21} , \right. P_{22} , \left. P_{23} \right)^{T} \\
\text{div}\left( P_{31} , \right. P_{32} , \left. P_{33} \right)^{T}
\end{array}
\right) \, .
\end{equation}
With these definitions, we have consistently $\text{Curl}\,\text{D}u=0$.
The cross product between a second order tensor and a vector is also needed and is defined row-wise as follows
\begin{equation}
m \times b =
\left(
\begin{array}{ccc}
(b \times (m_{11},m_{12},m_{13}))^T \\
(b \times (m_{21},m_{22},m_{23}))^T \\
(b \times (m_{31},m_{32},m_{33}))^T \\
\end{array}
\right) =
m \cdot \epsilon \cdot b =
m_{ik} \, \epsilon_{kjh} \, b_{h} 
\, ,
\end{equation}
where $m \in \mathbb{R}^{3\times 3}$, $b \in \mathbb{R}^{3}$, and $\epsilon$ is the Levi-Civita tensor.
%
%
%
%
%
%
\section{The isotropic relaxed micromorphic model}
The isotropic relaxed micromorphic model has the kinematics of the classical Eringen-Mindlin micromorphic isotropic model \cite{Mindlin1964,Eringen1999}, i.e. the displacement $u \in \mathbb{R}^{3}$ and the non-symmetric micro-distortion $P \in \mathbb{R}^{3\times3}$ as independent fields.
The strain energy density reads
\begin{align}
W \left(\text{D}u, P,\text{Curl}\,P\right) = &
\, \mu_{\rm e} \left\lVert \text{sym} (\text{D}u - P ) \right\rVert^{2}
+ \mu_{\rm c} \left\lVert \text{skew} (\text{D}u - P ) \right\rVert^{2}
+ \frac{\lambda_{\rm e}}{2} \text{tr}^2 (\text{D}u - P )
\notag
\\*
&
+ \mu_{\rm micro} \left\lVert \text{sym}\,P \right\rVert^{2}
+ \frac{\lambda_{\rm micro}}{2} \text{tr}^2 \left(P \right)
\label{eq:energy_RM_3}
\\*
&
+ \frac{ \mu_{\rm macro} L_{\rm c}^2}{2} \,
\left(
a_1 \, \left\lVert \text{dev sym} \, \text{Curl} \, P\right\rVert^2 +
a_2 \, \left\lVert \text{skew} \,  \text{Curl} \, P\right\rVert^2 +
\frac{a_3}{3} \, \text{tr}^2 \left(\text{Curl} \, P\right)
\right),
\notag
\end{align}
while the two equilibrium equations are
\begin{align}
\text{Div} \, \sigma = f \, ,
\qquad\qquad\qquad\qquad
\sigma - \sigma_{\rm micro}- \text{Curl}\,m &= M \, ,
\label{eq:equi_RM}
\end{align}
with
\begin{align}
\label{constEq}
\sigma & \coloneqq
2\mu_{\rm e}\,\text{sym}  (\text{D}u - P )
+ 2\mu_{\rm c}\,\text{skew} (\text{D}u - P )
+ \lambda_{\rm e} \text{tr} (\text{D}u - P ) \mathbbm{1}
\, ,
\notag
\\*
\sigma_{\rm micro} & \coloneqq
2 \mu_{\rm micro}\,\text{sym}\,P
+ \lambda_{\rm micro} \text{tr} \left( P \right) \mathbbm{1}
\, ,
\\*
\notag
m & \coloneqq
\mu_{\rm macro}\, L_{\rm c}^2 \left( a_1 \, \text{dev sym} \, \text{Curl} \, P +
a_2 \, \text{skew} \, \text{Curl} \, P +
\frac{a_3}{3} \, \text{tr} \left(\text{Curl} \, P\right)\mathbbm{1}\right)
\, ,
\end{align}
where $\sigma$ is the non-symmetric elastic force stress tensor, $m$ is the non-symmetric moment tensor, $f$ is the standard body force vector and $M$ is the body volume couple tensor.
The homogeneous Neumann and the Dirichlet boundary conditions are  
\begin{align}
&
\text{Neumann:}
&&&&
t
\coloneqq
\sigma \, n
= 
0
\, ,
&&
\text{and}
&&
\eta
\coloneqq
m \times n = 
0
\, ,
&
\label{eq:BC_RM_Neu}
\\*
&
\text{Dirichlet:}
&&&&
u
=
\overline{u}
\, ,
&&
\text{and}
&&
\overline{Q}=P \times n 
\, ,
&
\label{eq:BC_RM_Dir}
\end{align}
where the higher-order Dirichlet boundary conditions in  (\ref{eq:BC_RM_Dir}) can be particularised to
\begin{align}
P \times n
=
\overline{Q}
=
\text{D}u \times n
\, ,
\label{eq:CCBC_RM}
\end{align}
formally called ``consistent coupling boundary conditions'' \cite{d2022consistent}.
The following additional (but not independent) equilibrium equation can be derived combining the two equilibrium equations~(\ref{eq:equi_RM}) based on the fundamental property of differential operators $\text{Div}\,\text{Curl}\,(\cdot) =0$
\begin{align}
\text{Div}\,\sigma_{\rm micro} = f-\text{Div}\,M \, .
\end{align}
A similar additional equilibrium equation for $\sigma_{\rm micro}$ does not exist in the classical Eringen-Mindlin micromorphic model or the Cosserat model.
\section{The isotropic relaxed micromorphic model in plane-strain}

Under the plane-strain hypothesis only the in-plane components of the kinematic fields are different from zero and they only depend on $(x_1,x_2)$.
The structure of the kinematic fields ($\widetilde{u}$,$\widetilde{P}$) are \cite{iecsan2001plane}
\begin{align}
\widetilde{u}
=
\left(
\begin{array}{ccc}
u_{1} \\ 
u_{2} \\
0
\end{array}
\right)
\, ,
&&
\widetilde{u}^{\sharp}
=
\left(
\begin{array}{ccc}
u_{1} \\ 
u_{2}
\end{array}
\right)
\, ,
&&
\widetilde{P}
=
\left(
\begin{array}{ccc}
P_{11} & P_{12} & 0 \\ 
P_{21} & P_{22} & 0 \\ 
0 & 0 & 0
\end{array}
\right)
\, ,
&&
\widetilde{P}^{\sharp}
=
\left(
\begin{array}{ccc}
P_{11} & P_{12}\\
P_{21} & P_{22}
\end{array}
\right)
\, ,
\label{eq:gen_kine}
\end{align}
while the structures of the gradient of the displacement $\text{D}\widetilde{u}$, of the Curl of the micro distortion tensor $\text{Curl}\,\widetilde{P}$, and of the double Curl of the micro distortion tensor $\text{Curl}\,\text{Curl}\,\widetilde{P}$ are
\begin{equation}
\begin{gathered}
\text{D}\widetilde{u}
=
\left(
\begin{array}{ccc}
u_{1,1} & u_{1,2} & 0 \\
u_{2,1} & u_{2,2} & 0 \\
0 & 0 & 0 \\
\end{array}
\right)
\, ,
\qquad
\text{D}\widetilde{u}^{\sharp}
=
\left(
\begin{array}{ccc}
u_{1,1} & u_{1,2} \\
u_{2,1} & u_{2,2}
\end{array}
\right)
\, ,
\\*[5pt]
\text{Curl} \,\widetilde{P}
= 
\left(
\begin{array}{cc|c}
0 & 0 & P_{12,1}-P_{11,2} \\
0 & 0 & P_{22,1}-P_{21,2} \\
\hline
0 & 0 & 0 \\
\end{array}
\right)
= 
\left(
\begin{array}{cc|c}
0 & 0 & \\
0 & 0 & \multicolumn{1}{c}{\smash{\raisebox{.5\normalbaselineskip}{$\text{Curl}_{\text{2D}}\,\widetilde{P}^{\sharp}$}}} \\
\hline
0 & 0 & 0 \\
\end{array}
\right)
\, , \qquad 
\text{Curl}_{\text{2D}} \,\widetilde{P}^{\sharp}
\coloneqq
\left(
\begin{array}{ccc}
P_{12,1}-P_{11,2} \\
P_{22,1}-P_{21,2}
\end{array}
\right)
\, ,
\\*[5pt]
\text{Curl} \,\text{Curl} \,\widetilde{P}
= 
\left(
\begin{array}{cc|c}
P_{12,12}-P_{11,22} &
P_{11,12}-P_{12,11} & 0 \\
P_{22,12}-P_{21,22} &
P_{21,12}-P_{22,11} & 0 \\
\hline
0 & 0 & 0 \\
\end{array}
\right)
= 
\left(
\begin{array}{cc|c}
&   & 0 \\
\multicolumn{2}{c|}{\smash{\raisebox{.5\normalbaselineskip}{$\text{Curl} \, \text{Curl}_{\text{2D}}\,\widetilde{P}^{\sharp}$}}} & 0 \\
\hline
0 & 0 & 0 \\
\end{array}
\right)
\, ,
\\*[5pt]
\text{Curl} \, \text{Curl}_{\text{2D}}\,\widetilde{P}^{\sharp}\coloneqq
\left(
\begin{array}{cc}
P_{12,12}-P_{11,22} &
P_{11,12}-P_{12,11} \\
P_{22,12}-P_{21,22} &
P_{21,12}-P_{22,11}
\end{array}
\right)
\, .
\label{eq:gen_kine_deriv}
\end{gathered}
\end{equation}
The operator $\text{Curl}_{\text{2D}}$ is a rotated divergence and 
\begin{align}
    \text{Curl}_{\text{2D}} \widetilde{P}^{\sharp} &= 
    \text{Div} (\widetilde{P}^{\sharp} R^T) \, , &&  R = \begin{bmatrix}
        0 & 1 \\
        -1 & 0
    \end{bmatrix} \, , \notag \\
    \text{Curl}\text{Curl}_{\text{2D}} \widetilde{P}^{\sharp} &=  (\text{D} \text{Curl}_{\text{2D}} \widetilde{P}^{\sharp}) R^T = (\text{D} \text{Div} [\widetilde{P}^{\sharp} R^T]) R^T \, .    
\end{align}

Because of the nature of the Curl operator, it is noted that $\text{Curl} \,P$ just has out of plane components that depend on the in-plane components of $P$, while $\text{Curl} \,\text{Curl} \,\widetilde{P}$ fully preserves the in-plane structure.
Moreover, since 
\begin{gather}
\text{tr}(\text{Curl} \,\widetilde{P})=0,
\\*[5pt]
\left\lVert \text{dev sym} \, \text{Curl} \, \widetilde{P}\right\rVert^2
=\left\lVert \text{sym} \, \text{Curl} \, \widetilde{P}\right\rVert^2
=\left\lVert \text{skew} \,  \text{Curl} \, \widetilde{P}\right\rVert^2
=\frac{1}{2}\left\lVert\text{Curl} \, \widetilde{P}\right\rVert^2
=\frac{1}{2}\left\lVert\text{Curl}_{\text{2D}} \, \widetilde{P}^{\sharp}\right\rVert^2
\, ,
\notag
\end{gather}
the plane strain isotropic relaxed micromorphic model will just depend on one cumulative dimensionless parameter $\widetilde{a}\coloneqq\frac{a_1+a_2}{2}$, and the strain energy density in  (\ref{eq:energy_RM_3}) reduces to
\begin{align}
W (\text{D}\widetilde{u}, \widetilde{P},\text{Curl}\,\widetilde{P})
= &
\, \mu_{\rm e} \left\lVert \text{sym} (\text{D}\widetilde{u} - \widetilde{P} ) \right\rVert^{2}
+ \mu_{\rm c} \left\lVert \text{skew} (\text{D}\widetilde{u} - \widetilde{P} ) \right\rVert^{2}
+ \frac{\lambda_{\rm e}}{2} \text{tr}^2 (\text{D}\widetilde{u} - \widetilde{P} )
\notag\\*
&
+ \mu_{\rm m} \left\lVert \text{sym}\,\widetilde{P} \right\rVert^{2}
+ \frac{\lambda_{\rm m}}{2} \text{tr}^2 (\widetilde{P})
+ \frac{\mu_{\rm M}L_{\rm c}^2 }{2} \, \widetilde{a} \,
\left\lVert \text{Curl} \, \widetilde{P}\right\rVert^2 
\, ,
\notag
\\*
= &
\, \mu_{\rm e} \left\lVert \text{sym} (\text{D}\widetilde{u}^{\sharp} - \widetilde{P}^{\sharp} ) \right\rVert^{2}
+ \mu_{\rm c} \left\lVert \text{skew} (\text{D}\widetilde{u}^{\sharp} - \widetilde{P}^{\sharp} ) \right\rVert^{2}
+ \frac{\lambda_{\rm e}}{2} \text{tr}^2 (\text{D}\widetilde{u}^{\sharp} - \widetilde{P}^{\sharp} )
\label{eq:energy_RM_2D}
\\*
&
+ \mu_{\rm m} \left\lVert \text{sym}\,\widetilde{P}^{\sharp} \right\rVert^{2}
+ \frac{\lambda_{\rm m}}{2} \text{tr}^2 (\widetilde{P}^{\sharp})
+ \frac{\mu_{\rm M}L_{\rm c}^2 }{2} \, \widetilde{a} \,
\left\lVert \text{Curl}_{\text{2D}} \, \widetilde{P}^{\sharp}\right\rVert^2 
\notag\\
= &
\, \mu_{\rm e} \left\lVert \text{dev}_2\,\text{sym} (\text{D}\widetilde{u}^{\sharp} - \widetilde{P}^{\sharp} ) \right\rVert^{2}
+ \mu_{\rm c} \left\lVert \text{skew} (\text{D}\widetilde{u}^{\sharp} - \widetilde{P}^{\sharp} ) \right\rVert^{2}
+ \frac{\kappa_{\rm e}}{2} \text{tr}^2 (\text{D}\widetilde{u}^{\sharp} - \widetilde{P}^{\sharp} )
\notag\\*
&
+ \mu_{\rm m} \left\lVert \text{dev}_2\,\text{sym}\,\widetilde{P}^{\sharp} \right\rVert^{2}
+ \frac{\kappa_{\rm m}}{2} \text{tr}^2 (\widetilde{P}^{\sharp})
+ \frac{\mu_{\rm M}\, L_{\rm c}^2 }{2} \, \widetilde{a} \,
\left\lVert \text{Curl}_{\text{2D}} \, \widetilde{P}^{\sharp}\right\rVert^2 \,,\notag
\end{align}
where $\text{dev}_2 X \coloneqq X -\frac12\text{tr}(X)\cdot\mathbbm{1}_2$. Also, for better readability we employ the following abbreviated forms: $\mu_{\rm M} \equiv \mu_{\rm macro}$, $\mu_{\rm m} \equiv \mu_{\rm micro}$,  $\lambda_{\textrm{m}} \equiv \lambda_{\textrm{micro}}$ and $\lambda_{\textrm{M}} \equiv \lambda_{\textrm{macro}}$. Moreover, 
under plane-strain conditions, the bulk micro-moduli $\kappa_{\rm e}$ and $\kappa_{\rm m} \equiv\kappa_{\rm micro}$ are related with the respective Lamé type micro-moduli through the 2D relations
\begin{equation}\label{eq:KAPPA}
\kappa_{\rm e}\coloneqq \lambda_{\rm e}+\mu_{\rm e}, \qquad \quad \kappa_{\rm m} \coloneqq\lambda_{\rm m}+\mu_{\rm m}\,.
\end{equation}
Accordingly, the relations between the macro moduli ($\mu_{\rm M},\lambda_{\rm M}, \kappa_{\rm M}$) and the micro-moduli in plane strain become (see Appendix A2)
\begin{align}
\label{defmod}
\mu_{\rm M} &\coloneqq \dfrac{\mu_{\rm e} \, \mu_{\rm m}}{\mu_{\rm e} + \mu_{\rm m}} 
\qquad\Leftrightarrow\qquad \frac{1}{\mu_{\rm M}}=\frac{1}{\mu_{\rm e}}+\frac{1}{\mu_{\rm m}}\,,\notag\\
\kappa_{\rm M} &\coloneqq \dfrac{\kappa_{\rm e} \, \kappa_{\rm m}}{\kappa_{\rm e} + \kappa_{\rm m}} \qquad\Leftrightarrow\qquad
 \frac{1}{\kappa_{\rm M}}=\frac{1}{\kappa_{\rm e}}+\frac{1}{\kappa_{\rm m}}\,,\\
\lambda_{\rm m} &\coloneqq
\dfrac{(\mu_{\rm e} + \lambda_{\rm e})(\mu_{\rm m} + \lambda_{\rm m})}{(\mu_{\rm e} + \lambda_{\rm e}) + (\mu_{\rm m} + \lambda_{\rm m})} -\dfrac{\mu_{\rm e} \, \mu_{\rm m}}{\mu_{\rm e} + \mu_{\rm m}}
\, ,
\notag
\end{align}
where $\kappa_{\rm M} \equiv \kappa_{\rm macro}$ with $\kappa_{\rm M}= \mu_{\rm M}+\lambda_{\rm M} $. The 3D relations for the macro and micro moduli are given in Appendix A. From here and onwards, unless otherwise stated, the macro and micro moduli will refer to the case of plane strain and will be defined through equations  \eqref{eq:KAPPA} and \eqref{defmod}.

Taking the first variation of the strain energy $I=\displaystyle\int_{\Omega} W \, \text{d} x$ under the plane strain hypothesis with respect to ($\widetilde{u}^{\sharp}$,$\widetilde{P}^{\sharp}$) leads to
\begin{align}
\delta I^{\widetilde{u}^{\sharp}}
&=
\int_{\Omega}
\hspace{-0.1cm}
\Big(
2\mu_{\rm e} \,
\langle \text{sym} (\text{D}\widetilde{u}^{\sharp} - \widetilde{P}^{\sharp} ),\text{D}\delta \widetilde{u}^{\sharp} \rangle
+ 2\mu_{\rm c} \,
\langle \text{skew} (\text{D}\widetilde{u}^{\sharp} - \widetilde{P}^{\sharp} ),\text{D}\delta \widetilde{u}^{\sharp} \rangle
+ \lambda_{\rm e} \langle \text{tr} (\text{D}\widetilde{u}^{\sharp} - \widetilde{P}^{\sharp} ) \mathbbm{1}_2, \text{D}\delta \widetilde{u}^{\sharp} \rangle
\Big)
\text{d}x
\, , \notag
\\*[5pt]
\delta I^{\widetilde{P}^{\sharp}}
&=
\int_{\Omega}
\hspace{-0.1cm}
\Big(
-2\mu_{\rm e} \,
\langle \text{sym} (\text{D}\widetilde{u}^{\sharp} - \widetilde{P}^{\sharp} ),\delta \widetilde{P}^{\sharp} \rangle
-2\mu_{\rm c} \,
\langle \text{skew} (\text{D}\widetilde{u}^{\sharp} - \widetilde{P}^{\sharp} ),\delta \widetilde{P}^{\sharp} \rangle
-\lambda_{\rm e} \langle \text{tr} (\text{D}\widetilde{u}^{\sharp} - \widetilde{P}^{\sharp} ) \mathbbm{1}_2, \delta \widetilde{P}^{\sharp} \rangle
\label{eq:first_varia_energy_RM_2D_P}
\\
&
\phantom{=\int_{\Omega}\Big(}
+
2\mu_{\rm m} \langle \text{sym}\,\widetilde{P}^{\sharp},\delta \widetilde{P}^{\sharp} \rangle
+ \lambda_{\rm m} \langle \text{tr} (\widetilde{P}^{\sharp} ) \mathbbm{1}_2, \delta \widetilde{P}^{\sharp} \rangle
+ \mu_{\rm M}\, L_{\rm c}^2 \, \widetilde{a} \,
\langle \text{Curl}_{\text{2D}} \, \widetilde{P}^{\sharp}, \text{Curl}_{\text{2D}} \, \delta \widetilde{P}^{\sharp} \rangle
\Big)
\text{d}x
\, .
\notag
\end{align}
The equilibrium equation are now obtained by requiring 
\begin{equation}
\delta I^{\widetilde{u}^{\sharp}}= \langle \widetilde{f} , \delta \widetilde{u}^{\sharp} \rangle \, ,
\quad
\forall
\,
\delta \widetilde{u}^{\sharp}
\qquad\qquad
\text{and}
\qquad\qquad
\delta I^{\widetilde{P}^{\sharp}}= \langle \widetilde{M} , \delta \widetilde{P}^{\sharp} \rangle \, ,
\quad
\forall
\,
\delta \widetilde{P}^{\sharp}  \, .
\label{eq:first_variation_2D}
\end{equation}
We define the following quantities
\begin{align}
\label{2Dconst}
\widetilde{\sigma}&\coloneqq
2\mu_{\rm e} \, \text{sym} (\text{D}\widetilde{u}^{\sharp} - \widetilde{P}^{\sharp} )
+ 2\mu_{\rm c} \, \text{skew} (\text{D}\widetilde{u}^{\sharp} - \widetilde{P}^{\sharp} )
+ \lambda_{\rm e} \, \text{tr} (\text{D}\widetilde{u}^{\sharp} - \widetilde{P}^{\sharp} ) \mathbbm{1}_2
\, ,
\notag
\\*
\widetilde{\sigma}_{\rm m}&\coloneqq
2\mu_{\rm m} \, \text{sym}\,\widetilde{P}^{\sharp}
+ \lambda_{\rm m} \, \text{tr} (\widetilde{P}^{\sharp}) \mathbbm{1}_2
\in \mathbb{R}^{2\times2}\, ,
\\*
\widetilde{m}&\coloneqq
\mu_{\rm M}\, L_{\rm c}^2\, \widetilde{a}\, \text{Curl}_{\text{2D}} \widetilde{P}^{\sharp}
\in \mathbb{R}^{2}\, ,
\notag
\end{align}
where we used the tilde $\widetilde{\sigma}$ and $\widetilde{\sigma}_{\rm m} \equiv \widetilde{\sigma}_{\rm micro}$ to emphasize that here we are only considering the in-plane components.
We can rewrite the first variation $\delta I^{\widetilde{u}}$ as
\begin{align}
\delta I^{\widetilde{u}^{\sharp}}
&=
\int_{\Omega}
\hspace{-0.1cm}
\langle
\widetilde{\sigma}
,
\text{D}\delta \widetilde{u}^{\sharp}
\rangle
\,
\text{d}x
=
\int_{\Omega}
\hspace{-0.1cm}
\text{div}
(
\widetilde{\sigma}^{\rm T}
\, \delta \widetilde{u}^{\sharp}
)
-
\langle 
\text{Div} \, 
\widetilde{\sigma}
,
\delta \widetilde{u}^{\sharp}
\rangle
\,
\text{d}x
=
\int_{\partial \Omega}
\hspace{-0.1cm}
\langle
\widetilde{\sigma}^{\rm T}
\, \delta \widetilde{u}^{\sharp}
,
n^{\sharp}
\rangle
\,
\text{d}s
-
\int_{\Omega}
\hspace{-0.1cm}
\langle 
\text{Div} \, 
\widetilde{\sigma}
,
\delta \widetilde{u}^{\sharp}
\rangle
\,
\text{d}x
\label{eq:first_varia_energy_RM_2D_u_2}
\\*[5pt]
&=
\int_{\partial \Omega}
\hspace{-0.1cm}
\langle
\widetilde{\sigma}
\, n^{\sharp}
,
\delta \widetilde{u}^{\sharp}
\rangle
\,
\text{d}s
-
\int_{\Omega}
\hspace{-0.1cm}
\langle 
\text{Div} \, 
\widetilde{\sigma}
,
\delta \widetilde{u}^{\sharp}
\rangle
\,
\text{d}x \, ,
\notag
\end{align}
which, because of the equation (\ref{eq:first_variation_2D}), and highlighting that $\widetilde{u}$ is orthogonal with respect to the out-of-plane displacement, implies that
\begin{align}
\text{Div} \, 
\widetilde{\sigma}
= \widetilde{f}
\quad
\text{in }
\,
\Omega \, ,
\qquad\qquad
\widetilde{\sigma} \, n^{\sharp} = 0
\quad
\text{on }
\,
\partial \Omega \, .
\label{eq:equa_var_u}
\end{align}
where the out-of-plane components of $\text{Div}\, \widetilde{\sigma}$ and $\widetilde{\sigma} \, n$ must not be considered, and $n^{\sharp}$ is the vector of the in-plane component of the normal to the boundary.
We recall that the tangent vector in two dimensions is unique and can be obtained as $t^\sharp = R^T n^\sharp$ where the following product rule holds:
\begin{align}
    \text{div}[(\delta \widetilde{P}^\sharp R^T)^T v] &= \langle\text{Div}(\delta \widetilde{P}^\sharp R^T), v \rangle + \langle \delta \widetilde{P}^\sharp R^T, \text{D}v \rangle
    = \langle\text{Curl}_\text{2D}\delta \widetilde{P}^\sharp, v \rangle + \langle  \delta \widetilde{P}^\sharp , (\text{D}v) R \rangle  \\
    &= \langle\text{Curl}_\text{2D}\delta \widetilde{P}^\sharp, v \rangle - \langle  \delta \widetilde{P}^\sharp , (\text{D}v) R^T \rangle \,, \notag
\end{align}
which implies 
\begin{align}
    \langle\text{Curl}_\text{2D}\delta \widetilde{P}^\sharp, v \rangle =  \langle \delta \widetilde{P}^\sharp , (\text{D}v) R^T \rangle  +  \text{div}(R (\delta \widetilde{P}^\sharp)^T v)\, ,
\end{align}
for some smooth vector field $v \in \mathbb{R}^2$.
Thus, there holds the Green-type identity 
\begin{align}
    \int_\Omega \langle\text{Curl}_\text{2D} \widetilde{P}^\sharp, \text{Curl}_\text{2D}\delta \widetilde{P}^\sharp \rangle \, \text{d}x &= \int_\Omega \langle  (\text{D} \text{Curl}_\text{2D} \widetilde{P}^\sharp) R^T, \delta \widetilde{P}^\sharp \rangle \, \text{d}x + \int_\Omega \text{div}[R (\delta \widetilde{P}^\sharp)^T \text{Curl}_\text{2D} P^\sharp] \, \text{d}x \notag \\
    &= \int_\Omega \langle  \text{Curl} \text{Curl}_\text{2D} P^\sharp, \delta \widetilde{P}^\sharp \rangle \, \text{d}x + \int_{\partial \Omega} \langle R (\delta \widetilde{P}^\sharp)^T \text{Curl}_\text{2D} \widetilde{P}^\sharp, n^\sharp \rangle \, \text{d}s  \\
    &= \int_\Omega \langle  \text{Curl} \text{Curl}_\text{2D} \widetilde{P}^\sharp, \delta \widetilde{P}^\sharp \rangle \, \text{d}x + \int_{\partial \Omega} \langle  \text{Curl}_\text{2D} \widetilde{P}^\sharp, \delta \widetilde{P}^\sharp R^T n^\sharp \rangle \, \text{d}s \notag \\
    &= \int_\Omega \langle  \text{Curl} \text{Curl}_\text{2D} \widetilde{P}^\sharp, \delta \widetilde{P}^\sharp \rangle \, \text{d}x + \int_{\partial \Omega} \langle  \text{Curl}_\text{2D} \widetilde{P}^\sharp, \delta \widetilde{P}^\sharp t^\sharp \rangle \, \text{d}s \, , \notag
\end{align}
where we simply substituted $v$ for $\text{Curl} \widetilde{P}^\sharp$ in the above product rule and applied the divergence theorem. Using the latter, we can now rewrite the first variation $\delta I^{\widetilde{P}^{\sharp}}$ as
\begin{align}
\delta I^{\widetilde{P}^{\sharp}}
&=
\int_{\Omega}
\hspace{-0.1cm}
-
\langle \widetilde{\sigma},\delta \widetilde{P}^{\sharp} \rangle
+
\langle \widetilde{\sigma}_{\rm m}, \delta \widetilde{P}^{\sharp} \rangle
+
\langle \mu_{\rm M}\, L_{\rm c}^2\widetilde{a}\,\text{Curl}_{\text{2D}}\,\widetilde{P}^{\sharp}, \text{Curl}_{\text{2D}} \, \delta \widetilde{P}^{\sharp} \rangle
\,
\text{d}x
\notag
\\*
&=
\int_{\Omega}
\hspace{-0.1cm}
\langle
-
\widetilde{\sigma}
+
\widetilde{\sigma}_{\rm m}
,
\delta \widetilde{P}^{\sharp}
\rangle
+
\langle \mu_{\rm M}\, L_{\rm c}^2\widetilde{a}\,\text{Curl}_{\text{2D}}\,\widetilde{P}^{\sharp}, \text{Curl}_{\text{2D}} \, \delta \widetilde{P}{\sharp} \rangle
\,
\text{d}x
\notag
\\*
&=
\int_{\Omega}
\hspace{-0.1cm}
\langle
-
\widetilde{\sigma}
+
\widetilde{\sigma}_{\rm m}
,
\delta \widetilde{P}^{\sharp}
\rangle
+
\langle \mu_{\rm M}\, L_{\rm c}^2\widetilde{a}\,\text{Curl} \, \text{Curl}_{\text{2D}} \, \widetilde{P}^{\sharp}, \delta \widetilde{P}^{\sharp} \rangle
+
\text{div} [ R (\delta \widetilde{P}^{\sharp})^T  (\mu_{\rm M}L_{\rm c}^2\widetilde{a}\,\text{Curl}_{\text{2D}}\,\widetilde{P}^{\sharp}) ] 
\,
\text{d}x
\label{eq:first_varia_energy_RM_2D_P_2}
\\*
&=
\int_{\Omega}
\hspace{-0.1cm}
\langle
-
\widetilde{\sigma}
+
\widetilde{\sigma}_{\rm m}
+
\mu_{\rm M}\, L_{\rm c}^2\widetilde{a}\,\text{Curl} \, \text{Curl}_{\text{2D}}\,\widetilde{P}^{\sharp}
,
\delta \widetilde{P}^{\sharp}
\rangle
+
\int_{\partial \Omega}
\hspace{-0.1cm}
\langle R  (\delta \widetilde{P}^{\sharp})^T   \mu_{\rm M}L_{\rm c}^2\widetilde{a}\,\text{Curl}_{\text{2D}}\,\widetilde{P}^{\sharp},    n^{\sharp} \rangle 
\,
\text{d}s
\notag
\\*
&=
\int_{\Omega}
\hspace{-0.1cm}
\langle
-
\widetilde{\sigma}
+
\widetilde{\sigma}_{\rm m}
+
\mu_{\rm M}\, L_{\rm c}^2\widetilde{a}\,\text{Curl} \, \text{Curl}_{\text{2D}}\,\widetilde{P}^{\sharp}
,
\delta \widetilde{P}^{\sharp}
\rangle
+
\int_{\partial \Omega}
\hspace{-0.1cm}
\langle  \mu_{\rm M}L_{\rm c}^2\widetilde{a}\,\text{Curl}_{\text{2D}}\,\widetilde{P}^{\sharp},   \delta \widetilde{P}^{\sharp} t^{\sharp} \rangle 
\,
\text{d}s
\, ,
\notag
\end{align}
which, because of (\ref{eq:first_variation_2D}),
implies
\begin{align}
\widetilde{\sigma}
-
\widetilde{\sigma}_{\rm m}
-
\mu_{\rm M} L_{\rm c}^2\widetilde{a}\,\text{Curl} \, \text{Curl}_{\text{2D}} \, \widetilde{P}^{\sharp} = \widetilde{M}
\quad
\text{in }
\,
\Omega \, ,
\qquad\qquad
\displaystyle \mu_{\rm M}\, L_{\rm c}^2\widetilde{a}\,\text{Curl}_{\text{2D}}\,\widetilde{P}^{\sharp} = 0
\quad
\text{on }
\,
\partial \Omega \, .
\label{eq:equa_var_P}
\end{align}
where the out-of-plane components of (\ref{eq:equa_var_P})$_1$ and (\ref{eq:equa_var_P})$_2$ must not be considered.
We can now collect all the homogeneous equilibrium equations obtained and the homogeneous Neumann boundary conditions
\begin{align}
\left.
\begin{array}{rrr}
\text{Div}\, \widetilde{\sigma} = \widetilde{f}
\\*
\widetilde{\sigma} - \widetilde{\sigma}_{\rm m} - \mu_{\rm M} L_{\rm c}^2\widetilde{a}\,\text{Curl} \, \text{Curl}_{\text{2D}} \, \widetilde{P}^\sharp = \widetilde{M}
\end{array}
\right\}
\quad
\text{in }
\,
\Omega
\, ,
\qquad\qquad
\left.
\begin{array}{rrr}
\widetilde{\sigma} \, n^{\sharp} = 0
\\*
\mu_{\rm M} L_{\rm c}^2\widetilde{a}\,\text{Curl}_{\text{2D}}\,\widetilde{P}^\sharp = 0
\end{array}
\right\}
\quad
\text{on }
\,
\partial \Omega
\, .
\label{eq:equi_RM_2D}
\end{align}
Since $\text{Div} \, (\mu_{\rm M}\, L_{\rm c}^2\widetilde{a}\,\text{Curl} \, \text{Curl}_{\text{2D}}\,\widetilde{P}^\AS{\sharp}) =0$, combining the two equations in (\ref{eq:equi_RM_2D})$_1$ gives rise to another equilibrium equation that depends only on $\text{sym}\,\widetilde{P}$
\begin{align}
\left.
\begin{array}{rrr}
\text{Div}\, \widetilde{\sigma} = \widetilde{f}
\\*
\widetilde{\sigma} - \widetilde{\sigma}_{\rm m} - \mu_{\rm M}\, L_{\rm c}^2\widetilde{a}\,\text{Curl} \, \text{Curl}_{\text{2D}} \, \widetilde{P}^\AS{\sharp} = \widetilde{M}
\\*
\{
\text{Div}\, \widetilde{\sigma}_{\rm m} = \widetilde{f} + \text{Div}\, \widetilde{M}
\}
\end{array}
\right\}
\quad
\text{in }
\,
\Omega
\, ,
\qquad\qquad
\left.
\begin{array}{rrr}
\widetilde{\sigma} \, n^{\sharp} = 0
\\*
\mu_{\rm M} L_{\rm c}^2\widetilde{a}\,\text{Curl}_{\text{2D}}\,\widetilde{P}^\sharp = 0
\end{array}
\right\}
\quad
\text{on }
\,
\partial \Omega
\, .
\label{eq:equi_RM_2D_2}
\end{align}
It should be noted however that the additional equation $\text{Div}\, \widetilde{\sigma}_{\rm m} = \widetilde{f} + \text{Div}\, \widetilde{M}$ is not independent with respect to the other two, and any  solution of (\ref{eq:equi_RM_2D})$_1$ will automatically satisfy it.
The governing equilibrium equations~(\ref{eq:equi_RM_2D})$_1$ in components become then
\begin{align}\begin{split}
\label{pstrain}
(\lambda_{\rm e}+2 \mu_{\rm e})
\left(u_{1,11}-P_{11,1}\right)
+\lambda_{\rm e} \left(u_{2,12}-P_{22,1}\right)
\quad
&
\\*
+\mu_{\rm c} \left(u_{1,22}-u_{2,12}-P_{12,2}+P_{21,2}\right)
+\mu_{\rm e} \left(u_{1,22}+u_{2,12}-P_{12,2}-P_{21,2}\right)
&=f_1 \, ,
\\*[5pt]
(\lambda_{\rm e}+2 \mu_{\rm e}) \left(u_{2,22}-P_{22,2}\right)
+\lambda_{\rm e} \left(u_{1,12}-P_{11,2}\right)
\quad
&
\\*
+\mu_{\rm c} \left(u_{2,11}-u_{1,12}-P_{21,1}+P_{12,1}\right)
+\mu_{\rm e} \left(u_{2,11}+u_{1,12}-P_{21,1}-P_{12,1}\right)
&=f_2 \, ,
\\*[10pt]
\mu_{\rm M}\, L_{\rm c}^2\, \widetilde{a} \left(P_{11,22}-P_{12,12}\right)
\quad
&
\\*
-P_{11} (\lambda_{\rm e}+\lambda_{\rm m}+2 (\mu_{\rm e}+\mu_{\rm m}))
-(\lambda_{\rm e}+\lambda_{\rm m}) P_{22}
+(\lambda_{\rm e}+2 \mu_{\rm e}) u_{1,1}+
\lambda_{\rm e} u_{2,2}
&=M_{11} \, ,
\\*[5pt]
-\mu_{\rm M}\, L_{\rm c}^2\, \widetilde{a} \left(P_{11,12}-P_{12,11}\right)
\quad
&
\\*
-(\mu_{\rm c}+\mu_{\rm e}+\mu_{\rm m}) P_{12}
+(\mu_{\rm c}-\mu_{\rm e}-\mu_{\rm m}) P_{21}
+(\mu_{\rm c}+\mu_{\rm e}) u_{1,2}
+(\mu_{\rm e}-\mu_{\rm c}) u_{2,1}
&=M_{12} \, ,
\\*[5pt]
\mu_{\rm M}\, L_{\rm c}^2\, \widetilde{a} \left(P_{21,22}-P_{22,12}\right)
\quad
&
\\*
+(\mu_{\rm c}-\mu_{\rm e}-\mu_{\rm m}) P_{12}
-(\mu_{\rm c}+\mu_{\rm e}+\mu_{\rm m}) P_{21}
+(\mu_{\rm e}-\mu_{\rm c}) u_{1,2}
+(\mu_{\rm c}+\mu_{\rm e}) u_{2,1}
&=M_{21} \, ,
\\*[5pt]
-\mu_{\rm M}\, L_{\rm c}^2\, \widetilde{a} \left(P_{21,12}-P_{22,11}\right)
\quad
&
\\*
-P_{22} (\lambda_{\rm e}+\lambda_{\rm m}+2 (\mu_{\rm e}+\mu_{\rm m}))
-(\lambda_{\rm e}+\lambda_{\rm m}) P_{11}
+(\lambda_{\rm e}+2 \mu_{\rm e}) u_{2,2}
+\lambda_{\rm e} u_{1,1}
&=M_{22} \, .
\end{split}
\end{align}
which, according to \eqref{constEq} or \eqref{2Dconst}, are accompanied by the following constitutive plane strain equations
\begin{align}\begin{split}\label{const}
\sigma_{11}&=(\lambda_{\rm e}+2\mu_{\rm e})u_{1,1}+\lambda_{\rm e} u_{2,2}-(\lambda_{\rm e}+2\mu_{\rm e})P_{11}-\lambda_{\rm e} P_{22}\, ,
\\[5pt]
\sigma_{22}&=(\lambda_{\rm e}+2\mu_{\rm e})u_{2,2}+\lambda_{\rm e} u_{1,1}-(\lambda_{\rm e}+2\mu_{\rm e})P_{22}-\lambda_{\rm e} P_{11}\, ,
\\*[5pt]
\sigma_{12}&=(\mu_{\rm e}+\mu_{\rm c})u_{1,2}+(\mu_{\rm e}-\mu_{\rm c})u_{2,1}-(\mu_{\rm e}+\mu_{\rm c})P_{12}-(\mu_{\rm e}-\mu_{\rm c})P_{21}\, ,
\\*[5pt]
\sigma_{21}&=(\mu_{\rm e}+\mu_{\rm c})u_{2,1}+(\mu_{\rm e}-\mu_{\rm c})u_{1,2}-(\mu_{\rm e}+\mu_{\rm c})P_{21}-(\mu_{\rm e}-\mu_{\rm c})P_{12}\, ,
\\*[5pt]
m_{13}&=-\mu_{\rm M}\, L_{\rm c}^2\, \widetilde{a}\, (P_{11,2}-P_{12,1})\, ,
\\*[5pt]
m_{23}&=-\mu_{\rm M}\, L_{\rm c}^2\, \widetilde{a}\,(P_{21,2}-P_{22,1}) \, .
\end{split}
\end{align}
Note that according to Eqs \eqref{constEq}, the out-of-plane stress $\sigma_{33}$ is given as $\sigma_{33}=\tfrac{\lambda_{\rm e}}{2(\lambda_{\rm e}+\mu_{\rm e})}(\sigma_{11}+\sigma_{22})$, and the out-of-plane moment stresses $m_{31}$ and $m_{32}$ are given as: $m_{31}=\tfrac{(a_1-a_2)}{(a_1+a_2)} m_{13}$ and $m_{32}=\tfrac{(a_1-a_2)}{(a_1+a_2)} m_{23}$. These components  however do not affect the plane strain equilibrium equations \eqref{pstrain}.
%
%
%

\section{Fundamental solutions for the relaxed micromorphic continuum under plane strain conditions}

\subsection{Concentrated force: The Kelvin problem}

The  Kelvin problem \cite{Kelvin} provides the solution of a point  force acting in the interior of an infinite elastic medium \cite{timoshenko1970}. The solution is of fundamental importance since it provides the plane strain Green's function for the relaxed micromorphic theory. Lord Kelvin (William Thompson, 1824-1907) solved the problem for classical isotropic linear elasticity that was later named after him in 1848.

We consider a body occupying the full plane ($-\infty<x_1<\infty,\ -\infty<x_2<\infty$) under plane-strain conditions. The body is acted upon by a concentrated line force situated at the origin of the coordinate system. There is no loss of generality if we assume that the direction of the line force coincides with the
$x_2$-axis of the coordinate system due to isotropy. In this case, we have that 
\begin{align}
f=\left(
\begin{array}{c}
0 \\
-1
\end{array}
\right) \delta(x_1)\delta(x_2)\,, 
\qquad M=
\left(
\begin{array}{cc}
0 & 0\\
0 & 0
\end{array}
\right),  
\end{align}
with $\delta(x)$ being the Dirac delta function. 

For the solution of the problem the 2D Fourier transform will be employed.  The direct (FT) and inverse (FT$^{-1}$) double Fourier transforms are defined, respectively, as
\begin{align}
\widehat{y}(\xi)=\text{FT}\{y(x)\}=\int_{x\in\mathbb{R}^2} y(x)\, e^{i \, \langle x,\xi \rangle}\, \text{d}x,
\qquad
y(x)=\text{FT}^{-1}\{\widehat{y}(\xi)\}=\frac{1}{(2\pi)^2} \int_{\xi\in\mathbb{R}^2} \widehat{y}(\xi)\,e^{-i \, \langle x,\xi \rangle}\, \text{d}\xi,
\end{align}
where  $\xi=(\xi_1,\xi_2)$ is the 2D Fourier vector with $\| \xi \| \equiv \upxi=\sqrt{\xi_1^2+\xi_2^2}$ and $i$ is the imaginary unit \cite{debnath2014integral}. 
Applying the Fourier transform on the equilibrium equations \eqref{pstrain} and noting that $\text{FT}\{\delta(x_1)\delta(x_2)\}=1$, yields
\begin{equation}  
 \begin{split}\label{FTPL23}
-\left((\lambda_{\rm e} +2\mu _e) \xi_1^2+(\mu_{\rm c}+\mu_{\rm e})\xi_2^2 \right)\widehat{u}_1 - \left(\lambda_{\rm e}+\mu_{\rm e}-\mu_{\rm c}\right)\xi_1 \xi_2 \, \widehat{u}_2+i (\lambda_{\rm e}+2\mu_{\rm e})\xi_1 \widehat{P}_{11}
\\*[3pt]
+i \lambda_{\rm e} \xi_1 \widehat{P}_{22}+i\left(\mu_{\rm c}+\mu_{\rm e}\right) \xi_2 \widehat{P}_{12}
-i \left(\mu_{\rm c}-\mu_{\rm e}\right) \xi_2 \widehat{P}_{21}
&=0
\, ,
\\[5pt]
-\left(\lambda_{\rm e}+\mu_{\rm e}-\mu_{\rm c}\right) \xi_1 \xi_2  \widehat{u}_1-\left((\lambda_{\rm e} +2\mu _e) \xi_2^2+(\mu_{\rm c}+\mu_{\rm e})\xi_1^2 \right)\widehat{u}_2+i \lambda_{\rm e} \xi_2 \widehat{P}_{11}
\\*[3pt]
+i(\lambda_{\rm e}+2\mu_{\rm e})\xi_2 \widehat{P}_{22}+i(\mu_{\rm e}-\mu_{\rm c})\xi_1 \widehat{P}_{12}+i \left(\mu_{\rm c}+\mu_{\rm e}\right)\xi_1 \widehat{P}_{21}
&=-1
\, ,
\\[5pt]
-i(\lambda_{\rm e} +2\mu_{\rm e}) \xi_1 \widehat{u}_1-i \lambda_{\rm e} \xi_2 \widehat{u}_2-\left(\lambda_{\rm e}+2\mu_{\rm e}+\lambda_{\rm m}+2\mu_{\rm m}+\widetilde{a} \mu_{\rm M}\, L_{\rm c}^2 \xi_2^2\right) \widehat{P}_{11}
\\*[3pt]
-(\lambda_{\rm e}+\lambda_{\rm m})\widehat{P}_{22}+\widetilde{a} \mu_{\rm M}\, L_{\rm c}^2 \xi_1 \xi_2 \widehat{P}_{12}
&=0
\, ,
\\[5pt]
-i \lambda_{\rm e} \xi_1 \widehat{u}_1-i(\lambda_{\rm e} +2\mu_{\rm e}) \xi_2 \widehat{u}_2-\left(\lambda_{\rm e}+2\mu_{\rm e}+\lambda_{\rm m}+2\mu_{\rm m}+\widetilde{a} \mu_{\rm M}\, L_{\rm c}^2 \xi_1^2\right) \widehat{P}_{22}
\\*[3pt]
-(\lambda_{\rm e}+\lambda_{\rm m})\widehat{P}_{11}+\widetilde{a} \mu_{\rm M}\, L_{\rm c}^2 \xi_1 \xi_2 \widehat{P}_{21}
&=0
\, ,
\\[5pt]
-i \left(\mu_{\rm c}+\mu_{\rm e}\right)\xi_2 \widehat{u}_1-i \xi_1 \left(\mu_{\rm e}-\mu_{\rm c}\right) \widehat{u}_2+\widetilde{a} \mu_{\rm M}\, L_{\rm c}^2 \xi_1 \xi_2 \widehat{P}_{11}
\\*[3pt]
-\left((\mu_{\rm c}+\mu_{\rm e}+\mu_{\rm m})+\widetilde{a} \mu_{\rm M}\, L_{\rm c}^2 \xi_1^2 \right) \widehat{P}_{12}-(\mu_{\rm e}+\mu_{\rm m}-\mu_{\rm c}) \widehat{P}_{21}
&=0
\, ,
\\[5pt]
-i \left(\mu_{\rm e}-\mu_{\rm c}\right) \xi_2 \widehat{u}_1-i \left(\mu_{\rm e}+\mu_{\rm c}\right) \xi_1 \widehat{u}_2+\widetilde{a}  \mu_{\rm M}\, L_{\rm c}^2 \xi_1\xi_2 \widehat{P}_{22}
\\*[3pt]
-\left((\mu_{\rm c}+\mu_{\rm e}+\mu_{\rm m})+\widetilde{a} \mu_{\rm M}\, L_{\rm c}^2 \xi_2^2 \right) \widehat{P}_{21}-(\mu_{\rm e}+\mu_{\rm m}-\mu_{\rm c}) \widehat{P}_{12}
&=0,
 \end{split}
\end{equation}
where we recall  that $\widetilde{a}\coloneqq(a_1+a_2)/2>0$. 
The algebraic system can be written in the following form
\begin{equation}\label{eq:diesesSystem}
\mathbb{A}({\bf \xi})\,\widehat{u}= \widehat{v},
\end{equation}
where $\widehat{u}=\{\widehat{u}_1,\widehat{u}_2,\widehat{P}_{11},\widehat{P}_{12},\widehat{P}_{21},\widehat{P}_{22}\}^T$, $\widehat{v}=\{0,-1,0,0,0,0\}^T$, and the symmetric Fourier matrix $\mathbb{A}$ is given as
\begin{align}
\footnotesize
&\mathbb{A}(\xi)=\\
\notag
&\hspace{-2.5em}\resizebox{1.15\textwidth}{!}{
$
\left(
\begin{array}{cccccc}
-\xi_2^2 \left(\mu_c+\mu_e\right)-\left(\xi_1^2 \left(\lambda_e+2 \mu_e\right)\right) & -\left(\xi_1 \xi_2
\left(-\mu_c+\lambda_e+\mu_e\right)\right) & i \xi_1 \left(\lambda_e+2 \mu_e\right) & i \xi_2 \left(\mu_c+\mu_e\right) & -i \xi_2 \left(\mu_c-\mu_e\right) & i \xi_1 \lambda_e \\
-\left(\xi_1 \xi_2 \left(-\mu_c+\lambda_e+\mu_e\right)\right) & -\left(\xi_1^2 \left(\mu_c+\mu_e\right)\right)-\xi_2^2 \left(\lambda_e+2 \mu_e\right) & i \xi_2 \lambda_e & -i \xi_1 \left(\mu_c-\mu_e\right)
& i \xi_1 \left(\mu_c+\mu_e\right) & i \xi_2 \left(\lambda_e+2 \mu_e\right) \\
i \xi_1 \left(\lambda_e+2 \mu_e\right) & i \xi_2 \lambda_e &  \widetilde{a} \mu_{\rm M}\, L_{\rm c}^2 \xi_2^2+\lambda_e+2 \left(\mu_e+\mu_m\right)+\lambda_m & - \widetilde{a} \mu_{\rm M}\, L_{\rm c}^2 \xi_1 \xi_2 & 0 & \lambda_e+\lambda_m \\
i \xi_2 \left(\mu_c+\mu_e\right) & i \xi_1 \left(\mu_e-\mu_c\right) & -\widetilde{a} \mu_{\rm M}\, L_{\rm c}^2 \xi_1 \xi_2 & \widetilde{a} \mu_{\rm M}\, L_{\rm c}^2 \xi_1^2+\mu_c+\mu_e+\mu_m & -\mu_c+\mu_e+\mu_m & 0 \\
-i \xi_2 \left(\mu_c-\mu_e\right) & i \xi_1 \left(\mu_c+\mu_e\right) & 0 & -\mu_c+\mu_e+\mu_m & \widetilde{a} \mu_{\rm M}\, L_{\rm c}^2 \xi_2^2+\mu_c+\mu_e+\mu_m & -\widetilde{a} \mu_{\rm M}\, L_{\rm c}^2 \xi_1 \xi_2 \\
i \xi_1 \lambda_e & i \xi_2 \left(\lambda_e+2 \mu_e\right) & \lambda_e+\lambda_m & 0 & -\widetilde{a} \mu_{\rm M}\, L_{\rm c}^2 \xi_1 \xi_2 & \widetilde{a} \mu_{\rm M}\, L_{\rm c}^2 \xi_1^2+\lambda_e+2 \left(\mu_e+\mu_m\right)+\lambda_m \\
\end{array}
\right)
$}
\end{align}
The determinant of the Fourier matrix $\mathbb{A}(\xi)$ becomes
\begin{equation}\label{eq:determinante}
\det \mathbb{A}(\xi)=\begin{cases}\widetilde{a}^2 L_{\rm c}^4\,\mu_{\rm M}^2\, \mu_{\rm m}\, (\mu_{\rm e}+\mu_{\rm c})(\lambda_{\rm e}+2\mu_{\rm e})(\lambda_{\rm m}+2\mu_{\rm m})(\ell_1^{-2}+\upxi^2)(\ell_2^{-2}+\upxi^2)\upxi^4, & \quad\mu_{\rm c}>0,\\
                      \widetilde{a}^2 L_{\rm c}^4\,\mu_{\rm M}^2\, \mu_{\rm m}\, \mu_{\rm e}(\lambda_{\rm e}+2\mu_{\rm e})(\lambda_{\rm m}+2\mu_{\rm m})(\ell_1^{-2}+\upxi^2)\upxi^6, & \quad\mu_{\rm c}=0,
                     \end{cases}
\end{equation}
where $\ell_1$ and $\ell_2$ are two characteristic lengths related with the internal length $L_{\rm c}$ as 
\begin{equation}
\label{chlengths}
\ell_1=L_{\rm c} \sqrt{\frac{{\widetilde{a}\,\beta}\mu_{\rm M}}{4(\kappa_{\rm M}+\mu_{\rm M})}},
\qquad\qquad
\ell_2=L_{\rm c}\sqrt{\frac{{\widetilde{a}\,\mu_{\rm M} (\mu_{\rm e}+\mu_{\rm c})}}{4\mu_{\rm c}\, \mu_{\rm e}}}.
\end{equation}
We recall also that the macroscopic  moduli ($\lambda_{\rm M}, \mu_{\rm M}, \kappa_{\rm M}$)  are related to microscopic-moduli of the relaxed micromorphic medium through equations \eqref{eq:KAPPA} and \eqref{defmod}. Further, the dimensionless parameter $\beta$ is defined as
\begin{equation}
\beta:=\frac{(\kappa_e+\mu_{\rm e})(\kappa_{\rm m}+\mu_{\rm m})}{(\kappa_e+\kappa_{\rm m})(\mu_{\rm e}+\mu_{\rm m})}>0.
\end{equation}
It is interesting to note that $\det \mathbb{A}(\xi)$ is an 8$^{th}$-order polynomial of $\upxi$ with corresponding terms $\{\upxi^8,\upxi^6,\upxi^4\}$, whereas in classical isotropic linear elasticity the Fourier determinant assumes the form
\begin{equation}
\det \mathbb{A}_{\rm lin.elast}(\xi)=\mu_{\rm M}(\lambda_{\rm M}+2\mu_{\rm M})\upxi^4.
\end{equation}

The positive definiteness conditions read simply
\begin{equation}
\label{PD}
\mu_{\rm m}>0, \qquad \mu_{\rm c} \ge 0, \qquad \mu_{\rm e}>0, \qquad \kappa_{\rm m}>0, \qquad \kappa_{\rm e}>0,\qquad \widetilde{a}\,L_{\rm c}^2>0,
\end{equation}
which according to \eqref{defmod}, imply that $\mu_{\rm e}>\mu_{\rm M}>0$ and $\kappa_{\rm e}>\kappa_{\rm M}>0$.

From \eqref{eq:determinante}, the plane(-strain) ellipticity conditions can be readily obtained as (cf. \cite{NeffRealWave})
\begin{align}
 \mu_{\rm M}>0, \quad \mu_{\rm m}>0,\quad \mu_{\rm e} + \mu_{\rm c}>0, \quad \mu_{\rm c}\geq0,  \quad 2\mu_{\rm e}+\lambda_{\rm e}>0, \quad 2\mu_{\rm m} + \lambda_{\rm m}>0,\quad \widetilde{a}\,L_{\rm c}^2>0\,.\end{align}
From the solution of the above non-homogeneous system \eqref{eq:diesesSystem} we derive the solutions for the transformed field variables. These can be written in the following form which is amenable for analytical treatment:
\begin{align}
\begin{split}
\widehat{u}_1
&=
-\frac{\kappa_{\rm M}}{\mu_{\rm M}(\kappa_{\rm M}+\mu_{\rm M})} \frac{\xi_1 \xi_2 }{\upxi^4}-\frac{\widetilde{a} \mu_{\rm M}\, L_{\rm c}^2}{4} \left(\frac{\zeta}{\kappa_{\rm M}+\mu_{\rm M}}\right)^2 \xi_1 \xi_2 \, \phi_1(\upxi)+\frac{\widetilde{a} \mu_{\rm M}\, L_{\rm c}^2}{4\mu_{\rm e}^2}\xi_1 \xi_2 \, \phi_2(\upxi) \, ,
\\*
\widehat{u}_2
&=
\frac{1}{\mu_{\rm M} \xi^2}-\frac{\kappa_{\rm M}}{\mu_{\rm M}(\kappa_{\rm M}+\mu_{\rm M})}\frac{\xi_2^2}{\upxi^4}-\frac{\widetilde{a} \mu_{\rm M}\, L_{\rm c}^2}{4} \left(\frac{\zeta}{\kappa_{\rm M}+\mu_{\rm M}}\right)^2 \xi_2^2 \, \phi_1(\upxi)-\frac{\widetilde{a} \mu_{\rm M}\, L_{\rm c}^2}{4\mu_{\rm e}^2}\xi_1^2 \, \phi_2(\upxi)
\, ,\\
\widehat{P}_{11}
&=
i \frac{\kappa_{\rm M}}{\mu_{\rm m}(\kappa_{\rm M}+\mu_{\rm M})}\frac{\xi_1^2 \xi_2}{\upxi^4}+\frac{i \zeta \xi_2 \left(\varepsilon \widetilde{a} \mu_{\rm M}\, L_{\rm c}^2 \xi_1^2 +2 (\kappa_{\rm m}+\mu_{\rm m}) \right) }{4 (\kappa_{\rm M}+\mu_{\rm M})(\kappa_{\rm m}+\mu_{\rm m})} \phi_1(\upxi)
\, ,
\\*
\widehat{P}_{12}
&=
i \frac{\kappa_{\rm M}}{\mu_{\rm m}(\kappa_{\rm M}+\mu_{\rm M}) }\frac{\xi_1 \xi_2^2}{\upxi^4}+\frac{i \zeta \varepsilon \widetilde{a} \mu_{\rm M}\, L_{\rm c}^2 \, \xi_1 \xi_2^2}{4 (\kappa_{\rm M}+\mu_{\rm M})(\kappa_{\rm m}+\mu_{\rm m})} \phi_1(\upxi)+\frac{i \xi_1}{2 \mu_{\rm e}} \phi_2(\upxi)
\, ,
\\*
\widehat{P}_{21}
&=
-\frac{i \xi_1 \left(\left(\kappa_{\rm M}+\mu_{\rm M}\right)\xi_1^2+\mu_{\rm M} \xi_2^2 \right)}{
\mu_{\rm m} \left(\kappa_{\rm M}+\mu_{\rm M}\right) \upxi^4}+\frac{i \zeta \varepsilon \widetilde{a} \mu_{\rm M}\, L_{\rm c}^2\, \xi_1 \xi_2^2}{4 (\kappa_{\rm M}+\mu_{\rm M})(\kappa_{\rm m}+\mu_{\rm m})} \phi_1(\upxi)-\frac{i \xi_1}{2 \mu_{\rm e}} \phi_2(\upxi)
\, ,
\\*
\widehat{P}_{22}
&=
-i \frac{\kappa_{\rm M}}{\mu_{\rm m}(\kappa_{\rm M}+\mu_{\rm M}) }\frac{\xi_1^2 \xi_2}{\upxi^4}-\frac{i \xi_2}{\left(\kappa_{\rm m}+\mu_{\rm m}\right)\upxi^2}-\frac{i \zeta \xi_2 \left(\varepsilon \widetilde{a} \mu_{\rm M}\, L_{\rm c}^2 \xi_1^2 +2 (\kappa_{\rm m}-\mu_{\rm m}) \right)}{4 (\kappa_{\rm M}+\mu_{\rm M})(\kappa_{\rm m}+\mu_{\rm m})} \phi_1(\upxi)
\, ,
\end{split}
\end{align}
where the transformed functions $\phi_j(\upxi)$ ($j=1,2$) and dimensionless parameters ($\zeta$, $\varepsilon$) are defined as
\begin{equation}
\phi_j(\upxi)=\frac{1}{\upxi^2}-\frac{\ell_j^{2}}{1+\ell_j^{2}\upxi^2}
\, ,
\qquad\qquad
\zeta=\frac{\mu_{\rm M}}{\mu_{\rm m}}-\frac{\kappa_M}{\kappa_{\rm m}}
\, ,
\qquad\qquad
\varepsilon=\frac{\kappa_{\rm m}}{\kappa_{\rm M}+\mu_{\rm M}}\, \beta.
\end{equation}
We employ now some useful classical results (see e.g. \cite{gradshteyn2014table}, \cite{Now74}):
\begin{align}
\label{FPs}
I_1&=\text{FT}^{-1}\{\left(\xi _1^2+\xi _2^2\right)^{-1}\}=-\frac{1}{2\pi} (b+\ln r) \, ,
\notag
\\
I_2&=\text{FT}^{-1}\{\left(\xi _1^2+\xi _2^2\right)^{-2}\}=\frac{1}{8\pi} r^2(b+\ln r) \, ,
\\
I_3&=\text{FT}^{-1}\{\left(\ell^{-2}+\xi _1^2+\xi _2^2\right)^{-1}\}=\frac{1}{2\pi} K_0\left[\frac{r}{\ell}\right] \, ,
\notag
\end{align}
and
\begin{equation}
\partial^m_{x_1}\partial^n_{x_2}I_j=(-i \xi_1)^m (- i\xi_2)^n I_j, \qquad (m,n=0,1,2,...), \enskip (j=1,2,3)
\end{equation}
where $r=\sqrt{x_1^2+x_2^2}$, $K_n[\cdot]$ is the $n$-th order second kind modified Bessel functions and $b=0.57...$ is Euler's constant \cite{gradshteyn2014table}. It should be noted that the first two integrals in \eqref{FPs} are defined as the finite part integrals\footnote{The concept of a finite-part integral has been first introduced by Hadamard \cite{had23} in 1923. These integrals have stronger singularities than Cauchy principal value integrals and they exist in the finite part sense \cite{kutt1975,monegato2009}.}.

Using the above results, the definitions of the characteristic lengths $\ell_1$, $\ell_2$ and ignoring rigid body motions in the displacement field, we obtain after some rather extensive algebra the following expressions for the displacement and micro-distortion fields

\begin{align}\label{Green}
\begin{split}
u_1&=
\frac{\kappa_{\rm M}\, x_1 x_2}{4 \pi\, \mu_{\rm M} \left(\kappa_{\rm M}+\mu_{\rm M}\right)  r^2}+\frac{\zeta^2 x_1 x_2}{2 \pi  \beta  \left(\kappa_{\rm M}+\mu_{\rm M}\right)r^2}\Phi_1-\frac{\mu_{\rm c} x_1 x_2}{2 \pi\, \mu_{\rm e} \left(\mu_{\rm c}+\mu_{\rm e}\right)r^2} \Phi_2
\, ,
\\[5pt]
u_2&=
\frac{\kappa_{\rm M} x_2^2}{4 \pi  \mu_{\rm M} \left(\kappa_{\rm M}+\mu_{\rm M}\right)r^2}
-\frac{(\kappa_{\rm M}+2 \mu_{\rm M})}{4\pi\,\mu_{\rm M} (\kappa_{\rm M}+\mu_{\rm M})}\ln r
-\frac{\zeta^2}{4 \pi  \beta  \left(\kappa_{\rm M}+ \mu_{\rm M}\right)}\left(\frac{(x_1^2-x_2^2)}{r^2} \Phi_1+K_0\left[\frac{r}{\ell_1}\right]\right)
\\*[3pt]
&\phantom{=}
+\frac{\mu_{\rm c}}{4 \pi\, \mu_{\rm e} \left(\mu_{\rm c}+\mu_{\rm e}\right)}\left(\frac{(x_1^2-x_2^2)}{r^2} \Phi_2-K_0\left[\frac{r}{\ell_2}\right]\right)
\, ,
\\[5pt]
P_{11}&=-
\frac{\kappa_{\rm M} x_2 \left(x_1^2-x_2^2\right)}{4 \pi\, \mu_{\rm m} (\kappa_{\rm M}+\mu_{\rm M}) r^4}
+\frac{\zeta  x_2}{4 \pi \left(\kappa_{\rm M}+\mu_{\rm M}\right) r}\Psi_1
-\frac{\zeta \varepsilon x_2 \left(x_1^2-x_2^2\right)}{2 \pi \left(\kappa_{\rm m}+\mu_{\rm m}\right)\beta r^4} \Phi_1
+\frac{\zeta \varepsilon x_1 x_2}{2 \pi \left(\kappa_{\rm m}+\mu_{\rm m}\right)\beta r^2}\, \partial_{x_1}\Phi_1
\, ,
\\[5pt]
P_{12}&=
\frac{\kappa_{\rm M} x_1 \left(x_1^2-x_2^2\right)}{4 \pi\, \mu_{\rm m} \left(\kappa_{\rm M}+\mu_{\rm M}\right)r^4}
+\frac{x_1}{4 \pi\, \mu_{\rm e} r}\Psi_2+\frac{\zeta \varepsilon x_1 \left(x_1^2-x_2^2\right)}{2 \pi \left(\kappa_{\rm m}+\mu_{\rm m}\right)\beta r^4} \Phi_1
+\frac{\zeta \varepsilon x_1 x_2}{2 \pi \left(\kappa_{\rm m}+\mu_{\rm m}\right)\beta r^2}\, \partial_{x_2} \Phi_1
\, ,
\\[5pt]
P_{21}&=
\frac{\kappa_{\rm M} x_1 \left(x_1^2-x_2^2\right)}{4 \pi\, \mu_{\rm m} \left(\kappa_{\rm M}+\mu_{\rm M}\right)r^4}
-\frac{x_1}{4 \pi\, \mu_{\rm e} r}\Psi_2+\frac{\zeta \varepsilon x_1 \left(x_1^2-x_2^2\right)}{2 \pi \left(\kappa_{\rm m}+\mu_{\rm m}\right)\beta r^4} \Phi_1+\frac{\zeta \varepsilon x_1 x_2}{2 \pi \left(\kappa_{\rm m}+\mu_{\rm m}\right)\beta r^2}\, \partial_{x_2} \Phi_1
\\*[3pt]
&\phantom{={}}
-\frac{x_1}{2\pi\, \mu_{\rm m} r^2}
\, ,
\\[5pt]
P_{22}&=
-\frac{x_2}{2\pi (\kappa_{\rm m}+\mu_{\rm m})r^2}
+\frac{\kappa_{\rm M} x_2 \left(x_1^2-x_2^2\right)}{4 \pi\, \mu_{\rm m} \left(\kappa_{\rm M}+\mu_{\rm M}\right)r^4}-\frac{\zeta (\kappa_{\rm m}-\mu_{\rm m}) x_2}{4 \pi \left(\kappa_{\rm M}+\mu_{\rm M}\right) (\kappa_{\rm m}+\mu_{\rm m}) r}\Psi_1
\\*[3pt]
&\phantom{={}}
+\frac{\zeta \varepsilon x_2 \left(x_1^2-x_2^2\right)}{2 \pi \left(\kappa_{\rm m}+\mu_{\rm m}\right)\beta r^4} \Phi_1
-\frac{\zeta \varepsilon x_1 x_2}{2 \pi \left(\kappa_{\rm m}+\mu_{\rm m}\right)\beta r^2}\, \partial_{x_1}\Phi_1,
\end{split}
\end{align}
where the functions $\Phi_j$ and $\Psi_j$ ($j=1,2$) are defined as
\begin{align}
\Phi_j \equiv \Phi_j(r)=\frac{2 \ell_j^2}{r^2}-K_2\left[\frac{r}{\ell_j}\right] \, ,
\qquad\qquad\qquad
\Psi_j \equiv \Psi_j(r)=\frac{1}{r}\left(1-\frac{r}{\ell_j}K_1\left[\frac{r}{\ell_j}\right]\right) \, .
\end{align}
Some useful relations and limits for the second kind modified Bessel functions that have been used for the derivation of our equations can be found in Appendix \ref{Bessel}.

Equations \eqref{Green} are the basic results of this paper and constitute the Green's functions for the general relaxed isotropic micromorphic continuum under plane strain conditions for the case of a concentrated force acting in the $x_2$-direction. The Green's functions for the case where the concentrated force acts in the $x_1$-direction can be readily derived from the above solution by interchanging the indices $1 \leftrightarrow 2$.

The micro-rotation for the relaxed micromorphic medium in the case of plane strain is defined as the skew-symmetric part of $P$ (see \eqref{skewP})
\begin{equation}
\label{microrot}
\vartheta_3=\frac{1}{2}(P_{21}-P_{12})=-\frac{x_1}{4\pi r^2}\left(\frac{1}{\mu_{\rm M}}-\frac{1}{\mu_{\rm e}}\frac{r}{\ell_2}K_1\left[\frac{r}{\ell_2}\right]\right).
\end{equation}
Finally, it is noted that the stresses and higher order stresses can be derived from the constitutive relations \eqref{const}.

Using now the asymptotic properties of the second kind modified Bessel functions as $z \to 0$ \cite{gradshteyn2014table}
\begin{equation}
K_n[z] \sim
\begin{cases}
-\ln\tfrac{z}{2}-b, \enskip &\text{for} \enskip n=0\,,\\
\tfrac{\Gamma[n]}{2} \left(\tfrac{2}{z}\right)^n \enskip 
&\text{for} \enskip n>0\,,
\end{cases}    
\end{equation}
we may readily deduce that as $r \to 0$ the displacement field becomes logarithmically unbounded as in the classical linear elastic theory and the micro-distortion field $P$ exhibits an $r^{-1}$ singularity consistent with the additive coupling $\text{D}u-P$.    This in turn implies that, according to \eqref{constEq}$_1$, the components of the stress tensor $\sigma$ behave also as $O(r^{-1})$ as  $r \to 0$. The same singular asymptotic behaviour is exhibited by the micro-rotation $\vartheta_3$. In particular, the second term in \eqref{microrot} is bounded as $r \to 0$ but the first term behaves as $r^{-1}$. It is worth noting that the micro-rotation becomes bounded if $\mu_{\rm m} \to \infty$ which is the case of micro-stretch, micropolar and couple stress elasticity as we shall see next. Interestingly, it turns out that the components of $\text{Curl}P$ have at most a logarithmic singularity which implies, according to \eqref{constEq}$_3$, that the moment stresses $m$ exhibit also a $O(\log r)$ behaviour as  $r \to 0$.

The above results corroborate uniqueness for our solutions. Indeed, for a unique solution of the concentrated load problem the conventional and higher order stress singularities must behave at most as  $O(r^{-1})$ when $r \to 0$, where $r$ is the distance from the point of application of the concentrated loads (see  Hartranft and Shi \cite{hart68} and Sternberg \cite{stern68} for the case of couple stress elasticity and Eubanks and Sternberg \cite{sternberg1955concept} for the classical elasticity case). This is due to the fact that the tractions on a circle surrounding and separating the
concentrated load point from the rest of the medium must be statically equivalent to the concentrated force at that point. This is a general requirement and is independent of the elasticity theory that is employed.

\subsubsection{The relaxed micromorphic continuum with zero micro and macro Poisson's ratio}
A simpler case arises for zero micro and macro Poisson's ratio so that $\lambda_{\rm e}=\lambda_{\rm m}=0$ which implies $\lambda_{\rm M}=0$ and $\zeta=0$. In this case, we derive

\begin{align}
\label{zerop}
u_1&=\frac{x_1 x_2}{8 \pi r^2}\left(\frac{1}{\mu_{\rm M}}-\frac{4\mu_{\rm c}}{\mu_{\rm e} (\mu_{\rm e}+\mu_{\rm c})}\Phi_2\right) \, ,
\notag
\\*
u_2&=\frac{x_2^2}{8 \pi\, \mu_{\rm M}\, r^2}-\frac{3}{8 \pi\, \mu_{\rm M}} \ln r+\frac{\mu_{\rm c}}{4 \pi\, \mu_{\rm e} \left(\mu_{\rm c}+\mu_{\rm e}\right)}\left(\frac{(x_1^2-x_2^2)}{r^2} \Phi_2-K_0\left[\frac{r}{\ell_2}\right]\right) \, ,
\\*
P_{11}&=-\frac{x_2 \left(x_1^2-x_2^2\right)}{8 \pi\, \mu_{\rm m}  r^4} \, ,
\qquad\qquad \,\,\,
P_{12}=\frac{x_1 \left(x_1^2-x_2^2\right)}{8 \pi\, \mu_{\rm m} r^4}+\frac{x_1}{4 \pi\, \mu_{\rm e} r} \, \Psi_2\, ,
\notag
\\*
P_{22}&=-\frac{x_2 \left(x_1^2+3 x_2^2\right)}{8 \pi\, \mu_{\rm m} r^4} \, ,
\qquad\qquad
P_{21}=-\frac{x_1 \left(3x_1^2+5x_2^2\right)}{8 \pi\, \mu_{\rm m} r^4}-\frac{x_1}{4 \pi\, \mu_{\rm e} r} \Psi_2 \,.
\notag
\end{align}
It is evident that $u_2$ retains the logarithmic singularity but the detailed field is different, in particular
\begin{align}
u_2=-\frac{\left(3 \mu_{\rm c} \mu_{\rm e}+\mu_{\rm c} \mu_{\rm m}+3 \mu_{\rm e}^2+3 \mu_{\rm e} \mu_{\rm m}\right)}{8 \pi  \mu_{\rm e} \mu_{\rm m} \left(\mu_{\rm c}+\mu_{\rm e}\right)} \ln r, \qquad \text{as} \quad r \to 0.
\end{align}

\subsubsection{The pure relaxed micromorphic continuum with symmetric force stress tensor}
Another special case of interest is the pure relaxed micromorphic continuum with symmetric force stress tensor $\sigma$. In this case we have that the Cosserat modulus $\mu_{\rm c} = 0$ (which implies that $\ell_2 \to \infty$) and accordingly (see Appendix \ref{Bessel})
\begin{align}
\label{limitpure}
\lim_{\mu_{\rm c} \to 0} \mu_{\rm c}\,\Phi_2=0
\, ,
\qquad\qquad\qquad
\lim_{\mu_{\rm c} \to 0}   \Psi_2=0
\, ,
\qquad\qquad\qquad
\lim_{\mu_{\rm c} \to 0} \mu_{\rm c}\, K_0\left[\frac{r}{\ell_2}\right]=0
\, ,
\end{align}
and we derive
\begin{align}\begin{split}
u_1&=
\frac{\kappa_{\rm M}\, x_1 x_2}{4 \pi\, \mu_{\rm M} \left(\kappa_{\rm M}+\mu_{\rm M}\right)  r^2}+\frac{\zeta^2 x_1 x_2}{2 \pi  \beta  \left(\kappa_{\rm M}+\mu_{\rm M}\right)r^2}\Phi_1
\, ,
\\[5pt]
u_2&=
\frac{\kappa_{\rm M} x_2^2}{4 \pi  \mu_{\rm M} \left(\kappa_{\rm M}+\mu_{\rm M}\right)r^2}
-\frac{(\kappa_{\rm M}+2 \mu_{\rm M})}{4\pi\,\mu_{\rm M} (\kappa_{\rm M}+\mu_{\rm M})}\ln r
-\frac{\zeta^2}{4 \pi  \beta  \left(\kappa_{\rm M}+ \mu_{\rm M}\right)}\left(\frac{(x_1^2-x_2^2)}{r^2} \Phi_1+K_0\left[\frac{r}{\ell_1}\right]\right)
\, ,
\\[5pt]
P_{11}&=
-\frac{\kappa_{\rm M} x_2 \left(x_1^2-x_2^2\right)}{4 \pi\, \mu_{\rm m} (\kappa_{\rm M}+\mu_{\rm M}) r^4}
+\frac{\zeta  x_2}{4 \pi \left(\kappa_{\rm M}+\mu_{\rm M}\right) r}\Psi_1
-\frac{\zeta \varepsilon x_2 \left(x_1^2-x_2^2\right)}{2 \pi \left(\kappa_{\rm m}+\mu_{\rm m}\right)\beta r^4} \Phi_1
+\frac{\zeta \varepsilon x_1 x_2}{2 \pi \left(\kappa_{\rm m}+\mu_{\rm m}\right)\beta r^2}\, \partial_{x_1}\Phi_1
\, ,
\\[5pt]
P_{12}&=
\frac{\kappa_{\rm M} x_1 \left(x_1^2-x_2^2\right)}{4 \pi\, \mu_{\rm m} \left(\kappa_{\rm M}+\mu_{\rm M}\right)r^4}
+\frac{\zeta \varepsilon x_1 \left(x_1^2-x_2^2\right)}{2 \pi \left(\kappa_{\rm m}+\mu_{\rm m}\right)\beta r^4} \Phi_1
+\frac{\zeta \varepsilon x_1 x_2}{2 \pi \left(\kappa_{\rm m}+\mu_{\rm m}\right)\beta r^2}\, \partial_{x_2} \Phi_1
\, ,
\\[5pt]
P_{21}&=
\frac{\kappa_{\rm M} x_1 \left(x_1^2-x_2^2\right)}{4 \pi\, \mu_{\rm m} \left(\kappa_{\rm M}+\mu_{\rm M}\right)r^4}
+\frac{\zeta \varepsilon x_1 \left(x_1^2-x_2^2\right)}{2 \pi \left(\kappa_{\rm m}+\mu_{\rm m}\right)\beta r^4} \Phi_1+\frac{\zeta \varepsilon x_1 x_2}{2 \pi \left(\kappa_{\rm m}+\mu_{\rm m}\right)\beta r^2}\, \partial_{x_2} \Phi_1
-\frac{x_1}{2\pi\, \mu_{\rm m} r^2}
\, ,
\\[5pt]
P_{22}&=
-\frac{x_2}{2\pi (\kappa_{\rm m}+\mu_{\rm m})r^2}
+\frac{\kappa_{\rm M} x_2 \left(x_1^2-x_2^2\right)}{4 \pi\, \mu_{\rm m} \left(\kappa_{\rm M}+\mu_{\rm M}\right)r^4}-\frac{\zeta (\kappa_{\rm m}-\mu_{\rm m}) x_2}{4 \pi \left(\kappa_{\rm M}+\mu_{\rm M}\right) (\kappa_{\rm m}+\mu_{\rm m}) r}\Psi_1
\\*[3pt]
&\phantom{={}}
+\frac{\zeta \varepsilon x_2 \left(x_1^2-x_2^2\right)}{2 \pi \left(\kappa_{\rm m}+\mu_{\rm m}\right)\beta r^4} \Phi_1
-\frac{\zeta \varepsilon x_1 x_2}{2 \pi \left(\kappa_{\rm m}+\mu_{\rm m}\right)\beta r^2}\, \partial_{x_1}\Phi_1.
\end{split}
\end{align}

\subsubsection{Limiting cases}
It is shown here that the fundamental solutions of several well-known generalized continua can be obtained as singular limiting cases of the general relaxed micromorphic fundamental solution for a concentrated force.

\paragraph{Micro-stretch elasticity}
In order to pass from the general relaxed micromorphic continua to the micro-stretch continua we let $\mu_{\rm m} \to \infty$ which, according to \eqref{defmod}, implies that: $\mu_{\rm e} \to \mu_{\rm M}$, and
\begin{equation}
\label{limstretch}
\lim_{\mu_{\rm m} \to \infty} \zeta=\frac{\kappa_{\rm M}-\kappa_{\rm e}}{\kappa_{\rm e}} \,, \quad \qquad \lim_{\mu_{\rm m} \to \infty} \beta=\frac{(\kappa_{\rm e}-\kappa_{\rm M})(\kappa_{\rm e}+\mu_{\rm M})}{\kappa_{\rm e}^2}.
\end{equation}
In this case, the kinematical fields read
\begin{align}
\begin{split}
u_1=&
\frac{\kappa_{\rm M}\, x_1 x_2}{4 \pi\, \mu_{\rm M} \left(\kappa_{\rm M}+\mu_{\rm M}\right)  r^2}+\frac{\kappa_{\rm e}-\kappa_{\rm M}}{2 \pi (\kappa_{\rm e}+\mu_{\rm M})(\kappa_{\rm M}+\mu_{\rm M})}\frac{x_1 x_2}{r^2}\Phi_1-\frac{\mu_{\rm c}}{2 \pi\, \mu_{\rm M} \left(\mu_{\rm c}+\mu_{\rm M}\right)} \frac{x_1 x_2}{r^2} \Phi_2
\, ,
\\[5pt]
u_2=&
\frac{\kappa_{\rm M} x_2^2}{4 \pi  \mu_{\rm M} \left(\kappa_{\rm M}+\mu_{\rm M}\right)r^2}
-\frac{(\kappa_{\rm M}+2 \mu_{\rm M})}{4\pi\,\mu_{\rm M} (\kappa_{\rm M}+\mu_{\rm M})}\ln r
+\frac{\mu_{\rm c}}{4 \pi\, \mu_{\rm M} \left(\mu_{\rm c}+\mu_{\rm M}\right)}\left(\frac{(x_1^2-x_2^2)}{r^2} \Phi_2-K_0\left[\frac{r}{\ell_2}\right]\right)
\\*[3pt]
&-\frac{\kappa_{\rm e}-\kappa_{\rm M}}{4 \pi (\kappa_{\rm e}+\mu_{\rm M})(\kappa_{\rm M}+\mu_{\rm M})}\left(\frac{(x_1^2-x_2^2)}{r^2} \Phi_1+K_0\left[\frac{r}{\ell_1}\right]\right)
\, ,
\\[5pt]
P_{11}=&P_{22}=-\frac{(\kappa_{\rm e}-\kappa_{\rm M})x_2}{4 \pi \kappa_{\rm e}(\kappa_{\rm M}+\mu_{\rm M})r} \Psi_1 \, , \qquad  P_{12}=-P_{21}=\frac{x_1}{4 \pi\, \mu_{\rm M} r} \Psi_2
\, ,
\end{split}
\end{align}
and the micro-rotation is given as 
\begin{equation}
\label{microrot2}
\vartheta_3=\frac{1}{2}(P_{21}-P_{12})=-\frac{x_1}{4\pi\,\mu_{\rm M}\, r^2}\left(1-\frac{r}{\ell_2}K_1\left[\frac{r}{\ell_2}\right]\right),
\end{equation}
where the characteristic lengths are now defined as
\begin{equation}
\label{l1stretch}
\ell_1=L_{\rm c} \sqrt{\frac{\widetilde{a}\,\mu_{\rm M}(\kappa_{\rm e}+\mu_{\rm M})}{4(\kappa_{\rm M}+\mu_{\rm M})\,(\kappa_{\rm e}+\kappa_{\rm m})}} \, , \quad\qquad \ell_2=L_{\rm c}\sqrt{\frac{{\widetilde{a}\,(\mu_{\rm M}+\mu_{\rm c})}}{4\mu_{\rm c}}}.
\end{equation}
We note again that $\ell_{2}\to\infty$ as $\mu_{\rm c}\to0$.

\paragraph{Cosserat (micropolar) elasticity}
As $(\mu_{\rm m},{\kappa_{\rm m}}) \to \infty$ we have that: $\mu_{\rm e} \to \mu_{\rm M}$, $\kappa_{\rm e} \to \kappa_{\rm M}$, $\lambda_{\rm e} \to \lambda_{\rm M}$, and also $\zeta \to 0$, $\beta \to 0$ which implies further that $\ell_1 \to 0$. Furthermore, by recalling that $\kappa_{\rm M}=\lambda_{\rm M}+\mu_{\rm M}$, and identifying (using Nowacki's notation \cite{Now72}) $\mu_{\rm c}=\alpha$, $ a_1 \mu_{\rm M}\, L_{\rm c}^2= 2\gamma$, $a_2 \mu_{\rm M}\, L_{\rm c}^2=2 \varepsilon$, the relaxed micromorphic solution degenerates to the known micropolar solution (\cite{LiaHuang96}, \cite{Dyszlewicz2004})\footnote{It should be noted that in \cite{Dyszlewicz2004} there is a misprint in the plane strain fundamental solution (3.78). In particular, the term $(1-\nu)$ should be replaced with $(1-\nu)^{-1}$. Also, the solution in \cite{Dyszlewicz2004} is for a horizontal force which can be transformed to the solution for a vertical force solution as in the present case by interchanging the indices $1 \leftrightarrow 2$.}
\begin{align}
\begin{split}
u_1&=
\frac{\left(\lambda_{\rm M}+\mu_{\rm M}\right)}{4 \pi\, \mu_{\rm M} \left(\lambda_{\rm M}+2 \mu_{\rm M}\right)}\frac{x_1 x_2}{r^2}-\frac{\alpha}{2 \pi\, \mu_{\rm M} (\mu_{\rm M}+\alpha)}\frac{x_1 x_2}{r^2}\left(\frac{2\ell^2}{r^2}- K_2\left[\frac{r}{\ell}\right] \right)
\, ,
\\[5pt]
u_2&=
\frac{(\lambda_{\rm M}+\mu_{\rm M})}{4 \pi  \mu_{\rm M} \left(\lambda_{\rm M}+2 \mu_{\rm M}\right)}\frac{x_2^2}{r^2}
-\frac{(\lambda_{\rm M}+3 \mu_{\rm M})}{4\pi\,\mu_{\rm M} (\lambda_{\rm M}+2\mu_{\rm M})}\ln r-\frac{\alpha}{4 \pi\, \mu_{\rm M} \left(\alpha+\mu_{\rm M}\right)} K_0\left[\frac{r}{\ell}\right]
\\*[3pt]
&\phantom{={}}
+\frac{\alpha}{4 \pi\, \mu_{\rm M} \left(\alpha+\mu_{\rm M}\right)}\frac{(x_1^2-x_2^2)}{r^2} \left(\frac{2\ell^2}{r^2}- K_2\left[\frac{r}{\ell}\right] \right)
\, ,
\\[5pt]
P_{11}&=
P_{22}=0, 
\qquad\qquad
A_{12}=P_{12}=-P_{21}=-A_{21}=\frac{x_1}{4 \pi\, \mu_{\rm M}\, r^2}\left(1-\frac{r}{\ell}K_1\left[\frac{r}{\ell}\right]\right)
\end{split}
\end{align}
with the micro-rotation $\vartheta_3$ given as
\begin{equation}
\vartheta_3=\frac{1}{2}(P_{21}-P_{12})=-\frac{x_1}{4 \pi\, \mu_{\rm M}\, r^2}\left(1-\frac{r}{\ell}K_1\left[\frac{r}{\ell}\right]\right)
\end{equation}
where 
\begin{align}
\label{lmp}
\ell \equiv \ell_2=\sqrt{\frac{(\gamma+\varepsilon)(\mu_{\rm M}+\alpha)}{4\, \alpha\, \mu_{\rm M}}} \, =L_{\rm c}\sqrt{\frac{{\widetilde{a}\,(\mu_{\rm M}+\mu_{\rm c})}}{4\mu_{\rm c}}},
\end{align}
is the known characteristic length of the Cosserat (micropolar) theory.
\paragraph{Couple stress elasticity}
As $(\mu_{\rm m},\kappa_{\rm m},\mu_{\rm c}) \to \infty$ we have that: $\mu_{\rm e} \to \mu_{\rm M}$, $\lambda_{\rm e} \to \lambda_{\rm M}$, and also $\zeta \to 0$, $\beta \to 0$ which implies further that $\ell_1 \to 0$. In this case, we pass to Mindlin's \cite{mindlin1962effects} and Koiter's \cite{koiter1964couple} theory of couple stress elasticity (see also \cite{Ghiba2017indeterminate,ANewView,zisis2014some, WhyYang,gourgiotis2019hertz}). It is well-known that the spherical part of the couple stress tensor remains indeterminate within the frame of this theory, as Neff et al. \cite{neff2016some} pointed out this is not inconsistent, indeed, like the pressure in an incompressible body, it is indeterminate in the local constitutive law but can be found a posteriori from the solution of the boundary value problem. For a different take on the couple stress model which allows to determine locally the spherical part see \cite{soldatos2023determination}.

Now, identifying $a_1 \mu_{\rm M}\, L_{\rm c}^2=a_2  \mu_{\rm M}L_{\rm c}^2=4 \eta$, we derive the fundamental solution in couple stress theory \cite{hattori2023isogeometric} which assumes the following form
\begin{align}
\begin{split}
u_1&=
\frac{\left(\lambda_{\rm M}+\mu_{\rm M}\right)}{4 \pi\, \mu_{\rm M} \left(\lambda_{\rm M}+2 \mu_{\rm M}\right)}\frac{x_1 x_2}{r^2}-\frac{1}{2 \pi\, \mu_{\rm M}}\frac{x_1 x_2}{r^2}\left(\frac{2\ell^2}{r^2}- K_2\left[\frac{r}{\ell}\right] \right)
\, ,
\\[5pt]
u_2&=
\frac{(\lambda_{\rm M}+\mu_{\rm M})}{4 \pi  \mu_{\rm M} \left(\lambda_{\rm M}+2 \mu_{\rm M}\right)}\frac{x_2^2}{r^2}-\frac{(\lambda_{\rm M}+3 \mu_{\rm M})}{4\pi\,\mu_{\rm M} (\lambda_{\rm M}+2\mu_{\rm M})}\ln r-\frac{1}{4 \pi\, \mu_{\rm M}} K_0\left[\frac{r}{\ell}\right]
\, 
\\*[3pt]
&\phantom{={}}
+\frac{1}{4 \pi\, \mu_{\rm M}}\frac{(x_1^2-x_2^2)}{r^2} \left(\frac{2\ell^2}{r^2}- K_2\left[\frac{r}{\ell}\right] \right)
\, ,
\\[5pt]
P_{11}&=P_{22}=0,
\qquad\qquad
P_{12}=-P_{21}=\frac{x_1}{4 \pi\, \mu_{\rm M}\, r^2}\left(1-\frac{r}{\ell}K_1\left[\frac{r}{\ell}\right]\right),
\end{split}
\end{align}
where the characteristic length of the couple stress elasticity model is defined as
\begin{align}
\label{lcs}
\ell \equiv \ell_2=\sqrt{\frac{\eta}{\mu_{\rm M}}} =L_{\rm c}\sqrt{\frac{\widetilde{a}}{4}}.
\end{align}
As expected, the continuum-rotation $\overline{\vartheta}_3$ coincides with the skew symmetric part of $P$ (i.e. the micro-rotation $\vartheta_3$). Indeed,
\begin{equation}
\overline{\vartheta}_3=\frac{1}{2}\left( \frac{\partial u_2}{\partial x_1}-\frac{\partial u_1}{\partial x_2} \right)=\frac{1}{2}(P_{21}-P_{12})=-\frac{x_1}{4 \pi\, \mu_{\rm M}\, r^2}\left(1-\frac{r}{\ell}K_1\left[\frac{r}{\ell}\right]\right).
\end{equation}
Fundamental solutions for anisotropic couple stress materials under static and dynamic conditions can be found in \cite{GB16a,GB16b,BG16,GB17}.
\paragraph{Classical linear elasticity  ($L_{\rm c} \to 0$) - lower bound macroscopic stiffness}
As $L_{\rm c} \to 0$ we have also that $\ell_j \to 0$ ($j=1,2$) if $\mu_{\rm c}>0$, and in this case we obtain that (see Appendix \ref{Bessel})
\begin{align}
\label{lmclF}
\lim_{\ell_j \to 0} \Phi_j=0 \, ,
\qquad\quad
\lim_{\ell_j \to 0} \partial_{x_i}\Phi_j=0 \, ,
\quad (i=1,2)\,,
\qquad\quad
\lim_{\ell_j \to 0}   \Psi_j=\frac{1}{r} \, ,
\qquad\quad
\lim_{\ell_j \to 0} K_0\left[\frac{r}{\ell_j}\right]=0 \, .
\end{align}
Moreover, by using $\kappa_{\rm M}=\lambda_{\rm M}+\mu_{\rm M}$, we finally derive
\begin{align}
\label{CLF}
u_1=
\frac{\left(\lambda_{\rm M}+\mu_{\rm M}\right)}{4 \pi\, \mu_{\rm M} \left(\lambda_{\rm M}+2 \mu_{\rm M}\right)} \frac{x_1 x_2}{r^2} \, , \qquad
u_2=
\frac{(\lambda_{\rm M}+\mu_{\rm M})}{4 \pi  \mu_{\rm M} \left(\lambda_{\rm M}+2 \mu_{\rm M}\right)}\frac{x_2^2}{r^2}-\frac{(\lambda_{\rm M}+3 \mu_{\rm M})}{4\pi\, \mu (\lambda_{\rm M}+2\mu_{\rm M})}\ln r \, ,
\end{align}
which is the standard classical linear elasticity fundamental solution for the displacements \cite{timoshenko1970}. Moreover, the continuum rotation is given as 
\begin{equation}
\overline{\vartheta}_3=-\frac{x_1}{4\, \pi\,\mu_{\rm M}\, r^2}.
\end{equation}

In addition,

\begin{align}
\begin{split}
P_{11}&=
\frac{\zeta  x_2}{4 \pi \left(\lambda_{\rm M}+2 \mu_{\rm M}\right) r^2}-\frac{(\lambda_{\rm M}+\mu_{\rm M})x_2 \left(x_1^2-x_2^2\right)}{4 \pi\, \mu_{\rm m} (\lambda_{\rm M}+2\mu_{\rm M}) r^4}
\, ,
\\[5pt]
P_{12}&=
\frac{x_1}{4 \pi\, \mu_{\rm e} r^2}+\frac{(\lambda_{\rm M}+\mu_{\rm M}) x_1 \left(x_1^2-x_2^2\right)}{4 \pi\, \mu_{\rm m} \left(\lambda_{\rm M}+2 \mu_{\rm M}\right)r^4}
\, ,
\\[5pt]
P_{21}&=
-\frac{x_1}{4\pi\,\mu_{\rm e} r^2}-\frac{x_1}{2\pi\, \mu_{\rm m} r^2}+\frac{(\lambda_{\rm M}+\mu_{\rm M}) x_1 \left(x_1^2-x_2^2\right)}{4 \pi\, \mu_{\rm m} \left(\lambda_{\rm M}+2 \mu_{\rm M}\right)r^4}
\, ,
\\[5pt]
P_{22}&=
-\frac{(\zeta  \lambda_{\rm m}+2 \left(\lambda_{\rm M}+2 \mu_{\rm M}\right))x_2}{4 \pi  \left(\lambda_{\rm m}+2 \mu_{\rm m}\right) \left(\lambda_{\rm M}+2 \mu
_M\right)r^2}+\frac{(\lambda_{\rm M}+\mu_{\rm M}) x_2 \left(x_1^2-x_2^2\right)}{4 \pi\, \mu_{\rm m} \left(\lambda_{\rm M}+2 \mu_{\rm M}\right)r^4} \, .
\end{split}
\end{align}

\paragraph{Classical linear elasticity  ($L_{\rm c} \to \infty$) - upper bound microscopic stiffness}
As $L_{\rm c} \to \infty$ we have also that $\ell_j \to \infty$ ($j=1,2$), and in this case we obtain that (see Appendix \ref{Bessel})
\begin{align}
\label{CFLL}
u_1=&
\frac{\left(\lambda_{\rm m}+\mu_{\rm m}\right)}{4 \pi\, \mu_{\rm m} \left(\lambda_{\rm m}+2 \mu_{\rm m}\right)} \frac{x_1 x_2}{r^2} + \frac{(\kappa_{\rm e}-\mu_{\rm c})}{4 \pi  \left(\mu_{\rm c}+\mu_{\rm e}\right) \left(\kappa_{\rm e}+\mu_{\rm e}\right)} \frac{x_1 x_2}{r^2}\, , \notag\\
u_2=&
\frac{(\lambda_{\rm m}+\mu_{\rm m})}{4 \pi  \mu_{\rm m} \left(\lambda_{\rm m}+2 \mu_{\rm m}\right)}\frac{x_2^2}{r^2}-\frac{(\lambda_{\rm m}+3 \mu_{\rm m})}{4\pi\, \mu (\lambda_{\rm m}+2\mu_{\rm m})}\ln r +
 \frac{(\kappa_{\rm e}-\mu_{\rm c})}{4 \pi  \left(\mu_{\rm c}+\mu_{\rm e}\right) \left(\kappa_{\rm e}+\mu_{\rm e}\right)}\frac{x_2^2}{r^2}
 \\*[3pt]
&-\frac{\mu_{\rm c}+\kappa_{\rm e}+2 \mu_{\rm e}}{4 \pi  \left(\mu_{\rm c}+\mu_{\rm e}\right) \left(\kappa_{\rm e}+\mu_{\rm e}\right)} \, \ln r \, . \notag
\end{align}
The first two terms in the displacements \eqref{CFLL} are the classical linear elasticity terms (see \eqref{CLF}) but with the micro Lamé moduli ($\mu_{\rm m}, \kappa_{\rm m}$) instead of the macro ones. The other two terms depend also upon the rest of the parameters. 

Furthermore, we obtain the components of the micro-distortion tensor $P$ depending only on the microscopic moduli $(\mu_{\rm m}, \kappa_{\rm m})$ as
\begin{equation}
\begin{aligned}
&P_{11}=\frac{\kappa_{\rm m}x_2(x_2^2-x_1^2)}{4 \pi\, \mu_{\rm m}(\kappa_{\rm m}+\mu_{\rm m})r^4} \, , \qquad \qquad P_{12}=-\frac{\kappa_{\rm m} x_1(x_2^2-x_1^2)}{4 \pi\, \mu_{\rm m}(\kappa_{\rm m}+\mu_{\rm m})r^4}\, ,\\
&P_{21}=-\frac{x_1(x_1^2 \left(\kappa_{\rm m}+2 \mu_{\rm m}\right)+x_2^2 \left(3 \kappa_{\rm m}+2 \mu_{\rm m}\right))}{4 \pi\, \mu_{\rm m}(\kappa_{\rm m}+\mu_{\rm m})r^4}\, , \qquad  P_{22}=\frac{x_2 \left(x_1^2 \left(\kappa_{\rm m}-2 \mu_{\rm m}\right)-x_2^2 \left(\kappa_{\rm m}+2 \mu_{\rm m}\right)\right)}{4 \pi 
    \mu_{\rm m} \left(\kappa_{\rm m}+\mu_{\rm m}  \right)r^4}\,.
\end{aligned}
\end{equation}
It should be noted that the displacement solution does not depend not only upon the micro-moduli since the body force is not zero in this case which corroborates with the findings in the Appendix A.2.2.

\subsection{Concentrated couple}

We consider again a body occupying the full plane under plane-strain conditions. The body is now acted upon by a concentrated line unit couple situated at the origin of the coordinate system. In this case, we have 
\begin{align}
\label{fM}
f=\left(
\begin{array}{c}
0 \\
0
\end{array}
\right), 
\qquad M=
\left(
\begin{array}{cc}
0 & 1/2\\
-1/2 & 0
\end{array}
\right)\delta(x_1)\delta(x_2)\,,
\end{align}
such that $M_{12}-M_{21}=1 \cdot \delta(x_1)\delta(x_2)$. Note that the couple is defined in the relaxed micromorphic as the skew symmetric part of $M$. The diagonal components $M$ do not contribute to the couple as they are self equilibrated double forces (see also Mindlin \cite{Mindlin1964}, Section 4).

Applying the Fourier transform on the equilibrium equations \eqref{pstrain} and solving the non-homogeneous algebraic system yields the following solutions for the transformed field variables
\begin{align}
\widehat{u}_1&=-\frac{i \xi_2}{2\mu_{\rm M} \upxi^2}+\frac{i \xi_2}{2\mu_{\rm e} \left(\ell_2^{-2}+\upxi^2\right)}, \qquad \quad
\widehat{u}_2=\frac{i \xi_1}{2\mu_{\rm M} \upxi^2}-\frac{i \xi_1}{2\mu_{\rm e} \left(\ell_2^{-2}+\upxi^2\right)},\\
\widehat{P}_{11}&=-\widehat{P}_{22}=-\frac{\xi_1 \xi_2}{2\mu_{\rm m} \upxi^2}, \quad
\widehat{P}_{12}=-\frac{\xi_2^2}{2\mu_{\rm m} \upxi^2}-\frac{1}{\widetilde{a} \mu_{\rm M}\, L_{\rm c}^2 \left(\ell_2^{-2}+\upxi^2\right)},\quad
\widehat{P}_{21}=\frac{\xi_1^2}{2\mu_{\rm m} \upxi^2}+\frac{1}{\widetilde{a} \mu_{\rm M}\, L_{\rm c}^2 \left(\ell_2^{-2}+\upxi^2\right)}. \notag
\end{align}
Note that the solution does not depend upon the parameters $\lambda_{\rm e}$ and $\lambda_{\rm m}$, which is to be expected due to the dominant shear character of the loading. Inverting the transformed fields we obtain the following solution for the kinematical fields 

\begin{align}
\label{CMu}
u_1&=-\frac{x_2}{4 \pi \mu_{\rm M} r^2} \left(1-\frac{\mu_{\rm M}}{\mu_{\rm e}} \frac{r}{\ell_2} K_1\left[\frac{r}{\ell_2}\right] \right),&
u_1&=\frac{x_1}{4 \pi \mu_{\rm M} r^2} \left(1-\frac{\mu_{\rm M}}{\mu_{\rm e}} \frac{r}{\ell_2} K_1\left[\frac{r}{\ell_2}\right] \right),\notag\\
P_{11}&=-P_{22}=\frac{x_1 x_2}{2 \pi\, \mu_{\rm m}  r^4},&&\\
P_{12}&=\frac{x_2^2-x_1^2}{4 \pi\, \mu_{\rm m}  r^4}-\frac{1}{2 \pi\,\widetilde{a}\,\mu_{\rm M}L_{\rm c}^2}\, K_0\left[\frac{r}{\ell_2}\right], & 
P_{21}&=\frac{x_2^2-x_1^2}{4 \pi\, \mu_{\rm m}  r^4}+\frac{1}{2 \pi\,\widetilde{a}\,\mu_{\rm M}L_{\rm c}^2}\, K_0\left[\frac{r}{\ell_2}\right].\notag
\end{align}
The micro-rotation is given as 
\begin{equation}
\label{microrotM}
\vartheta_3=\frac{1}{2}(P_{21}-P_{12})=\frac{1}{2 \pi\,\widetilde{a}\,\mu_{\rm M}L_{\rm c}^2}\,   K_0\left[\frac{r}{\ell_2}\right] .
\end{equation}
The stresses and higher order stresses can be derived from the constitutive relations \eqref{const}.

Regarding the asymptotic behaviour of the kinematical fields, we remark that as $r \to 0$ the displacements behave as $r^{-1}$, the micro-distortions $P$ behave as $r^{-2}$, and the micro-rotation exhibits a logarithmic singularity due to the $K_0$-Bessel function. In particular, the modulus of the displacement vector depends (in all theories) only upon the radial distance $r$ and there is no angular dependence (see Figure \ref{fig1z} and Figure \ref{fig2z}). Interestingly, according to the equations \eqref{const}, the stress components $(\sigma_{11},\sigma_{22})$ are bounded at the point of application of the concentrated couple whereas the shear stresses ($\sigma_{12}, \sigma_{21}$) exhibit a logarithmic singularity as $r \to 0$. Finally, the higher order moment stresses $(m_{13},m_{23})$ behave as $O(r^{-1})$
at the origin. All quantities converge to the classical linear elasticity solution (c.f. section 4.2.2.3) as we move away from the concentrated load.

Finally, we note that there is an alternative indirect way to induce the concentrated couple in the relaxed micromorhic medium by superimposing four unit forces (double dipole) in a rotational manner with infinitesimal lever arms. This is the  way to induce the concentrated couple in  classical elasticity  (see e.g. Timoshenko \cite{timoshenko1970}, p. 131). Here however we chose the natuaral way to induce the couple through the skew symmetric part of $M$.

\subsubsection{The pure relaxed micromorphic continuum with symmetric force stress tensor}
The special case of a pure relaxed micromorphic continuum with symmetric force stress tensor is derived by setting $\mu_c =0$ ($\ell_2 \to \infty$). In this case, we have according to \eqref{chlengths} that (see Appendix \ref{Bessel})
\begin{equation}
\label{limitpure2}
\lim_{\mu_{\rm c} \to 0}\frac{1}{\ell_2}=\lim_{\mu_{\rm c} \to 0}\sqrt{\frac{4 \mu_{\rm e} \mu_{\rm c}}{\widetilde{a} \mu_{\rm M}L_{\rm c}^2 (\mu_{\rm e}+\mu_{\rm c})}}=0, \qquad\quad \lim_{\mu_{\rm c} \to 0}\frac{1}{\ell_2}K_1\left[\frac{r}{\ell_2}\right]=\frac{1}{r},
\end{equation}
since $\lim\limits_{z\to 0}z\,K_1(z)=1$ (cf. \eqref{asymptsmall}) and employing \eqref{CMu} together with \eqref{defmod}, we finally derive
\begin{align}\label{pureCM}
u_1&=-\frac{x_2}{4 \pi\, \mu_{\rm m}\, r^2}, \qquad u_2=\frac{x_1}{4 \pi\, \mu_{\rm m}\, r^2}\, , \qquad \text{tr}(\text{D}u)=\text{div}\,u=0 \,,\\
P_{11}&=-P_{22}=\frac{x_1 x_2}{2 \pi\, \mu_{\rm m}  r^4},\quad
P_{12}=\frac{x_2^2-x_1^2}{4 \pi\, \mu_{\rm m}  r^4}+\frac{1}{2\pi\,\widetilde{a}\,\mu_{\rm M}\, L_{\rm c}^2} (\ln r+b),\quad
P_{21}=\frac{x_2^2-x_1^2}{4 \pi\, \mu_{\rm m}  r^4}-\frac{1}{2\pi\,\widetilde{a}\,\mu_{\rm M}\, L_{\rm c}^2} (\ln r+b), \notag
\end{align}
where the last two expressions for $P_{12}$ and $P_{21}$ were derived by taking the limit $\mu_{\rm c} \to 0$ directly in the transformed expressions of the pertinent field variables: Indeed, in the case of a concentrated couple \eqref{fM}, the Fourier system \eqref{FTPL23} has a solution of the form:
\begin{equation}
\begin{aligned}
&\widehat{u}_1=-\frac{i \xi_2}{2 \mu_{\rm m} \upxi^2}\, , \qquad \widehat{u}_2=\frac{i \xi_1}{2 \mu_{\rm m} \upxi^2}\, ,
\\
&\widehat{P}_{11}=-\widehat{P}_{22}=-\frac{\xi_1 \xi_2}{2 \mu_{\rm m} \upxi^2}, \qquad \widehat{P}_{12}=-\frac{2 \mu_{\rm m}+ \mu_{\rm M}\, \widetilde{a} L_{\rm c}^2  \xi_2^2}{2  \mu_{\rm m}\mu_{\rm M} \widetilde{a} L_{\rm c}^2 \upxi^2}\, , \qquad  \widehat{P}_{21}=\frac{2 \mu_{\rm m}+ \mu_{\rm M}\,\widetilde{a} L_{\rm c}^2  \xi_1^2}{2  \mu_{\rm m}\mu_{\rm M} \widetilde{a} L_{\rm c}^2 \upxi^2}.
\end{aligned}
\end{equation}
Using the results in \eqref{FPs} we can readily invert the above expressions and obtain the results in \eqref{pureCM}.

Finally, the micro-rotation is given as 
\begin{equation}
\label{pureCMtheta}
\vartheta_3=-\frac{1}{2\pi\,\widetilde{a}\,\mu_{\rm M}\, L_{\rm c}^2} (\ln r+b).
\end{equation}
According to \eqref{pureCMtheta}, the constant term related to the Euler's constant $b$ corresponds to a constant (rigid) micro-rotation and does not affect the stresses or higher order stresses in \eqref{const}, therefore it can be ignored. It is interesting to note that the displacement field in the pure relaxed micromorphic case \eqref{pureCM} does not converge to the classical macroscopic elasticity one (see \eqref{CMclu}) far away from the concentrated couple. Indeed, the former has in the denominator $\mu_{\rm m} $ and the latter $\mu_{\rm M}$ which means that limits are different as $r \to \infty$. This is not the case however with the complete relaxed micromorphic model (with $\mu_{\rm c}>0$) where, as $r \to \infty$ the Bessel functions in \eqref{pureCM}$_1$ and \eqref{pureCM}$_2$ tend to zero and the classical linear elasticity solution is restored.

It is intriguing to see that setting $\mu_{\rm c}=0$ in the concentrated couple problem acts like a zoom into the microstructure and activates the microscale shear modulus $\mu_{\rm m}$ in the displacement solution, which is not the case in the concentrated force problem.

\subsubsection{Limiting cases}
From the general relaxed micromorphic solution we can derive the fundamental solutions in other generalized continua as singular limiting cases.

\paragraph{Micro-stretch, micropolar and couple stress elasticity}
As $\mu_{\rm m} \to \infty$ we have that: $\mu_{\rm e} \to \mu_{\rm M}$ and  
\begin{equation}\begin{split}
u_1&=-\frac{x_2}{4 \pi\, \mu_{\rm M}\, r^2}\left(1-\frac{r}{\ell_2}K_1\left[\frac{r}{\ell_2}\right] \right),\qquad \quad
u_2=\frac{x_1}{4 \pi\, \mu_{\rm M}\, r^2}\left(1-\frac{r}{\ell_2}K_1\left[\frac{r}{\ell_2}\right] \right),\\
P_{11}&=P_{22}=0, \qquad P_{12}=-P_{21}=-\frac{\mu_{\rm c}+\mu_{\rm e}}{8 \pi\, \mu_{\rm c} \mu_{\rm e} \ell_2^2} K_0\left[\frac{r}{\ell_2}\right].
\end{split}
\end{equation}
This is the micro-stretch solution. Further, if we identify $\mu_c=\alpha$ the solution transforms to the micropolar solution with the characteristic length given by 
\eqref{lmp} \cite{Now72, Dyszlewicz2004}. Next, taking $\mu_{ \rm c} \to \infty$ we derive the couple stress solution \cite{weitsman1968two, hattori2023isogeometric} which is identical in form with the micro-stretch/micropolar solution but with the characteristic length given by \eqref{lcs}. It is worth noting that in the micro-stretch, micropolar, and couple stress theories the displacement field remains bounded and in particular becomes zero at the point of application of the concentrated couple (i.e. $ r \to 0$) which is in marked contrast with the respective relaxed micromorphic behaviour. As $\ell_2 \to \infty$ all the fields become null. Finally, the micro-rotation is given by  \eqref{microrotM} in all cases and exhibits a logarithmic singularity at the origin. As we move away from the load all solutions converge to the classical elasticity solution (section 4.2.2.2). 

\begin{figure}
\begin{minipage}[b]{0.45\textwidth}
\centering
\includegraphics[width=14em]{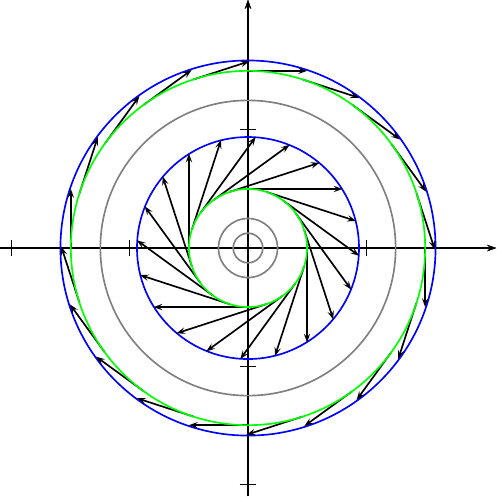}
 \caption{Inhomogeneous displacement solution for the concentrated couple. Circles are rotated and expanded by the deformation around zero.}\label{fig1z}
\end{minipage}
\hfill
\begin{minipage}[b]{0.45\textwidth}
 \centering
 \includegraphics[width=14em]{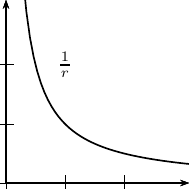}
 \caption{ $\lVert u \rVert$ behaves like $\frac{1}{r}$ in the relaxed micromorphic theory for the case of a concentrated couple.}
 \label{fig2z}
\end{minipage}
\end{figure}

\paragraph{Classical linear elasticity ($L_{\rm c} \to 0$) - lower bound macroscopic stiffness}
As $L_{\rm c} \to 0$ at $\mu_{\rm c}>0$ we have that $\ell_2 \to 0$,  and also (see Appendix \ref{Bessel})
\begin{equation}
\label{limclCM}
\lim_{\ell_2 \to 0} \ell_2^{-2} K_0\left[\frac{r}{\ell_2}\right]=0 \, , \quad \qquad \lim_{\ell_2 \to 0}\frac{1}{\ell_2} K_1\left[\frac{r}{\ell_2}\right]=0 .
\end{equation}
Accordingly, from  \eqref{CMu}, we obtain the standard classical elasticity result for the displacements\footnote{Timoshenko and Goodier \cite[p. 131]{timoshenko1970}; Love \cite[p. 214]{Love1927}.}
\begin{equation}
\label{CMclu}
u_1=-\frac{x_2}{4 \pi\, \mu_{\rm M}\, r^2}, \quad \qquad u_2=\frac{x_1}{4 \pi\, \mu_{\rm M}\, r^2},
\end{equation}
see Fig. \ref{fig1z}. In addition,
\begin{equation}
P_{11}=-P_{22}=\frac{x_1 x_2}{2 \pi\, \mu_{\rm m}  r^4}, \qquad \,\, P_{12}=P_{21}=\frac{x_2^2-x_1^2}{4 \pi\, \mu_{\rm m}  r^4}.
\end{equation}

\paragraph{Classical linear elasticity ($L_{\rm c} \to \infty$) - upper bound microscopic stiffness}
As $L_{\rm c} \to \infty$ ($\ell_2 \to \infty$) we have that
\begin{equation}
\label{limclCM1}
\lim_{\ell_2 \to \infty} \ell_2^{-2} K_0\left[\frac{r}{\ell_2}\right]=0 \, , \quad \qquad \lim_{\ell_2 \to \infty}\frac{1}{\ell_2} K_1\left[\frac{r}{\ell_2}\right]=\frac{1}{r}.
\end{equation}
Accordingly, from  \eqref{CMu}, we obtain the classical elasticity solution for the displacements but now with $\mu_{\rm m}$ instead of $\mu_{\rm M}$
\begin{equation}
\label{CMclu2}
u_1=-\frac{x_2}{4 \pi\, \mu_{\rm m} r^2}, \qquad u_2=\frac{x_1}{4 \pi\, \mu_{\rm m} r^2}.
\end{equation}
In addition, we derive again
\begin{equation}
P_{11}=-P_{22}=\frac{x_1 x_2}{2 \pi\, \mu_{\rm m}  r^4}, \qquad \,\, P_{12}=P_{21}=\frac{x_2^2-x_1^2}{4 \pi\, \mu_{\rm m}  r^4}\,,
\end{equation}
also only depending on the microscopic modulus $\mu_{\rm m}$. Note that since the body force is zero in this problem, the solution when $L_{\rm c} \to \infty$ depends only on the micromoduli as is expected (see Appendix A.2.2).

\section{Fundamental solution for an isotropic gauge-invariant incompatible elasticity model in plane strain}
We consider the gauge-invariant incompatible linear elasticity model \cite{knees2023global,neff2015relaxed,LazarAnastassiadis}  
\begin{align}
\mathbb{C}_{\rm e}\,\text{sym}\,e
+\mathbb{C}_{\rm c}\,\text{skew}\,e
+\mu_{\rm M}L_{\rm c}^2\,\text{Curl}\left(\mathbb{L}_{\rm c}\,\text{Curl}\,e\right)
=
M,
\qquad
e\times n\rvert_{\partial \Omega}=0 \,.
\label{eq:gouge_equili}
\end{align}
where $e\coloneqq{\rm D} u-P:\Omega\in\mathbb{R}^3\to\mathbb{R}^{3\times3}$ is the incompatible elastic distortion, and $\mathbb{C}_{\rm e},\mathbb{C}_{\rm c},\mathbb{L}$ are fourth order tensors as in (\ref{eq:energy_rmm}), while $M$ is similar as in (\ref{eq:equi_RM}).
Due to $\text{Div}\,\text{Curl}=0$, smooth solutions of (\ref{eq:gouge_equili}) satisfy the additional balance equation
\begin{align}
\text{Div}
\big(
\underbrace{
\mathbb{C}_{\rm e}\,\text{sym}\,e
+\mathbb{C}_{\rm c}\,\text{skew}\,e
}
_{\eqqcolon\sigma}
\big)
=
\text{Div} \, M
\eqqcolon \overline{f}.\label{eq:gouge_equili2}
\end{align}
Formally, \eqref{eq:gouge_equili} and \eqref{eq:gouge_equili2} appear as Euler-Lagrange equations of \eqref{eq:energy_rmm} with $\mathbb{C}_{\rm micro}\equiv0$. Substituting a compatible elastic distortion, $e=\text{D}u$, we retrieve from (\ref{eq:gouge_equili}) linear Cauchy elasticity with stiffness tensor $\mathbb{C}_{\rm e}$
\begin{align}
\text{Div}\,
\mathbb{C}_{\rm e}\,\text{sym}\,\text{D}u
=
\overline{f} \, ,
\qquad\qquad\qquad
\text{D}u\times n\rvert_{\partial \Omega}=0 \, .
\end{align}
Observe that the boundary value problem (\ref{eq:gouge_equili}) is still well-posed in terms of the elastic distortion $e$, due to the generalized incompatible Korn's inequality \cite{GLN2}. In the isotropic case (\ref{eq:gouge_equili}) reduces to
\begin{align}
2\mu_{\rm e}\,&\text{dev}\,\text{sym}\,e
+2\mu_{\rm c}\,\text{skew}\,e
+\kappa_{\rm e}\,\text{tr}(e)\,\mathbbm{1}_3
\label{eq:gouge_equili_iso}
\\*
&
+2\mu_{\rm M}L_{\rm c}^2\text{Curl}\left( a_1 \, \text{dev sym} \, \text{Curl} \, e +
a_2 \, \text{skew} \, \text{Curl} \, e +
\frac{a_3}{3} \, \text{tr} \left(\text{Curl} \, e\right)\mathbbm{1}\right)=M\,,
\notag
\end{align}
and this is the second balance equation from \eqref{eq:equi_RM}$_2$ for $\mu_{\rm m} \to 0$, $\kappa_{\rm m} \to 0$ and therefore $\sigma_{\rm micro}\equiv0$.\\
Fundamental solutions to (\ref{eq:gouge_equili_iso}) in the three-dimensional case have been obtained by Lazar \cite{lazar2009gauge} under the constitutive assumption of a strictly positive Cosserat couple modulus, $\mu_{\rm c}>0$.
The latter condition entails that 
\begin{align}
\sigma=
2\mu_{\rm e}\,\text{dev}\,\text{sym}\,e
+2\mu_{\rm c}\,\text{skew}\,e
+\kappa_{\rm e}\,\text{tr}(e)\,\mathbbm{1}_3
\,
\end{align}
can be algebraically inverted, i.e. we can express $e=\mathcal{G}(\sigma)$ if $\mu_{\rm e},\mu_{\rm c},\kappa_{\rm e}>0$, see (\cite{neff2015relaxed}).
Here, we will consider the fundamental solution to (\ref{eq:gouge_equili_iso}) in plane strain, but we allow for $\mu_{\rm c}\geq 0$.
The plane strain version of (\ref{eq:gouge_equili_iso}) is obtained by considering the following energy, connected to (\ref{eq:gouge_equili_iso}), namely
\begin{align}
&\int_{\Omega}
\, \mu_{\rm e} \left\lVert \text{dev} \,\text{sym} \,\widetilde{e} \right\rVert^{2}
+ \frac{\kappa_{\rm e}}{2} \text{tr}^2 \,(\widetilde{e})
\label{eq:energy_gouge}
\\*
&\phantom{\int_{\Omega}{}}
+
\frac{\mu_{\rm M}\, L_{\rm c}^2}{2}\left(
a_1 \, \left\lVert \text{dev sym} \, \text{Curl} \, \widetilde{e}\right\rVert^2 +
a_2 \, \left\lVert \text{skew} \,  \text{Curl} \, \widetilde{e}\right\rVert^2 +
\frac{a_3}{3} \, \text{tr}^2 \left(\text{Curl} \, \widetilde{e}\right)
-\langle M,\widetilde{e} \rangle\right)
\, \text{d}x
\to \text{min}\,\widetilde{e}
.
\notag
\end{align}
As can be seen, letting $L_{\rm c}\to \infty$ while assuming $a_1,a_2,a_3>0$ implies $\text{Curl }\widetilde{e}\equiv0$ and therefore $\widetilde{e}=\text{D}\widetilde{u}$ on contractible domains.
We will consider (\ref{eq:energy_gouge}) in an unbounded domain with given $M=\delta \times \widetilde{M}$.
Similarly, as in section 3, the plane strain energy becomes
\begin{align}
\int_{\Omega}
\, \mu_{\rm e} \left\lVert \text{dev}_2 \,\text{sym} \,\widetilde{e}^{\sharp} \right\rVert^{2}
+ \frac{\kappa_{\rm e}}{2} \text{tr}^2 \,(\widetilde{e}^{\sharp})
\label{eq:energy_gouge_2D}
+
\mu_{\rm M}L_{\rm c}^2\,\widetilde{a}\,
\lVert\text{Curl}_{\text{2D}} \, \widetilde{e}^{\sharp}\rVert^2
-\langle M,\widetilde{e}^{\sharp}\rangle
\, \text{d}x
\to \text{min}\,\widetilde{e}^{\sharp}.
\end{align}
and we obtain the plane strain equations in components
\begin{align}\begin{split}
-\mu_{\rm M}\, L_{\rm c}^2 \widetilde a \left(e_{11,22}-e_{12,12}\right)+(\lambda_{\rm e}+2 \mu_{\rm e}) e_{11}+\lambda_{\rm e} e_{22}
&=M_{11} \,,
\\*
\mu_{\rm M}\, L_{\rm c}^2 \widetilde a \left(e_{11,12}-e_{12,11}\right)+(\mu_{\rm c}+\mu_{\rm e}) e_{12}+(\mu_{\rm e}-\mu_{\rm c}) e_{21}
&=M_{12} \,,
\\*
- \mu_{\rm M}\, L_{\rm c}^2 \widetilde a \left(e_{21,22}-e_{22,12}\right)+(\mu_{\rm e}-\mu_{\rm c}) e_{12}+(\mu_{\rm c}+\mu_{\rm e}) e_{21}
&=M_{21} \,,
\\*
 \mu_{\rm M}\, L_{\rm c}^2 \widetilde a \left(e_{21,12}-e_{22,11}\right)+\lambda_{\rm e} e_{11}+(\lambda_{\rm e}+2 \mu_{\rm e}) e_{22}
&=M_{22} \,.\end{split}
\end{align}

We consider again the case of a concentrated line unit couple situated at the origin of the coordinate system. In this case, the components of the body volume moment $M$ are given by  \eqref{fM}. Following an analogous Fourier transform analysis as in the previous cases we derive the fundamental solution for a concentrated couple in gauge-invariant incompatible elasticity. The incompatible elastic distortions read then
\begin{equation}
\begin{split}
e_{11}&=-e_{22}=-\frac{x_1 x_2}{4 \pi (\mu_{\rm c}+\mu_{\rm e})\ell_2^2 \, r^2} K_2\left[\frac{r}{\ell_2}\right] \, ,\\
e_{12}&=\frac{\mu_{\rm e}}{8 \pi\, \mu_{\rm c} (\mu_{\rm c}+\mu_{\rm e})\ell_2^2} K_0\left[\frac{r}{\ell_2}\right]+\frac{x_1^2-x_2^2}{8 \pi (\mu_{\rm c}+\mu_{\rm e})\ell_2^2 \, r^2} K_2\left[\frac{r}{\ell_2}\right]\, ,\\
e_{21}&=-\frac{\mu_{\rm e}}{8 \pi\, \mu_{\rm c} (\mu_{\rm c}+\mu_{\rm e})\ell_2^2} K_0\left[\frac{r}{\ell_2}\right]+\frac{x_1^2-x_2^2}{8 \pi (\mu_{\rm c}+\mu_{\rm e})\ell_2^2 \, r^2} K_2\left[\frac{r}{\ell_2}\right]\, .
 \end{split}
\end{equation}
where $\ell_2$ is given by \eqref{chlengths}. It is interesting to note that the solution does not depend upon the elastic bulk modulus $\kappa_{\rm e}$ and that the elastic distortion tensor for the case of a concentrated couple is traceless (i.e. $\text{tr}(\widetilde{e}^{\sharp})=e_{11}+e_{22}=0$). 

\section{Numerical results and discussion}
We will now present some results regarding the behaviour of the relaxed micromorphic solution near the application of the applied loads. A comparison of the results with other well known generalized continua obtained as limiting cases of the general relaxed micromorphic model will also be performed.

The relaxed micromorphic continua under plane strain conditions can be fully described by four dimensionless parameters. In order to have a unified treatment for all the above cases, the following dimensionless quantities $g_i$ ($i=1,2,3,4$) are introduced:
\begin{equation}
\label{nondimpar}
\mu_{\rm e}=g_1\, \mu_{\rm M}, \qquad \mu_{\rm c}=g_2\, \mu_{\rm M}, \qquad \kappa_{\rm e}=g_3\, \mu_{\rm M}, \qquad \kappa_{\rm M}=g_4\, \mu_{\rm M}.
\end{equation}
In view of \eqref{PD}, we have that: $g_1>1$, $g_2\ge 0$, and  $g_3>g_4>0$. We also recall that $\lambda_i=\kappa_i-\mu_i$ with $i\in\{\rm{e}, \rm{m}, \rm{M}\}$ and using \eqref{nondimpar} that
\begin{equation}
\label{nondimpar2}
\mu_{\rm m}=\frac{g_1}{g_1-1}\mu_{\rm M}, \qquad \kappa_{\rm m}=\frac{g_3}{g_3-g_4}\kappa_{\rm M}. 
\end{equation} 
Further, for comparison purposes all distances from the origin are normalized with respect to the characteristic length $\ell_2$ of the relaxed micromorphic model. Results for the cases of a concentrated force and concentrated couple will be shown separately.

\subsection{Concentrated force}
\begin{figure}[!h]
\centering
\includegraphics[scale=0.4]{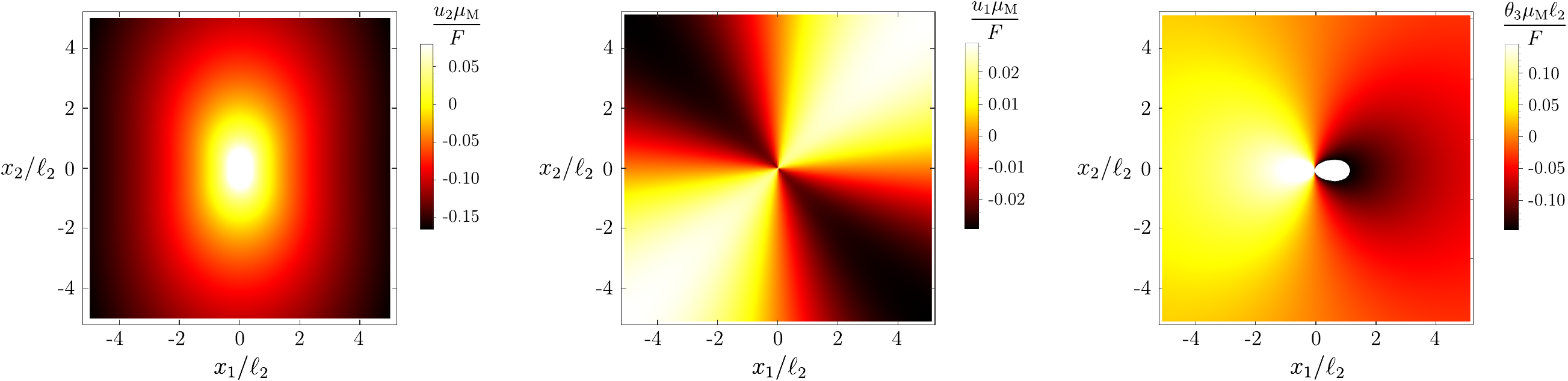}
\caption{Contours of the normalized displacements $\frac{u_i\, \mu_{\rm M}}{F}$ and micro-rotation $\frac{\vartheta_3\, \mu_{\rm M} \ell_2}{F}$ due to a concentrated unit line force ($F=1$) acting at the origin of relaxed micromorphic medium. The material is characterized by $g_1=1.2$, $g_2=3$, $g_3=5$ and $g_4=3$.}
\label{conts}
\end{figure}

\begin{figure}[!h]
\centering
\includegraphics[scale=0.45]{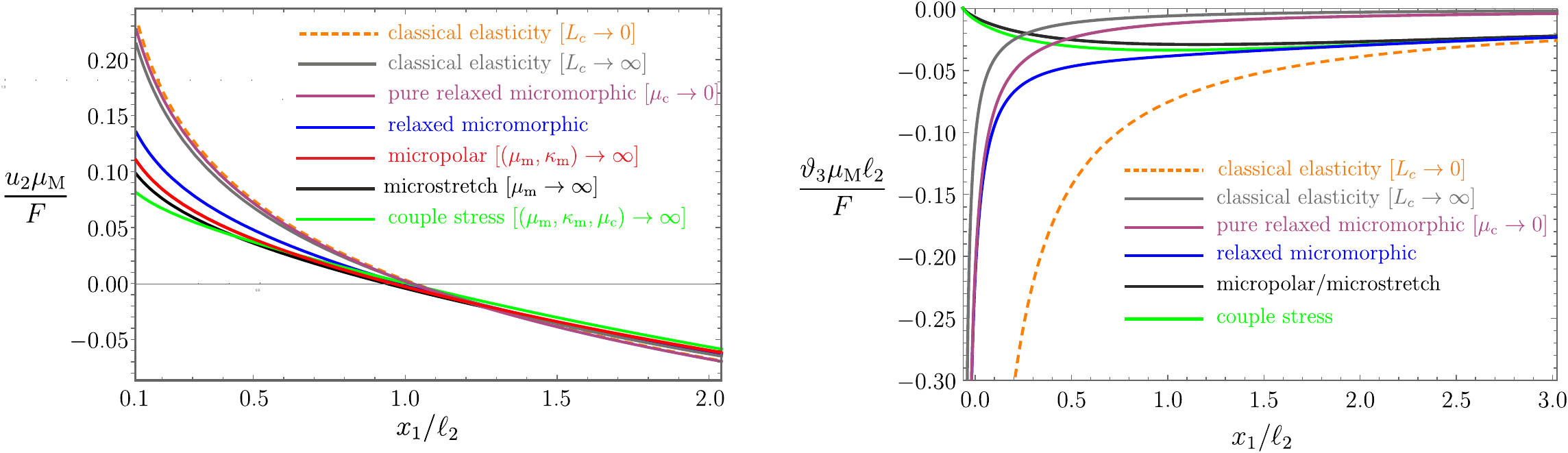}
\caption{Variation of the normalized displacement $\frac{u_2\mu_{\rm M}}{F}$ and the normalized micro-rotation $\frac{\vartheta_3 \mu_{\rm M} \ell_2}{F}$ along the positive $x_1$-axis due to a concentrated unit line force ($F=1$) in various generalized continuum theories. The relaxed micromorphic material is characterized by $g_1=1.2$, $g_2=3$, $g_3=5$ and $g_4=3$.}
\label{comp}
\end{figure}

Figure \ref{conts} shows contours of the normalized displacements and micro-rotation due to a concentrated line force acting at the origin for a relaxed micromorphic material characterized by  ($g_1=1.2,\, g_2=3, \, g_3=5, \, g_4=3$). This implies, according to \eqref{nondimpar2}, that $\mu_{\rm m}=6\mu_{\rm M}$ and $\kappa_{\rm m}=2.5 \kappa_{\rm M}$.   A comparison of the relaxed micromorphic continua with other generalized continua that can be obtained as limiting cases is shown in Figure \ref{comp}. In particular, in Fig. \ref{comp}, the normalized displacement $\frac{u_2\, \mu_{\rm M}}{F}$ and the normalized micro-rotation $\frac{\vartheta_3\, \mu_{\rm M} \ell_2}{F}$ ($F=1$) are plotted along the positive $x_1$-axis (i.e. for $x_2=0$) . The $u_2$ displacement has a logarithmic singularity at the origin in all theories. However, the singularity is eliminated in strain gradient elasticity \cite{gourgiotis2018concentrated}. It is observed that deviations from the classical elasticity solution (dashed line) are more noticeable within a range of $\lvert x_1 \rvert \le 2\ell_2$ from the point of application of the concentrated force. All solutions converge quickly to the classical elasticity solution as we move away from the origin. It is also shown that the classical elasticity and the couple stress elasticity serve as the upper and lower bounds for the solutions. In fact, couple stress elasticity predicts more pronounced size effects as compared to the other generalized continuum theories. The micropolar solution is in-between the classical and the couple stress solution. Also, we note that the relaxed micromorphic and the pure relaxed micromorphic are closer to the classical elasticity one. 

Regarding the behaviour of the micro-rotation we remark that the classical elasticity and the relaxed micromorphic elasticity predict unbounded micro-rotation at the origin which is in marked contrast with couple stress, micropolar, and micro-stretch theories that predict zero micro-rotation at the origin. In all theories the micro-rotation decays as $O(x_1^{-1})$ when $x_1 \to \infty$. However, as it can be seen from Figure 5b, in the pure relaxed micromorphic model and in the classical elasticity model with $L_{\rm c} \to \infty$ (upper bound microscopic stiffness) the solution does not converge in the standard classical elasticity solution ($L_{\rm c} \to 0$) as all  other theories do.

\subsection{Concentrated couple}

\begin{figure}[!h]
\centering
\includegraphics[scale=0.4]{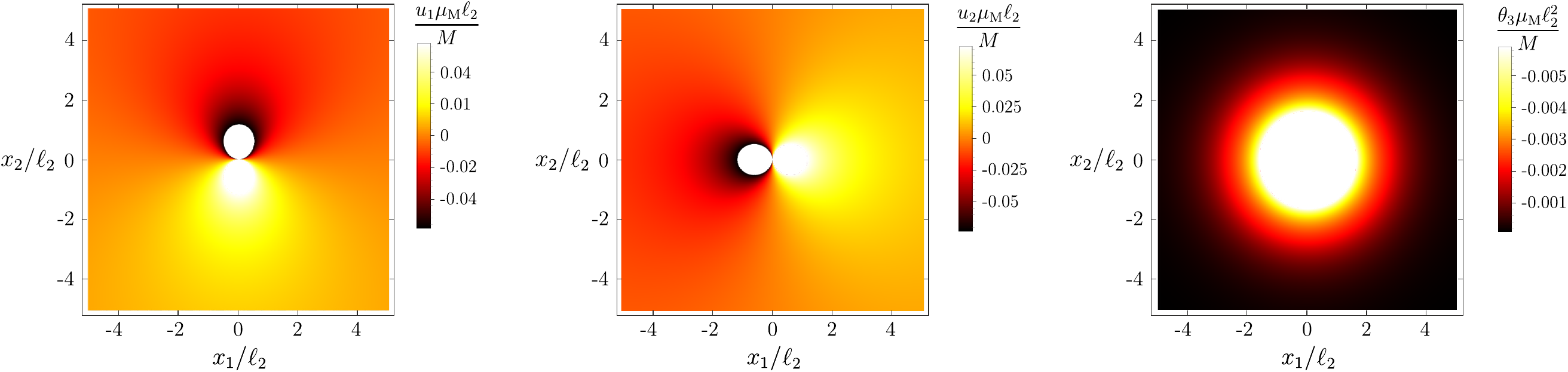}
\caption{Contours of the normalized displacements and micro-rotation due to a concentrated unit line couple ($M=1$) acting at the origin. The relaxed micromorphic material is characterized by $g_1=3$ and $g_2=2$.}
\label{conts2}
\end{figure}

\begin{figure}[!h]
\centering
\includegraphics[scale=0.6]{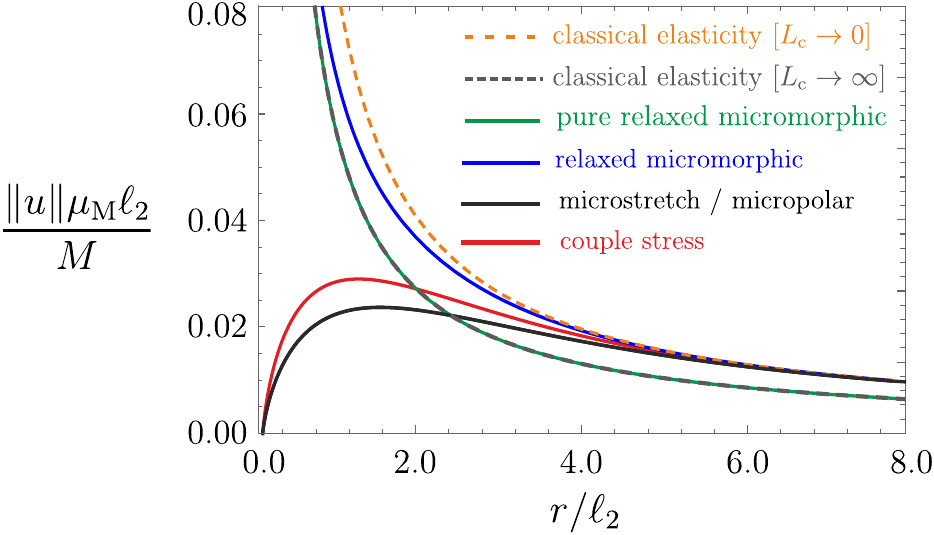}
\caption{Variation of the normalized modulus of the displacement vector $\frac{\lVert u \rVert \mu_{\rm M} \ell_2}{M}$ along the positive $x_1$-axis due to a concentrated unit line couple ($M=1$) in various generalized continuum theories. The relaxed micromorphic material is characterized by $g_1=3$ and $g_2=2$. The gap between the green and black line at the right is due to $\mu_{\rm m}>\mu_{\rm M}$. We note the fundamental qualitative difference between the relaxed micromorphic model and the other generalized continua (microstretch, micropolar, couple stress) in their behaviour near to the singularity.}
\label{comp2}
\end{figure}

Figure \ref{conts2} shows contours of the normalized displacements and micro-rotation for the case of a concentrated couple. In this case, only the parameters $g_1$ and $g_2$ need to be specified. A comparison of the relaxed micromorphic continua with other generalized continua  obtained as limiting cases is also shown in Figure \ref{comp2}. In particular, in Fig. \ref{comp2}, the normalized modulus of the  displacement vector $ \lVert u \rVert$ is plotted against the radial distance $r$. The material parameters for the relaxed micromorphic material are: $g_1=3$ and $g_2=2$ (which implies $\mu_{\rm m}=1.5 \mu_{\rm M}$). All distances are normalized with respect to characteristic length of the relaxed micromorphic theory $\ell_2$.

It is noted that $ \lVert u \rVert$ has a Cauchy type singularity $O(r^{-1})$ in the relaxed micromorphic theory, in the pure relaxed micromorphic, and in the classical elasticity theory ($L_{\rm c} \to 0$ and $L_{\rm c} \to \infty$) but the strengths of the singularities are different. In marked contrast, $ \lVert u \rVert$ is bounded and becomes zero at the origin in micro-stretch, micropolar and couple stress theory. As it was shown analytically (see sections 4.2.1 and 4.2.2.3), only the pure relaxed micromorphic solution and the classical elasticity solution with $L_{\rm c} \to \infty$ (green and dashed-gray lines in Fig. \ref{comp2}) do not converge to the standard classical elasticity ($L_{\rm c} \to 0$) as $r \to \infty$. This is to be expected since the latter solutions depend only upon the micro shear modulus $\mu_{\rm m}$. 

\begin{figure}[H]
\centering
\includegraphics[scale=0.4]{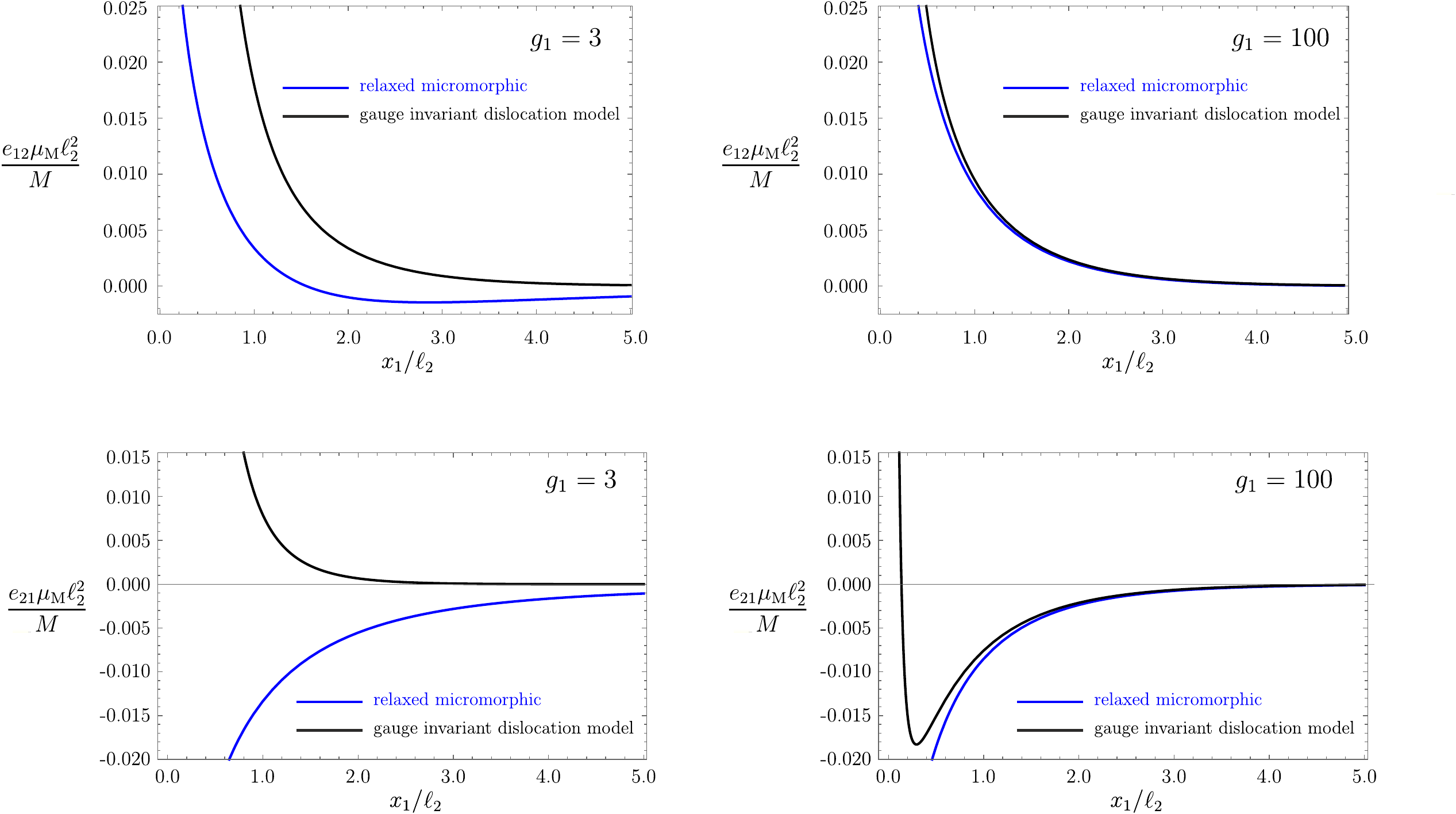}
\caption{Variation of the normalized incompatible elastic shear distortions along the positive axis $x_1$-axis due to a concentrated unit line couple ($M=1$) in the relaxed micromorphic theory and the gauge invariant dislocation model for $g_2=2$ and various values of the parameter $g_1$.}
\label{dis}
\end{figure}

Finally, a comparison of the incompatible elastic distortions $e_{12}=u_{1,2}-P_{12}$ and $e_{21}=u_{2,1}-P_{21}$ in the relaxed micromorphic theory and the gauge invariant dislocation model is shown in Figure \ref{dis}. It is observed that as $g_1$ increases as compared to $g_2$ (i.e. $\mu_{\rm e} \gg \mu_{\rm M}$ and $\mu_{\rm e} \gg \mu_{\rm c}$), the solutions for the gauge invariant dislocation model and the relaxed micromorphic model converge.

\section{Conclusions}

In the present work, the infinite plane 2D Green’s functions for a concentrated force and a concentrated couple have been derived in the context of the isotropic relaxed micromorphic theory. Our main concern here was to determine possible deviations from the
predictions of classical theory of elasticity but also from other generalized continuum theories that are extensively used nowadays for the prediction of scale effects. Closed form solutions were derived using a Fourier transform analysis and results from generalized functions.

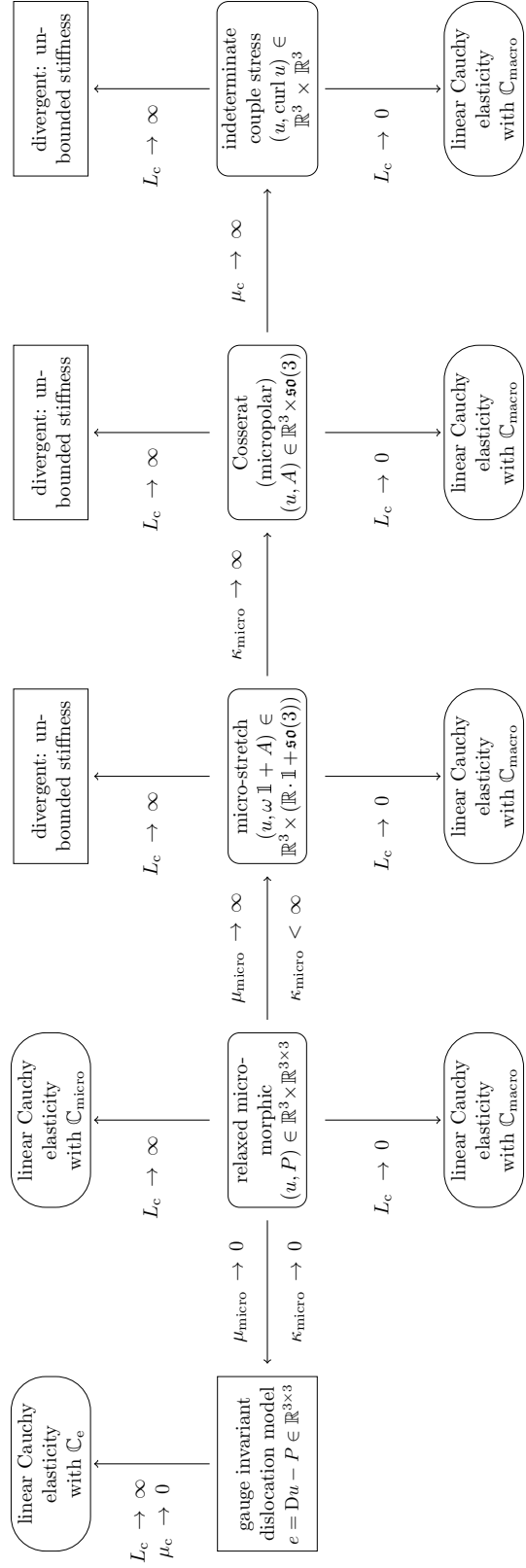
\begin{figure}[H]
\begin{sideways}
\scalebox{1.34}{\resizebox{\textwidth}{!}{%
\begin{tikzpicture}[x=1.0cm,y=1.0cm];
    \coordinate (A) at (0,0);
    \coordinate (Ap1) at ($(A)+(0,1)$);
    \coordinate (Ap2) at ($(Ap1)+(0,2)$);
    \coordinate (Ap3) at ($(Ap2)+(0,0.75)$);
    \coordinate (Ap4) at ($(Ap1)!0.5!(Ap2)+(-1,0)$);
    \coordinate (Am1) at ($(A)-(0,1)$);
    \coordinate (Am2) at ($(Am1)-(0,2)$);
    \coordinate (Am3) at ($(Am2)-(0,0.75)$);
    \coordinate (Am4) at ($(Am1)!0.5!(Am2)+(-1,0)$);
    \coordinate (B) at (6,0);
    \coordinate (Bp1) at ($(B)+(0,1)$);
    \coordinate (Bp2) at ($(Bp1)+(0,2)$);
    \coordinate (Bp3) at ($(Bp2)+(0,0.75)$);
    \coordinate (Bp4) at ($(Bp1)!0.5!(Bp2)+(-1,0)$);
    \coordinate (Bm1) at ($(B)-(0,1)$);
    \coordinate (Bm2) at ($(Bm1)-(0,2)$);
    \coordinate (Bm3) at ($(Bm2)-(0,0.75)$);
    \coordinate (Bm4) at ($(Bm1)!0.5!(Bm2)+(-1,0)$);
    \coordinate (C) at (12,0);
    \coordinate (Cp1) at ($(C)+(0,1)$);
    \coordinate (Cp2) at ($(Cp1)+(0,2)$);
    \coordinate (Cp3) at ($(Cp2)+(0,0.75)$);
    \coordinate (Cp4) at ($(Cp1)!0.5!(Cp2)+(-1,0)$);
    \coordinate (Cm1) at ($(C)-(0,1)$);
    \coordinate (Cm2) at ($(Cm1)-(0,2)$);
    \coordinate (Cm3) at ($(Cm2)-(0,0.75)$);
    \coordinate (Cm4) at ($(Cm1)!0.5!(Cm2)+(-1,0)$);
    \coordinate (D) at (18,0);
    \coordinate (Dp1) at ($(D)+(0,1)$);
    \coordinate (Dp2) at ($(Dp1)+(0,2)$);
    \coordinate (Dp3) at ($(Dp2)+(0,0.75)$);
    \coordinate (Dp4) at ($(Dp1)!0.5!(Dp2)+(-1,0)$);
    \coordinate (Dm1) at ($(D)-(0,1)$);
    \coordinate (Dm2) at ($(Dm1)-(0,2)$);
    \coordinate (Dm3) at ($(Dm2)-(0,0.75)$);
    \coordinate (Dm4) at ($(Dm1)!0.5!(Dm2)+(-1,0)$);
    \coordinate (E) at ($(Cm3)-(0,2.5)$);
    \coordinate (Ep1) at ($(E)+(-1.5,1)$);
    \coordinate (Ep2) at ($(Ep1)+(0,2)$);
    \coordinate (Ep3) at ($(Ep2)+(0,0.75)$);
    \coordinate (Ep4) at ($(Ep1)!0.5!(Ep2)+(-1,0)$);
    \coordinate (Em1) at ($(E)-(0,1)$);
    \coordinate (Em2) at ($(Em1)-(0,2)$);
    \coordinate (Em3) at ($(Em2)-(0,0.75)$);
    \coordinate (Em4) at ($(Em1)!0.5!(Em2)+(-1,0)$);
    \coordinate (Bmz) at ($(Bm1)+(0.5,0)$);
    \coordinate (Emz) at ($(Ep1)!0.5!(Bmz)+(0,1.05)$);
    \coordinate (Ez2) at ($(E)+(1.75,0)$);
    \coordinate (Ez3) at ($(Ez2)+(2.5,0)$);
    \coordinate (Ez4) at ($(Ez2)!0.5!(Ez3)+(0,0.5)$);
    \coordinate (Ez5) at ($(Ez3)+(1.75,0)$);
    \coordinate (F) at (-6,0);
    \coordinate (Fp1) at ($(F)+(0,1)$);
    \coordinate (Fp2) at ($(Fp1)+(0,2)$);
    \coordinate (Fp3) at ($(Fp2)+(0,0.75)$);
    \coordinate (Fp4) at ($(Fp1)!0.5!(Fp2)+(-1,0)$);
    \coordinate (Fm1) at ($(F)-(0,1)$);
    \coordinate (Fm2) at ($(Fm1)-(0,2)$);
    \coordinate (Fm3) at ($(Fm2)-(0,0.75)$);
    \coordinate (Fm4) at ($(Fm1)!0.5!(Fm2)+(-1,0)$);
    \node[draw, rounded corners=6mm, text width=8em, align=center] at (Ap3) {linear Cauchy elasticity with $\mathbb{C}_{\rm micro}$};
    \draw [->, line width=.5pt] (Ap1) -- (Ap2);
    \node[text width=8em, align=center] at (Ap4) {$L_{\rm c}\to\infty$};
    \node[draw, rounded corners=2mm, text width=8em, text height=6pt, text depth=26pt, text centered] at (A) {relaxed micro-morphic \\$(u,P)\in\mathbb{R}^3\times\mathbb{R}^{3\times 3}$};
    \node[text width=8em, align=center] at (Am4) {$L_{\rm c}\to 0$};
    \draw [->, line width=.5pt] (Am1) -- (Am2);
    \node[draw, rounded corners=6mm, text width=8em, align=center] at (Am3) {linear Cauchy elasticity with $\mathbb{C}_{\rm macro}$};
    \node[text width=8em, align=center] at ($(A)!0.5!(B)+(0,0.5)$) {$\mu_{\rm micro}\to \infty$};
    \draw [->, line width=.5pt] ($(A)+(1.75,0)$) -- ($(B)-(1.75,0)$);
    \node[text width=8em, align=center] at ($(A)!0.5!(B)-(0,0.5)$) {$\kappa_{\rm micro} < \infty$};
    \node[draw, text width=8em, text height=12.5pt, text depth=17.5pt, align=center] at (Bp3) {divergent: unbounded stiffness};
    \draw [->, line width=.5pt] (Bp1) -- (Bp2);
    \node[text width=8em, align=center] at (Bp4) {$L_{\rm c}\to\infty$};
    \node[draw, rounded corners=2mm, text width=8em, text height=6pt, text depth=26pt, text centered] at (B) {micro-stretch \\$(u,\omega\,\mathbbm{1}+A)\in\mathbb{R}^3\times(\mathbb{R}\cdot\mathbbm{1}+\mathfrak{so}(3))$};
    \node[text width=8em, align=center] at (Bm4) {$L_{\rm c}\to 0$};
    \draw [->, line width=.5pt] (Bm1) -- (Bm2);
    \node[draw, rounded corners=6mm, text width=8em, align=center] at (Bm3) {linear Cauchy elasticity with $\mathbb{C}_{\rm macro}$};
    \node[text width=8em, align=center] at ($(B)!0.5!(C)+(0,0.5)$) {$\kappa_{\rm micro}\to \infty$};
    \draw [->, line width=.5pt] ($(B)+(1.75,0)$) -- ($(C)-(1.75,0)$);
    \node[draw, text width=8em, text height=12.5pt, text depth=17.5pt, align=center] at (Cp3) {divergent: unbounded stiffness};
    \draw [->, line width=.5pt] (Cp1) -- (Cp2);
    \node[text width=8em, align=center] at (Cp4) {$L_{\rm c}\to\infty$};
    \node[draw, rounded corners=2mm, text width=8em, text height=6pt, text depth=26pt, text centered] at (C) {Cosserat\\(micropolar) \\$(u,A)\in\mathbb{R}^3\times \mathfrak{so}(3)$};
    \node[text width=8em, align=center] at (Cm4) {$L_{\rm c}\to 0$};
    \draw [->, line width=.5pt] (Cm1) -- (Cm2);
    \node[draw, rounded corners=6mm, text width=8em, align=center] at (Cm3) {linear Cauchy elasticity with $\mathbb{C}_{\rm macro}$};
    \node[text width=8em, align=center] at ($(C)!0.5!(D)+(0,0.5)$) {$\mu_{\rm c}\to \infty$};
    \draw [->, line width=.5pt] ($(C)+(1.75,0)$) -- ($(D)-(1.75,0)$);
    \node[draw, text width=8em, text height=12.5pt, text depth=17.5pt, align=center] at (Dp3) {divergent: unbounded stiffness};
    \draw [->, line width=.5pt] (Dp1) -- (Dp2);
    \node[text width=8em, align=center] at (Dp4) {$L_{\rm c}\to\infty$};
    \node[draw, rounded corners=2mm, text width=8em, text height=6pt, text depth=37pt, text centered, align=center] at (D) {indeterminate couple stress \\$(u,\text{curl}\,u)\in\mathbb{R}^3\times \mathbb{R}^3$};
    \node[text width=8em, align=center] at (Dm4) {$L_{\rm c}\to 0$};
    \draw [->, line width=.5pt] (Dm1) -- (Dm2);
    \node[draw, rounded corners=6mm, text width=8em, text centered] at (Dm3) {linear Cauchy elasticity with $\mathbb{C}_{\rm macro}$};
    \node[draw, rounded corners=6mm, text width=8em, align=center] at (Fp3) {linear Cauchy elasticity \\with $\mathbb{C}_{\rm e}$};
    \draw [->, line width=.5pt] (Fp1) -- (Fp2);
    \node[text width=8em, align=center] at (Fp4) {$L_{\rm c}\to\infty$\\$\mu_{\rm c}\to0$};
    \node[draw, text width=8em, text height=12.5pt, text depth=30pt, align=center] at (F) {gauge invariant dislocation model\\\resizebox{\textwidth}{!}{$e=\text{D}u-P \in\mathbb{R}^{3\times3}$}};
    \node[text width=8em, align=center] at ($(A)!0.5!(F)+(0,0.5)$) {$\mu_{\rm micro}\to 0$};
    \draw [->, line width=.5pt] ($(A)-(1.75,0)$) -- ($(F)+(1.75,0)$);
    \node[text width=8em, align=center] at ($(A)!0.5!(F)-(0,0.5)$) {$\kappa_{\rm micro} \to 0$};
\end{tikzpicture}
}}
\end{sideways}
\centering
\caption{Tree of the limit cases of the relaxed micromorphic model in statics.}
\label{tree}
\end{figure}

It is shown that the relaxed micromorphic solution is general enough to encompass several well known generalized continuum models which can be recovered as singular limiting cases. In particular, from the relaxed micromorphic solutions we may readily derive the couple-stress, Cosserat-micropolar, micro-stretch, and classical elasticity fundamental solutions  showing thus how versatile the relaxed micromorphic theory is. Yet, the model retains a crucial level of simplicity so that analytical solutions can be found. It has been shown that the relaxed micromorphic solutions are closer to the classical elasticity solutions showing thus milder size effects as compared to the predictions of other models such as the couple-stress, Cosserat-micropolar, and micro-stretch. Figure \ref{tree} shows a tree of the various generalized continua that can be derived as singular limits of the relaxed micromorphic continua and the paths that these can be obtained. 

Finally, we note that the present solutions  may be used for the solution of more general boundary problems and as fundamental solutions for the Boundary Element Method.

\pagebreak

{\scriptsize
\paragraph*{{\scriptsize Acknowledgements.}}
Angela Madeo and Gianluca Rizzi acknowledge support from the European Commission through the funding of the ERC Consolidator Grant META-LEGO, N$^\circ$ 101001759.
Patrizio Neff acknowledges support in the framework of the DFG-Priority Programme 2256 ``Variational Methods for Predicting Complex Phenomena in Engineering Structures and Materials'' and Neff 902/10-1, Project-No. 440935806.
}

\begingroup
\setstretch{1}
\setlength\bibitemsep{3pt}
\printbibliography
\endgroup

 {\footnotesize
\begin{alphasection}
\section{Appendix}
In this appendix, for the convenience of the reader, we exhibit the two-scale nature of the relaxed micromorphic model in three and two dimensions together with the form of equations and constitutive tensors in plane-strain. 

%
%
\subsection{A true two-scale model: the relaxed micromorphic limit for $L_{\rm c}\to 0$ and $L_{\rm c}\to \infty$ in three dimensions}
The relaxed micromorphic model reduces to a classical Cauchy elasticity model for both $L_{\rm c}\to 0$ and $L_{\rm c}\to \infty$ but with two different well-defined stiffnesses, $\mathbb{C}_{\rm macro}$ and $\mathbb{C}_{\rm micro}$, respectively.
The expressions of these stiffnesses in the isotropic case are presented in the next two sections for the convenience of the reader.
%
%
%
\subsubsection{Limit for $L_{\rm c}\to 0$: lower bound macroscopic stiffness $\mathbb{C}_{\rm macro}$}
For the limit $L_{\rm c} \to 0$, the equilibrium equations (\ref{eq:equi_RM}) reduce to
\begin{align}
\text{Div}
\big[
2\mu_{\rm e}\,\text{sym}  (\text{D}u - P )
+ 2\mu_{\rm c}\,\text{skew} (\text{D}u - P )
+ \lambda_{\rm e} \text{tr} (\text{D}u - P ) \mathbbm{1}
\big]
&=
f
\, ,
\label{eq:equi_RM_Lc_0_1}
\\*
2\mu_{\rm e}\,\text{sym}  (\text{D}u - P )
+ \lambda_{\rm e} \text{tr} (\text{D}u - P ) \mathbbm{1}
+ 2\mu_{\rm c}\,\text{skew} (\text{D}u - P )
- 2 \mu_{\rm micro}\,\text{sym}\,P
- \lambda_{\rm micro} \text{tr} \left( P \right) \mathbbm{1}
&=
M
\, .
\label{eq:equi_RM_Lc_0_2}
\end{align}
The equation (\ref{eq:equi_RM_Lc_0_2}) is now algebraic in $P$. Due to the orthogonality of the ``sym/skew'' decomposition, the equation (\ref{eq:equi_RM_Lc_0_2}) requires that
\begin{align}
2\mu_{\rm c} \, \text{skew} (\text{D}u - P ) 
&= \text{skew} \, M
\, ,
\\
2\mu_{\rm e}\,\text{sym}  (\text{D}u - P )
+ \lambda_{\rm e} \text{tr} (\text{D}u - P ) \mathbbm{1}
- 2 \mu_{\rm micro}\,\text{sym}\,P
- \lambda_{\rm micro} \text{tr} \left( P \right) \mathbbm{1}
&=
\text{sym} \, M
\, .
\notag
\end{align}
Since the ``sym'' operator is not orthogonal to the ``tr'' operator, we further decompose ``sym'' into ``dev sym'' and ``tr sym'' so that
\begin{align}
2\mu_{\rm c} \, \text{skew} (\text{D}u - P ) 
&= \text{skew} \, M
\, ,
\notag
\\
2\mu_{\rm e}\,\text{dev sym}  (\text{D}u - P )
+ \frac{2}{3}\mu_{\rm e}\,\text{tr}  (\text{D}u - P ) \mathbbm{1}
+ \lambda_{\rm e} \text{tr} (\text{D}u - P ) \mathbbm{1}
\quad
&
\label{eq:equi_RM_Lc_0_dev_sym}
\\*
- 2 \mu_{\rm micro}\,\text{dev sym}\,P
- \frac{2}{3} \mu_{\rm micro}\,\text{tr}\,(P) \mathbbm{1}
- \lambda_{\rm micro} \text{tr} \left( P \right) \mathbbm{1}
&=
\text{dev sym} \, M
+ \frac{1}{3}\text{tr}\,(M) \mathbbm{1}
\, ,
\notag
\end{align}
note that ``tr sym'' is the same as ``tr''.
We also recall here the definition of the volumetric part, the deviatoric part, and the skew-symmetric parts in the 3D case
\begin{align}\begin{split}
\text{3D volumetric part}
&\coloneqq
\frac{1}{3} \text{tr}(P)\mathbbm{1} \, ,
\\*
\text{3D deviatoric symmetric part}
&\coloneqq
\frac{P+P^{\rm T}}{2}-\frac{1}{3} \text{tr}(P)\mathbbm{1} \, ,
\\*
\text{3D skew symmetric part}
&\coloneqq
\frac{P-P^{\rm T}}{2} \, .\end{split}
\end{align}
With further manipulations and thanks to the orthogonality of the operator ``skew'', ``dev sym'', and ``tr'', the system (\ref{eq:equi_RM_Lc_0_dev_sym}) requires that
\begin{align}
2\mu_{\rm c} \, \text{skew} (\text{D}u - P ) 
&=
\text{skew}\,M
\, ,
\notag
\\
2\mu_{\rm e}\,\text{dev sym}  (\text{D}u - P )
- 2\mu_{\rm micro}\,\text{dev sym}\,P
&=
\text{dev sym}\,M
\label{eq:equi_RM_Lc_0_dev_sym_2}
\, ,
\\*
\left(\frac{2}{3}\mu_{\rm e}+ \lambda_{\rm e}\right)\,\text{tr}  (\text{D}u - P ) \mathbbm{1}
- \left(\frac{2}{3} \mu_{\rm micro} + \lambda_{\rm micro}\right)\,\text{tr}\,(P) \mathbbm{1}
&=
\frac{1}{3}\text{tr}\,(M) \mathbbm{1}
\, .
\notag
\end{align}
From equation (\ref{eq:equi_RM_Lc_0_dev_sym_2}) we can evaluate the expressions for $\text{skew} \, P$, $\text{dev sym} \, P$, and $\text{tr}(P)$ individually as
\begin{align}
\text{skew} \, \text{D}u
-
\frac{1}{2 \mu_{\rm c}}\text{skew}\,M
&=
\text{skew} \, P
\, ,
\notag
\\*
\frac{\mu_{\rm e}}{\mu_{\rm e}+\mu_{\rm micro}}
\text{dev sym}  \, \text{D}u
-
\frac{1}{2 (\mu_{\rm e}+\mu_{\rm micro})}\text{sym}\,M
&=
\text{dev sym}  \, P
\, ,
\label{eq:equi_RM_Lc_0_dev_sym_3}
\\*
\frac{\kappa_{\rm e}}{\kappa_{\rm e}+\kappa_{\rm micro}}
\text{tr}  \, \text{D}u
-
\frac{1}{3 (\kappa_{\rm e}+\kappa_{\rm micro})}\text{tr}(M)
&=
\text{tr}(P)
\, ,
\notag
\end{align}
where  $\kappa_{\rm e}=\frac{2\mu_{\rm e}+3\lambda_{\rm e}}{3}$ and $\kappa_{\rm micro}=\frac{2\mu_{\rm micro}+3\lambda_{\rm micro}}{3}$ are the 3D-elastic and micro bulk modulus, respectively. The contribution of the body volume moment $M$ can be incorporated in the classical body volume force $f^{*}$, but $f^{*}$ is now dependent on the elastic coefficients. 
Substituting back the relations (\ref{eq:equi_RM_Lc_0_dev_sym_3}) in the equation~(\ref{eq:equi_RM_Lc_0_1}) while also applying the ``dev sym'', and ``tr'' decomposition, allows us to write
\begin{align}
\text{Div}
\big[
2\mu_{\rm e} \, \text{dev sym} 
\left(\text{D}u - \left(\frac{\mu_{\rm e}}{\mu_{\rm e}+\mu_{\rm micro}}\text{D}u\right) \right)
+ \kappa_{\rm e} \, \text{tr}
\left(\text{D}u - \left(\frac{\kappa_{\rm e}}{\kappa_{\rm e}+\kappa_{\rm micro}} \, \text{D}u\right) \right) \mathbbm{1}
\big]
&=
f^{*}
\, ,
\notag
\\*[5pt]
\Longleftrightarrow
\qquad
\text{Div}
\big[
2\dfrac{\mu_{\rm e} \, \mu_{\rm micro}}{\mu_{\rm e} + \mu_{\rm micro}} \, \text{dev sym} 
\text{D}u
+ \dfrac{\kappa_{\rm e} \, \kappa_{\rm micro}}{\kappa_{\rm e} + \kappa_{\rm micro}} \, \text{tr}
\left(\text{D}u \right) \mathbbm{1}
\big]
&=
f^{*}
\, ,
\label{eq:equi_RM_Lc_0_3}
\\*[5pt]
\Longleftrightarrow
\qquad
\text{Div}
\big[
2\mu_{\rm macro} \, \text{dev sym} 
\text{D}u
+ \kappa_{\rm macro} \, \text{tr}
\left(\text{D}u \right) \mathbbm{1}
\big]
&=
f^{*}
\, .
\notag
\end{align}
where $f^*$ is defined as
\begin{align}
f^{*} \coloneqq 
f -
\text{Div}\big[
\frac{\mu_{\rm macro}}{\mu_{\rm micro}} \, \text{dev sym} \, M
+
\text{skew} \, M
+
\frac{1}{3}\frac{\kappa_{\rm macro}}{\kappa_{\rm micro}} \, \text{tr} (M) \mathbbm{1}
\big]
\, .
\end{align}
It is noted that $f^{*}$ depends on $\text{skew} \, M$ without any multiplicative elastic coefficient. This limit with a concentrated double body force may be instrumental in order to identify the \textit{micro} parameters. 
The equation (\ref{eq:equi_RM_Lc_0_3})$_3$ is the equilibrium equation for a classical isotropic linear elastic Cauchy continuum with stiffness $\mu_{\rm macro}$ and $\kappa_{\rm macro}$.
The relations for the macroscopic Lamé parameters ($\mu_{\rm macro},\lambda_{\rm macro}$) and the macroscopic bulk modulus ($\kappa_{\rm macro}$) are then
\begin{equation}
\begin{gathered}
\mu_{\rm macro} \coloneqq \dfrac{\mu_{\rm e} \, \mu_{\rm micro}}{\mu_{\rm e} + \mu_{\rm micro}} \, ,
\qquad\qquad\qquad
\kappa_{\rm macro} \coloneqq \dfrac{\kappa_{\rm e} \, \kappa_{\rm micro}}{\kappa_{\rm e} + \kappa_{\rm micro}} \, ,
\label{eq:Cmacro}
\\*[5pt]
\lambda_{\rm macro} \coloneqq \dfrac{1}{3}\left( 3 \kappa_{\rm macro}-2 \mu_{\rm macro} \right)\, \enskip  (\text{3D medium}) \, ,
\end{gathered}
\end{equation}
where  $\kappa_{\rm macro}$ is the macroscopic bulk modulus.
Relations (\ref{eq:Cmacro}) are the specialization of relation (\ref{eq:Ce_Cmac}) to the isotropic case \cite{barbagallo2017transparent}.
In order to have $\lambda_{\rm macro}=\lambda_{\rm micro}=0$, the only possible condition is $\lambda_{\rm micro}=\lambda_{\rm e}=0$.
Note that the macroscopic stiffness $\mathbb{C}_{\rm macro}$ (here $\mu_{\rm macro},\kappa_{\rm macro}$) is uniquely identified from classical homogenization techniques \cite{sarhil2022size}.
%
%
%
\subsubsection{Limit for $L_{\rm c}\to \infty$: upper bound microscopic stiffness $\mathbb{C}_{\rm micro}$}
The minimization of an energy functional that incorporates $\mu_{\rm M}L_{\rm c}^2 \, \| \text{Curl} P \|^2$, for the limit $L_{\rm c} \to \infty$, requires $\text{Curl} P=0$, and this implies that the micro-distortion tensor $P$ has to reduce to a gradient field $P \to \text{D}v$ on a simply connected domain such that
\begin{align}
\text{Curl} \, \text{D}v = 0 \qquad \forall \, v \in [\mathit{C}^\infty(\Omega)]^3 \, ,
\end{align}
thus asserting finite energies of the relaxed micromorphic model for arbitrarily large characteristic length values $L_{\rm c}$.
The corresponding strain energy density in terms of the reduced kinematics $\{ u, v \} : \Omega\to \mathbb{R}^{3}$ now reads
\begin{align}
W \left(\text{D}u, \text{D}v\right) = &
\, \mu_{\rm e} \left\lVert \text{sym} (\text{D}u - \text{D}v ) \right\rVert^{2}
+ \mu_{\rm c} \left\lVert \text{skew} (\text{D}u - \text{D}v ) \right\rVert^{2}
+ \frac{\lambda_{\rm e}}{2} \text{tr}^2 (\text{D}u - \text{D}v )
\label{eq:energy_RM_Lc_infty}
+ \mu_{\rm micro} \left\lVert \text{sym}\,\text{D}v \right\rVert^{2}
+ \frac{\lambda_{\rm micro}}{2} \text{tr}^2(\text{D}v)
\, .
\end{align}
The first variation of the strain energy $I=\displaystyle\int_{\Omega} W \, \text{d}x$ with respect to the two independent vector fields $u$ and $v$ leads to
\begin{align}
\delta I^{u}
&=
\int_{\Omega}
\hspace{-0.1cm}
\Big(
2\mu_{\rm e} \,
\langle \text{sym} (\text{D}u - \text{D}v ),\text{D}\delta u \rangle
+ 2\mu_{\rm c} \,
\langle \text{skew} (\text{D}u - \text{D}v ),\text{D}\delta u \rangle
+ \lambda_{\rm e} \langle \text{tr} (\text{D}u - \text{D}v ) \mathbbm{1}, \text{D}\delta u \rangle
+\langle f , \delta u \rangle
\Big)
\text{d}x
\, ,
\notag
\\*[5pt]
\delta I^{v}
&=
\int_{\Omega}
\hspace{-0.1cm}
\Big(
-2\mu_{\rm e} \,
\langle \text{sym} (\text{D}u - \text{D}v ),\text{D}\delta v \rangle
-2\mu_{\rm c} \,
\langle \text{skew} (\text{D}u - \text{D}v ),\text{D}\delta v \rangle
-\lambda_{\rm e} \langle \text{tr} (\text{D}u - \text{D}v ) \mathbbm{1}, \text{D}\delta v \rangle
\label{eq:first_varia_energy_RM_v_Lc_infy}
\\*
&
\phantom{=\int_{\Omega}\Big(}
+
2\mu_{\rm micro} \langle \text{sym}\,P,\text{D}\delta v \rangle
+ \lambda_{\rm micro} \langle \text{tr} (\text{D}v ) \mathbbm{1}, \text{D}\delta v \rangle
\Big)
\text{d}x
\, .
\notag
\end{align}
The equilibrium equations are now obtained by requiring 
\begin{equation}
\delta I^{u}= \langle f , \delta u \rangle \, ,
\quad
\forall
\,
\delta u
\qquad\qquad
\text{and}
\qquad\qquad
\delta I^{v}=
\langle M , \text{D}\delta v \rangle \, ,
\quad
\forall
\,
\delta v  \, .
\label{eq:first_variation_Lc_infy}
\end{equation}
where the contributions on the right sides are the virtual work of the external forces $f$ (classical body force) and $M$ (non-symmetric second order double body force tensor), and the equilibrium equations read
\begin{align}
\text{Div}
\big[
2\mu_{\rm e}\,\text{sym}  (\text{D}u - \text{D}v )
+ 2\mu_{\rm c}\,\text{skew} (\text{D}u - \text{D}v )
+ \lambda_{\rm e} \, \text{tr} (\text{D}u - \text{D}v ) \mathbbm{1}
\big]
&=
f
\, ,
\label{eq:equi_RM_Lc_infy_2_1}
\\*
-\text{Div}
\big[
2\mu_{\rm e}\,\text{sym}  (\text{D}u - \text{D}v )
+ 2\mu_{\rm c}\,\text{skew} (\text{D}u - \text{D}v )
+ \lambda_{\rm e} \,\text{tr} (\text{D}u - \text{D}v ) \mathbbm{1}
\big]
\quad
\notag
\\*
+
\text{Div}
\big[
2\mu_{\rm micro}\,\text{sym} \, \text{D}v
+ \lambda_{\rm micro} \,\text{tr} ( \text{D}v ) \mathbbm{1}
\big]
&=
\text{Div} \, M
\, ,
\label{eq:equi_RM_Lc_infy_2_2}
\end{align}
where the constraint $M \, n =0$ is required on the boundary, with $n$ the normal to the boundary.
The term on the left-hand side of equation (\ref{eq:equi_RM_Lc_infy_2_2}) can be substituted with the right-hand side of (\ref{eq:equi_RM_Lc_infy_2_1}) and, while keeping the equation (\ref{eq:equi_RM_Lc_infy_2_1}), we can re-write the system of equations (\ref{eq:equi_RM_Lc_infy_2_1})-(\ref{eq:equi_RM_Lc_infy_2_2}) as
\begin{align}
\text{Div}
\big[
2\mu_{\rm e}\,\text{sym}  (\text{D}u - \text{D}v )
+ 2\mu_{\rm c}\,\text{skew} (\text{D}u - \text{D}v )
+ \lambda_{\rm e} \,\text{tr} (\text{D}u - \text{D}v ) \mathbbm{1}
\big]
&=
f
\, ,
\label{eq:equi_RM_Lc_infy_3}
\\*
\text{Div}
\big[
2\mu_{\rm micro}\,\text{sym} \, \text{D}v
+ \lambda_{\rm micro} \text{tr} \, ( \text{D}v ) \mathbbm{1}
\big]
&=
f + \text{Div} \, M
\, ,
\notag
\end{align}
The only case in which $v=u$ is an admissible solution is if the classical body forces $f$ are zero.
In this case equations (\ref{eq:equi_RM_Lc_infy_3}) reduces to
\begin{align}
\text{Div} \, \sigma_{\rm micro}
=
\text{Div}
\big[
2\mu_{\rm micro}\,\text{sym} \, \text{D}u
+ \lambda_{\rm micro} \, \text{tr} ( \text{D}u ) \mathbbm{1}
\big]
=
\text{Div} \, M
\, ,
\label{eq:equi_RM_Lc_infy_4}
\end{align}
which is an equilibrium equation of the classical elasticity type with a microscopic stiffness given by $\mu_{\rm micro}$ and $\lambda_{\rm micro}$ and a body force vector equal to $\text{Div} \, M$.
%
%
%
\subsubsection{Limit for $\mathbb{C}_{\rm e} \to +\infty$ with $\mu_{\rm c}=0$: lower bound macroscopic stiffness $\mathbb{C}_{\rm macro}$ }
Due to the relations (\ref{eq:Cmacro}) we have formally $\mathbb{C}_{\rm micro}=\mathbb{C}_{\rm macro}$ as $\mathbb{C}_{\rm e} \to +\infty$.
The strain energy density  (\ref{eq:energy_RM_3}) is again reported here
\begin{align}
W \left(\text{D}u, P,\text{Curl}\,P\right) = &
\, \mu_{\rm e} \left\lVert \text{sym} (\text{D}u - P ) \right\rVert^{2}
+ \mu_{\rm c} \left\lVert \text{skew} (\text{D}u - P ) \right\rVert^{2}
+ \frac{\lambda_{\rm e}}{2} \text{tr}^2 (\text{D}u - P )
+ \mu_{\rm m} \left\lVert \text{sym}\,P \right\rVert^{2}
+ \frac{\lambda_{\rm m}}{2} \text{tr}^2 \left(P \right)
\notag
\\*
\label{eq:energy_RM_3_ce_infty}
&
+ \frac{\mu_{\rm M}\, L_{\rm c}^2}{2} \,
\left(
a_1 \, \left\lVert \text{dev sym} \, \text{Curl} \, P\right\rVert^2 +
a_2 \, \left\lVert \text{skew} \,  \text{Curl} \, P\right\rVert^2 +
\frac{a_3}{3} \, \text{tr}^2 \left(\text{Curl} \, P\right)
\right)
\to \text{min}(u,P) \,.
\end{align}
As $\mu_{\rm e},\lambda_{\rm e} \to \infty$, in order to remain with a bounded strain energy density, it is required that ${\text{sym}\,P=\text{sym}\,\text{D}u}$.
This, and $\mu_{\rm c}=0$, reduces the variational problem to 
\begin{align}
&\int_{\Omega}
\mu_{\rm m} \lVert \text{sym}\,\text{D}u \rVert^2 + \frac{\lambda_{\rm m}}{2}\text{tr}^2 (\text{sym}\,\text{D}u)
\label{eq:ene_lim_Ce_0}
\\*
&\phantom{\int_{\Omega}{}}
+ \frac{\mu_{\rm M}\, L_{\rm c}^2}{2} \,
\left(
a_1 \, \left\lVert \text{dev sym} \, \text{Curl} \, P\right\rVert^2 +
a_2 \, \left\lVert \text{skew} \,  \text{Curl} \, P\right\rVert^2 +
\frac{a_3}{3} \, \text{tr}^2 \left(\text{Curl} \, P\right)
\right)
\text{d}x
\quad\longrightarrow\text{min}(u,P).
\notag
\end{align}
The curvature part $\frac{\mu_{\rm M}\, L_{\rm c}^2}{2} \,\left(
a_1 \, \left\lVert \text{dev sym} \, \text{Curl} \, P\right\rVert^2 +
a_2 \, \left\lVert \text{skew} \,  \text{Curl} \, P\right\rVert^2 +
\frac{a_3}{3} \, \text{tr}^2 \left(\text{Curl} \, P\right)
\right)$
can be annihilated by choosing $\text{Curl}\,P=0$ which implies
\begin{align}
P=\text{D}\eta
\end{align}
on a simply connected domain.
Moreover, the remaining minimization in (\ref{eq:ene_lim_Ce_0}), using the consistent coupling condition delivers the unique solution $u$. Gathering, we have
\begin{align}
\text{sym}\,\text{D}u=\text{sym}\,\text{D}\eta
&
\quad\Longleftrightarrow\quad
\text{sym}(\text{D}(u-\eta))=0
\quad\Longleftrightarrow\quad
\text{D}(u-\eta)=A(x) \, , \quad A\in \mathfrak{so}(3)
\notag
\\*
&
\quad\Longrightarrow\quad
0=\text{Curl}\,\text{D}(u-\eta)=\text{Curl}A(x)
\notag
\\*
&
\quad\Longrightarrow\quad
A(x)=\overline{A}\quad \text{``rigidity''}\text{\cite{neff2008curl}}
\\*
&
\phantom{\quad\Longleftrightarrow\quad\,\,}
\text{D}u(x)-\text{D}\eta(x)=\overline{A} \in \mathfrak{so}(3)
\quad\Longrightarrow\quad
P=\text{D}\eta=\text{D}u-\overline{A}\quad \text{and} \quad \text{Curl}\,P=0
\,.
\notag
\end{align}
This leads to 
\begin{align}
I(u,P) & = \int_{\Omega}
\mu_{\rm M} \left\lVert \text{sym}\,(\text{D}u-\overline{A}) \right\rVert^{2}
+ \frac{\lambda_{\rm M}}{2} \text{tr}^2 \left(\text{D}u-\overline{A} \right)
+0
\, \text{d}x
\label{eq:energy_RM_3_ce_infy_2}
\\*
&= \int_{\Omega}
\mu_{\rm M} \left\lVert \text{sym}\,\text{D}u
\right\rVert^{2}
+ \frac{\lambda_{\rm M}}{2} \text{tr}^2 \left(\text{D}u \right)
\, \text{d}x
\qquad \to \text{min}\,u
\,.
\notag
\end{align}
Therefore $\mathbb{C}_{\rm e}\to +\infty$ gives size-independent linear elasticity with stiffness $\mathbb{C}_{\rm macro}$, as expected.
Note that, in contrast, the same limit of $\mathbb{C}_{\rm e}\to +\infty$ would lead to a gradient elasticity formulation for the classical Eringen-Mindlin micromorphic model \cite{d2022consistent} .
%
%

%
\subsection{A true two-scale model: the relaxed micromorphic model limit for $L_{\rm c}\to 0$ and $L_{\rm c}\to \infty$ in plane strain}
The relaxed micromorphic model reduces to a classical Cauchy model for both $L_{\rm c}\to 0$ and $L_{\rm c}\to \infty$ but with two different stiffnesses, $\mathbb{C}_{\rm macro}$ and $\mathbb{C}_{\rm micro}$, respectively.
The expressions of such stiffnesses are presented in the next two sections for the plane strain problem.
%
%
%
\subsubsection{Limit for $L_{\rm c}\to 0$: lower bound macroscopic stiffness $\mathbb{C}_{\rm macro}$}
For the limit $L_{\rm c} \to 0$, the equilibrium equations (\ref{eq:equi_RM_2D}) reduce to
\begin{align}
\text{Div}
\big[
2\mu_{\rm e}\,\text{sym}  (\text{D}\widetilde{u}^{\sharp} - \widetilde{P}^{\sharp} )
+ 2\mu_{\rm c}\,\text{skew} (\text{D}\widetilde{u}^{\sharp} - \widetilde{P}^{\sharp} )
+ \lambda_{\rm e} \text{tr} (\text{D}\widetilde{u}^{\sharp} - \widetilde{P}^{\sharp} ) \mathbbm{1}_2
\big]
&=
\widetilde{f}
\, ,
\label{eq:equi_RM_2D_Lc_0}
\\*
2\mu_{\rm e}\,\text{sym}  (\text{D}\widetilde{u}^{\sharp} - \widetilde{P}^{\sharp} )
+ 2\mu_{\rm c}\,\text{skew} (\text{D}\widetilde{u}^{\sharp} - \widetilde{P}^{\sharp} )
+ \lambda_{\rm e} \text{tr} (\text{D}\widetilde{u}^{\sharp} - \widetilde{P}^{\sharp} ) \mathbbm{1}_2
- 2 \mu_{\rm m}\,\text{sym}\,\widetilde{P}^{\sharp}
- \lambda_{\rm m} \text{tr} ( \widetilde{P}^{\sharp} ) \mathbbm{1}_2
&=
\widetilde{M}
\, .
\notag
\end{align}
The equation (\ref{eq:equi_RM_2D_Lc_0})$_2$ is now algebraic in $\widetilde{P}^{\sharp}$.
Due to the orthogonality of the ``sym/skew'' decomposition, the equation (\ref{eq:equi_RM_2D_Lc_0})$_2$ requires that
\begin{align}
2\mu_{\rm c} \, \text{skew} (\text{D}\widetilde{u}^{\sharp} - \widetilde{P}^{\sharp} ) 
&=
\text{sym} \, \widetilde{M}
\, ,
\\
2\mu_{\rm e}\,\text{sym}  (\text{D}\widetilde{u}^{\sharp} - \widetilde{P}^{\sharp} )
+ \lambda_{\rm e} \text{tr} (\text{D}\widetilde{u}^{\sharp} - \widetilde{P}^{\sharp} ) \mathbbm{1}_2
- 2 \mu_{\rm m}\,\text{sym}\,\widetilde{P}^{\sharp}
- \lambda_{\rm m} \text{tr} ( \widetilde{P}^{\sharp} ) \mathbbm{1}_2
&=
\text{skew} \, \widetilde{M}
\, .
\notag
\end{align}
Since the ``sym'' operator is not orthogonal to the ``tr'' operator, we further decompose ``sym'' into ``dev sym'' and ``tr sym'' so that
\begin{align}
2\mu_{\rm c} \, \text{skew} (\text{D}\widetilde{u}^{\sharp} - \widetilde{P}^{\sharp} ) 
&=
\text{skew} \, \widetilde{M}
\, ,
\quad
\notag
\\
2\mu_{\rm e}\,\text{dev$_2$ sym}  (\text{D}\widetilde{u}^{\sharp} - \widetilde{P}^{\sharp} )
+ \mu_{\rm e}\,\text{tr}  (\text{D}\widetilde{u}^{\sharp} - \widetilde{P}^{\sharp} ) \mathbbm{1}_2
+ \lambda_{\rm e} \text{tr} (\text{D}\widetilde{u}^{\sharp} - \widetilde{P}^{\sharp} ) \mathbbm{1}_2
\quad
&
\label{eq:equi_RM_2D_Lc_0_dev_sym}
\\*
- 2 \mu_{\rm m}\,\text{dev$_2$ sym}\,\widetilde{P}^{\sharp}
- \mu_{\rm m}\,\text{tr}\,(\widetilde{P}^{\sharp}) \mathbbm{1}_2
- \lambda_{\rm m} \text{tr} ( \widetilde{P}^{\sharp} ) \mathbbm{1}_2
&=
\text{sym} \, \widetilde{M}
\, .
\notag
\end{align}
note that ``tr sym'' is the same as ``tr''.
We also recall here the definition of the volumetric part, the deviatoric part, and the skew-symmetric parts in plane strain case
\begin{align}
\text{2D volumetric part}
&\coloneqq
\frac{1}{2} \text{tr}(\widetilde{P}^{\sharp})\mathbbm{1}_2 \, , \quad \mathbbm{1}_2=\begin{pmatrix} 1 & 0 \\ 0 & 1 \end{pmatrix} , \notag
\\*
\text{2D deviatoric symmetric part}
&\coloneqq
\frac{\widetilde{P}^{\sharp}+\widetilde{P}^{\sharp^{\rm T}}}{2}-\frac{1}{2} \text{tr}(\widetilde{P}^{\sharp})\mathbbm{1}_2=\text{dev$_2$ sym}\, \widetilde{P}^{\sharp} \, ,\notag
\\*
\label{skewP}
\text{2D skew symmetric part}
&\coloneqq
\frac{\widetilde{P}^{\sharp}-\widetilde{P}^{\sharp^{\rm T}}}{2} \, .
\end{align}
With further manipulations and due to the orthogonality of the operator ``skew'', ``dev sym'', and ``tr'', the system (\ref{eq:equi_RM_2D_Lc_0_dev_sym}) requires that
\begin{align}
2\mu_{\rm c} \, \text{skew} (\text{D}\widetilde{u}^{\sharp} - \widetilde{P}^{\sharp} ) 
&=
\text{skew} \, \widetilde{M}
\, ,
\notag
\\
\mu_{\rm e}\,\text{dev$_2$ sym}  (\text{D}\widetilde{u}^{\sharp} - \widetilde{P}^{\sharp} )
- \mu_{\rm m}\,\text{dev$_2$ sym}\,\widetilde{P}^{\sharp}
&=
\text{dev sym} \, \widetilde{M}
\label{eq:equi_RM_2D_Lc_0_dev_sym_2}
\, ,
\\*
\left(\mu_{\rm e}+ \lambda_{\rm e}\right)\,\text{tr}  (\text{D}\widetilde{u}^{\sharp} - \widetilde{P}^{\sharp} ) \mathbbm{1}_2
- \left(\mu_{\rm m} + \lambda_{\rm m}\right)\,\text{tr}\,(\widetilde{P}^{\sharp}) \mathbbm{1}_2
&=
\frac{1}{2}\text{tr} (\widetilde{M}) \mathbbm{1}_2
\, .
\notag
\end{align}
From equation (\ref{eq:equi_RM_2D_Lc_0_dev_sym_2}) we can evaluate the expressions for $\text{skew} \, \widetilde{P}^{\sharp}$, $\text{dev sym} \, \widetilde{P}^{\sharp}$, and $\text{tr}(\widetilde{P}^{\sharp})$ as
\begin{align}
\text{skew} \, \text{D}\widetilde{u}^{\sharp}
-
\frac{1}{2 \mu_{\rm c}}\text{skew}\,\widetilde{M}
&=
\text{skew} \, \widetilde{P}^{\sharp}
\, ,
\notag
\\*
\frac{\mu_{\rm e}}{\mu_{\rm e}+\mu_{\rm m}}
\text{dev$_2$ sym}  \, \text{D}\widetilde{u}^{\sharp}
-
\frac{1}{2 (\mu_{\rm e}+\mu_{\rm m})}\text{dev$_2$ sym}\,\widetilde{M}
&=
\text{dev$_2$ sym}  \, \widetilde{P}^{\sharp}
\, ,
\label{eq:equi_RM_2D_Lc_0_dev_sym_3}
\\*
\frac{\kappa_{\rm e}}{\kappa_{\rm e}+\kappa_{\rm m}}
\text{tr}  \, \text{D}\widetilde{u}^{\sharp}
-
\frac{1}{2 (\kappa_{\rm e}+\kappa_{\rm m})}\text{tr}\,\widetilde{M}
&=
\text{tr}(\widetilde{P}^{\sharp})
\, .
\notag
\end{align}
where  $\kappa_{\rm e}=\mu_{\rm e}+\lambda_{\rm e}$ and $\kappa_{\rm m}=\mu_{\rm m}+\lambda_{\rm m}$ are the plane strain bulk moduli.

Substituting back the relations (\ref{eq:equi_RM_2D_Lc_0_dev_sym_3}) in the equation~(\ref{eq:equi_RM_2D_Lc_0})$_1$ while also applying the ``dev sym'', and ``tr'' decomposition, we have
\begin{align}
\text{Div}
\big[
2\mu_{\rm e} \, \text{dev sym} 
\left(\text{D}\widetilde{u}^{\sharp} - \left(\frac{\mu_{\rm e}}{\mu_{\rm e}+\mu_{\rm m}}\text{D}\widetilde{u}^{\sharp}\right) \right)
+ \widetilde{\kappa}_{\rm e} \, \text{tr}
\left(\text{D}\widetilde{u}^{\sharp} - \left(\frac{\widetilde{\kappa}_{\rm e}}{\widetilde{\kappa}_{\rm e}+\widetilde{\kappa}_{\rm m}} \, \text{D}\widetilde{u}^{\sharp}\right) \right) \mathbbm{1}_2
\big]
&=
\widetilde{f}^{*}
\, ,
\notag
\\*[5pt]
\Longleftrightarrow
\qquad
\text{Div}
\big[
2\dfrac{\mu_{\rm e} \, \mu_{\rm m}}{\mu_{\rm e} + \mu_{\rm m}} \, \text{dev$_2$ sym} \,\text{D}\widetilde{u}^{\sharp}
+ \dfrac{\widetilde{\kappa}_{\rm e} \, \widetilde{\kappa}_{\rm m}}{\widetilde{\kappa}_{\rm e} + \widetilde{\kappa}_{\rm m}} \, \text{tr}
\left(\text{D}\widetilde{u}^{\sharp} \right) \mathbbm{1}_2
\big]
&=
\widetilde{f}^{*}
\, ,
\label{eq:equi_RM_2D_Lc_0_2}
\\*[5pt]
\Longleftrightarrow
\qquad
\text{Div}
\big[
2\mu_{\rm M} \, \text{dev$_2$ sym}\,\text{D}\widetilde{u}^{\sharp}
+ \widetilde{\kappa}_{\rm M} \, \text{tr}
(\text{D}\widetilde{u}^{\sharp}) \mathbbm{1}_2
\big]
&=
\widetilde{f}^{*}
\, .
\notag
\end{align}
where $\widetilde{f}^*$ is defined as
\begin{align}
\widetilde{f}^{*} \coloneqq
\widetilde{f} -
\text{Div}\big[
\frac{\mu_{\rm M}}{\mu_{\rm m}} \, \text{dev$_2$ sym} \, \widetilde{M}
+
\text{skew} \, \widetilde{M}
+
\frac{1}{2}\frac{\kappa_{\rm M}}{\kappa}_{\rm m} \, \text{tr} (\widetilde{M}) \mathbbm{1}_2
\big]
\, .
\end{align}
It is noted that $\widetilde{f}^{*}$ depends on $\text{skew} \, \widetilde{M}$ without any multiplicative elastic coefficient because of the choice of an isotropic constitutive law (an isotropic second order skew-symmetric tensor depends on one coefficient). This limit with a concentrated double body force may be instrumental in order to identify the \textit{micro} parameters. The equation (\ref{eq:equi_RM_2D_Lc_0_2})$_3$ is the equilibrium equation for a classical linear elastic isotropic Cauchy continuum with stiffness $\mu_{\rm macro}$ and $\kappa_{\rm macro}$. The relations for the macro Lamé parameters ($\mu_{\rm M},\lambda_{\rm M}$) and the macroscopic bulk modulus for plane strain are given in \eqref{defmod}. Note that in order to have $\lambda_{\rm macro}=\lambda_{\rm micro}=0$, the only possible condition is again $\lambda_{\rm micro}=\lambda_{\rm e}=0$.
%
%
%
\subsubsection{Limit for $L_{\rm c}\to \infty$: upper bound microscopic stiffness $\mathbb{C}_{\rm micro}$}
The minimization of an energy functional that incorporate $\mu_{\rm M}\, L_{\rm c}^2 \, \| \text{Curl} \widetilde{P}^{\sharp} \|^2$, for the limit $L_{\rm c} \to \infty$, requires $\text{Curl} \widetilde{P}^{\sharp}=0$, and this implies that the micro-distortion tensor $P$ has to reduce to a gradient field $\widetilde{P}^{\sharp} \to \text{D}\widetilde{v}^{\sharp}$ on a simply connected domain and
\begin{align}
\text{Curl} \, \text{D}\widetilde{v}^{\sharp} = 0 \qquad \forall \, \widetilde{v}^{\sharp} \in [\mathit{C}^\infty(\Omega)]^3 \, ,
\end{align}
thus asserting finite energies of the relaxed micromorphic model for arbitrarily large characteristic length values $L_{\rm c}$.
The corresponding strain energy density in terms of the reduced kinematics $\{ \widetilde{u}, \widetilde{v}^{\sharp} \} : \Omega\to \mathbb{R}^{3}$ now reads
\begin{align}
W \left(\text{D}\widetilde{u}, \text{D}\widetilde{v}^{\sharp}\right) = &
\, \mu_{\rm e} \left\lVert \text{sym} (\text{D}\widetilde{u}^{\sharp} - \text{D}\widetilde{v}^{\sharp} ) \right\rVert^{2}
+ \mu_{\rm c} \left\lVert \text{skew} (\text{D}\widetilde{u}^{\sharp} - \text{D}\widetilde{v}^{\sharp} ) \right\rVert^{2}
+ \frac{\lambda_{\rm e}}{2} \text{tr}^2 (\text{D}\widetilde{u}^{\sharp} - \text{D}\widetilde{v}^{\sharp} )
\label{eq:energy_RM_2D_Lc_infty}
\\*
&
+ \mu_{\rm m} \left\lVert \text{sym}\,\text{D}\widetilde{v}^{\sharp} \right\rVert^{2}
+ \frac{\lambda_{\rm m}}{2} \text{tr}^2 \left(\text{D}\widetilde{v}^{\sharp} \right)
\, .
\notag
\end{align}
The first variation of the strain energy $I=\displaystyle\int_{\Omega} W \, \text{d}x$ with respect to the two independent vector fields $\widetilde{u}^{\sharp}$ and $\widetilde{v}^{\sharp}$ leads to
\begin{align}
\delta I^{\widetilde{u}}
&=
\int_{\Omega}
\hspace{-0.1cm}
\Big(
2\mu_{\rm e} \,
\langle \text{sym} (\text{D}\widetilde{u}^{\sharp} - \text{D}\widetilde{v}^{\sharp} ),\text{D}\delta \widetilde{u}^{\sharp} \rangle
+ 2\mu_{\rm c} \,
\langle \text{skew} (\text{D}\widetilde{u}^{\sharp} - \text{D}\widetilde{v}^{\sharp} ),\text{D}\delta \widetilde{u}^{\sharp} \rangle
\label{eq:first_varia_energy_RM_2D_u_Lc_infy}
+ \lambda_{\rm e} \langle \text{tr} (\text{D}\widetilde{u}^{\sharp} - \text{D}\widetilde{v}^{\sharp} ) \mathbbm{1}_2, \text{D}\delta \widetilde{u}^{\sharp} \rangle
\Big)
\text{d}x
\, ,
\\*[5pt]
\delta I^{\widetilde{v}^{\sharp}}
&=
\int_{\Omega}
\hspace{-0.1cm}
\Big(
-2\mu_{\rm e} \,
\langle \text{sym} (\text{D}\widetilde{u}^{\sharp} - \text{D}\widetilde{v}^{\sharp} ),\text{D}\delta \widetilde{v}^{\sharp} \rangle
-2\mu_{\rm c} \,
\langle \text{skew} (\text{D}\widetilde{u}^{\sharp} - \text{D}\widetilde{v}^{\sharp} ),\text{D}\delta \widetilde{v}^{\sharp} \rangle
-\lambda_{\rm e} \langle \text{tr} (\text{D}\widetilde{u}^{\sharp} - \text{D}\widetilde{v}^{\sharp} ) \mathbbm{1}_2, \text{D}\delta \widetilde{v}^{\sharp} \rangle
\label{eq:first_varia_energy_RM_2D_v_Lc_infy}
\\*
&
\phantom{=\int_{\Omega}\Big(}
+
2\mu_{\rm m} \langle \text{sym}\,\text{D}\widetilde{v}^{\sharp},\text{D}\delta \widetilde{v}^{\sharp} \rangle
+ \lambda_{\rm m} \langle \text{tr} (\text{D}\widetilde{v}^{\sharp} ) \mathbbm{1}_2, \text{D}\delta \widetilde{v}^{\sharp} \rangle
\Big)
\text{d}x
\, .
\notag
\end{align}
The equilibrium equations are now obtained by requiring 
\begin{equation}
\delta I^{\widetilde{u}^{\sharp}}= \langle \widetilde{f} , \delta \widetilde{u}^{\sharp} \rangle \, ,
\quad
\forall
\,
\delta \widetilde{u}^{\sharp}
\qquad\qquad
\text{and}
\qquad\qquad
\delta I^{\widetilde{v}^{\sharp}}=
\langle \widetilde{M} , \text{D}\delta \widetilde{v}^{\sharp} \rangle \, ,
\quad
\forall
\,
\delta \widetilde{v}^{\sharp}  \, .
\label{eq:first_variation_2D_Lc_infy}
\end{equation}
where the contributions on the right sides are the virtual work of the external forces $\widetilde{f}$ (classical body force) and $\widetilde{M}$ (non-symmetric second order double body force tensor), and the equilibrium equations read
\begin{align}
\text{Div}
\big[
2\mu_{\rm e}\,\text{sym}  (\text{D}\widetilde{u}^{\sharp} - \text{D}\widetilde{v}^{\sharp} )
+ 2\mu_{\rm c}\,\text{skew} (\text{D}\widetilde{u}^{\sharp} - \text{D}\widetilde{v}^{\sharp} )
+ \lambda_{\rm e} \, \text{tr} (\text{D}\widetilde{u}^{\sharp} - \text{D}\widetilde{v}^{\sharp} ) \mathbbm{1}_2
\big]
&=
\widetilde{f}
\, ,
\label{eq:equi_RM_2D_Lc_infy_2}
\\*
-\text{Div}
\big[
2\mu_{\rm e}\,\text{sym}  (\text{D}\widetilde{u}^{\sharp} - \text{D}\widetilde{v}^{\sharp} )
+ 2\mu_{\rm c}\,\text{skew} (\text{D}\widetilde{u}^{\sharp} - \text{D}\widetilde{v}^{\sharp} )
+ \lambda_{\rm e} \,\text{tr} (\text{D}\widetilde{u}^{\sharp} - \text{D}\widetilde{v}^{\sharp} ) \mathbbm{1}_2
\big]
\quad
\notag
\\*
+
\text{Div}
\big[
2\mu_{\rm m}\,\text{sym}  \text{D}\widetilde{v}^{\sharp}
+ \lambda_{\rm m} \,\text{tr} ( \text{D}\widetilde{v}^{\sharp} ) \mathbbm{1}_2
\big]
&=
\text{Div} \, \widetilde{M}
\, ,
\notag
\end{align}
where the constraint $\widetilde{M} \, n =0$ is required on the boundary, with $n$ the normal to the boundary.
The term on the left-hand side of equation (\ref{eq:equi_RM_2D_Lc_infy_2})$_2$ can be substituted with the right-hand side of (\ref{eq:equi_RM_2D_Lc_infy_2})$_1$ and, while keeping the equation (\ref{eq:equi_RM_2D_Lc_infy_2})$_1$, we can re-write the system of equations (\ref{eq:equi_RM_2D_Lc_infy_2}) as
\begin{align}
\text{Div}
\big[
2\mu_{\rm e}\,\text{sym}  (\text{D}\widetilde{u}^{\sharp} - \text{D}\widetilde{v}^{\sharp} )
+ 2\mu_{\rm c}\,\text{skew} (\text{D}\widetilde{u}^{\sharp} - \text{D}\widetilde{v}^{\sharp} )
+ \lambda_{\rm e} \,\text{tr} (\text{D}\widetilde{u}^{\sharp} - \text{D}\widetilde{v}^{\sharp} ) \mathbbm{1}_2
\big]
&=
\widetilde{f}
\, ,
\label{eq:equi_RM_2D_Lc_infy_3}
\\*
\text{Div}
\big[
2\mu_{\rm m}\,\text{sym} \, \text{D}\widetilde{v}^{\sharp}
+ \lambda_{\rm m} \text{tr} \, ( \text{D}\widetilde{v}^{\sharp} ) \mathbbm{1}_2
\big]
&=
\widetilde{f} + \text{Div} \, \widetilde{M}
\, .
\notag
\end{align}
The only case in which $\widetilde{v}^{\sharp}=\widetilde{u}^{\sharp}$ is an admissible solution is if the classical body forces $\widetilde{f}$ are zero.
In this case  (\ref{eq:equi_RM_2D_Lc_infy_3}) reduces to
\begin{align}
\text{Div} \, \sigma_{\rm m}
=
\text{Div}
\big[
2\mu_{\rm m}\,\text{sym} \, \text{D}\widetilde{u}^{\sharp}
+ \lambda_{\rm m} \, \text{tr} ( \text{D}\widetilde{u}^{\sharp} ) \mathbbm{1}_2
\big]
=
\text{Div} \, \widetilde{M}
\, ,
\label{eq:equi_RM_2D_Lc_infy_4}
\end{align}
which is an equilibrium equation of the classical elasticity type with a micro stiffness given by $\mu_{\rm m}$ and $\lambda_{\rm m}$ and a body force vector equal to $\text{Div} \, \widetilde{M}$.
%
%
%
%
%
%
\subsection{Some particular cases of the relaxed micromorphic model}
\label{sec:zero_coss}
\subsubsection{The pure relaxed micromorphic equations}
If we set $\mu_{\rm c}=0$, the force stress tensor $\sigma$ becomes symmetric and the model reduces to
\begin{align}
\text{Div}
\big[
\overbrace{
2\mu_{\rm e}\,\text{sym}  (\text{D}\widetilde{u}^{\sharp} - \widetilde{P}^{\sharp} )
+ \lambda_{\rm e} \text{tr} (\text{D}\widetilde{u}^{\sharp} - \widetilde{P}^{\sharp} ) \mathbbm{1}_2
}^{
\mathlarger{\sigma}\coloneqq
}
\big]
&=
\widetilde{f}
\, ,
\notag
\\*
\sigma
- 2 \mu_{\rm m}\,\text{sym}\,\widetilde{P}^{\sharp}
- \lambda_{\rm m} \text{tr} ( \widetilde{P}^{\sharp} ) \mathbbm{1}_2
-\mu_{\rm M}  L_{\rm c}^2 \, \widetilde{a}\,
\text{Curl} \, \text{Curl}_{\text{2D}} \, \widetilde{P}^{\sharp}
&=
\widetilde{M}
\, ,
\label{eq:equi_RM_2D_CC_mc0}
\end{align}
\begin{align*}
\widetilde{M}
&=
\left(
\begin{array}{ccccccc}
M_{11} & M_{12} & 0\\
M_{21} & M_{22} & 0\\
0                  & 0                  & 0
\end{array}
\right)
\, ,
\qquad\qquad\qquad
\widetilde{f}=
\left(
\begin{array}{ccccccc}
f_{1} \\
f_{2} \\
0                  
\end{array}
\right)
\, .
\end{align*}
In components we have
\begin{align}
(\lambda_{\rm e}+2 \mu_{\rm e}) \left(u_{1,11}-P_{11,1}\right)+\lambda_{\rm e} \left(u_{2,12}-P_{22,1}\right)+\mu_{\rm e} \left(-P_{12,2}-P_{21,2}+u_{1,22}+u_{2,12}\right)
&=
f_1 \, ,
\notag
\\
(\lambda_{\rm e}+2 \mu_{\rm e}) \left(u_{2,22}-P_{22,2}\right)+\lambda_{\rm e} \left(u_{1,12}-P_{11,2}\right)+\mu_{\rm e} \left(-P_{12,1}-P_{21,1}+u_{1,12}+u_{2,11}\right)
&=
f_2 \, ,
\notag
\\
\mu_{\rm M}\,L_{\rm c}^2 \widetilde{a} \left(P_{11,22}-P_{12,12}\right)-P_{11} (\lambda_{\rm e}+\lambda_{\rm m}+2 (\mu_{\rm e}+\mu_{\rm m}))-(\lambda_{\rm e}+\lambda_{\rm m}) P_{22}+(\lambda_{\rm e}+2 \mu_{\rm e}) u_{1,1}+\lambda_{\rm e} u_{2,2}
&=
M_{11} \, ,
\\
-\mu_{\rm M}\,L_{\rm c}^2 \widetilde{a} \left(P_{11,12}-P_{12,11}\right)-(\mu_{\rm e}+\mu_{\rm m}) P_{12}-(\mu_{\rm e}+\mu_{\rm m}) P_{21}+\mu_{\rm e} \left(u_{1,2}+u_{2,1}\right)
&=
M_{12} \, ,
\notag
\\
\mu_{\rm M}\,L_{\rm c}^2 \widetilde{a} \left(P_{21,22}-P_{22,12}\right)-(\mu_{\rm e}+\mu_{\rm m}) P_{12}-(\mu_{\rm e}+\mu_{\rm m}) P_{21}+\mu_{\rm e} \left(u_{1,2}+u_{2,1}\right)
&=
M_{21} \, ,
\notag
\\
-\mu_{\rm M}\,L_{\rm c}^2 \widetilde{a} \left(P_{21,12}-P_{22,11}\right)-P_{22} (\lambda_{\rm e}+\lambda_{\rm m}+2 (\mu_{\rm e}+\mu_{\rm m}))-(\lambda_{\rm e}+\lambda_{\rm m}) P_{11}+(\lambda_{\rm e}+2 \mu_{\rm e}) u_{2,2}+\lambda_{\rm e} u_{1,1}
&=
M_{22} \, .
\notag
\end{align}
%
%
%
\subsubsection{The relaxed micromorphic model with zero micro and macro Poisson's ratio}
\label{sec:zero_poisson}
If we set $\lambda_{\rm m}=\lambda_{\rm e}=0$, which implies $\lambda_{\rm M}=0$, the equilibrium equations (\ref{eq:equi_RM_2D}) reduce to 
\begin{align}
\text{Div}
\big[
2\mu_{\rm e}\,\text{sym}  (\text{D}\widetilde{u}^{\sharp} - \widetilde{P}^{\sharp} )
+ 2\mu_{\rm c}\,\text{skew} (\text{D}\widetilde{u}^{\sharp} - \widetilde{P}^{\sharp} )
\big]
&=
\widetilde{f}
\, ,
\label{eq:equi_RM_2D_CC_l0}
\\*
2\mu_{\rm e}\,\text{sym}  (\text{D}\widetilde{u}^{\sharp} - \widetilde{P}^{\sharp} )
+ 2\mu_{\rm c}\,\text{skew} (\text{D}\widetilde{u}^{\sharp} - \widetilde{P}^{\sharp} )
- 2 \mu_{\rm m}\,\text{sym}\,\widetilde{P}^{\sharp}
- \mu_{\rm M}\, L_{\rm c}^2 \, \widetilde{a} \, \text{Curl} \, \text{Curl}_{\text{2D}} \, \widetilde{P}^{\sharp}
&=
\widetilde{M}
\, .
\notag
\end{align}
Componentwise, we have
\begin{align}\begin{split}
\mu_{\rm c} \left(u_{1,22}-u_{2,12}+P_{21,2}-P_{12,2}\right)
+\mu_{\rm e} \left(u_{1,22}+2 u_{1,11}+u_{2,12}-2 P_{11,1}-P_{12,2}-P_{21,2}\right)
&=f_1 \, ,
\\*[5pt]
\mu_{\rm c} \left(P_{12,1}-P_{21,1}-u_{1,12}+u_{2,11}\right)
+\mu_{\rm e} \left(u_{1,12}+2 u_{2,22}+u_{2,11}-P_{12,1}-P_{21,1}-2 P_{22,2}\right)
&=f_2 \, ,
\\*[10pt]
\widetilde{a} \mu_{\rm M}\, L_{\rm c}^2 \left(P_{11,22}-P_{12,12}\right)
+2\mu_{\rm e} \left(u_{1,1}-P_{11}\right)
-2 \mu_{\rm m} P_{11}
&=M_{11} \, ,
\\*[5pt]
\widetilde{a} \mu_{\rm M}\, L_{\rm c}^2 \left(P_{12,11}-P_{11,12}\right)+\mu_{\rm c} \left(u_{1,2}-u_{2,1}-P_{12}+P_{21}\right)
+\mu_{\rm e} \left(u_{1,2}+u_{2,1}-P_{12}-P_{21}\right)
-\mu_{\rm m} (P_{12}+P_{21})
&=M_{12} \, ,
\\*[5pt]
\widetilde{a} \mu_{\rm M}\, L_{\rm c}^2 \left(P_{21,22}-P_{22,12}\right)
+\mu_{\rm c} \left(u_{2,1}-u_{1,2}+P_{12}-P_{21}\right)
+\mu_{\rm e} \left(u_{1,2}+u_{2,1}-P_{12}-P_{21}\right)
-\mu_{\rm m} (P_{12}+P_{21}) 
&=M_{21} \, ,
\\*[5pt]
\widetilde{a} \mu_{\rm M}\, L_{\rm c}^2 \left(P_{22,11}-P_{21,12}\right)
+2\mu_{\rm e} \left(u_{2,2}-P_{22}\right)
-2 \mu_{\rm m} P_{22}
&=M_{22} \, .\end{split}
\end{align}

The conditions for existence and uniqueness for the model in  (\ref{eq:equi_RM_2D_CC_l0}) are
\begin{align}
\mu_{\rm e}>0 \, ,
\qquad\qquad
\mu_{\rm m}>0 \, ,
\qquad\qquad
\mu_{\rm M}\, L_{\rm c}^2\widetilde{a}>0 \, ,
\qquad\qquad
\mu_{\rm c}\geq0 \, .
\end{align}
For $\mu_{\rm c}\equiv0$, in order to guarantee existence and uniqueness, one needs tangential boundary conditions for $\widetilde{P}$, while for $\mu_{\rm c}>0$, one does not need boundary conditions for $\widetilde{P}$ in order to guarantee existence and uniqueness.
%
%
%
%
%
%
\subsubsection{The relaxed micromorphic model with one curvature parameter, a zero Cosserat couple modulus, and a zero micro and macro Poisson's ratio}
If in addition to the simplifications of Sec.~\ref{sec:zero_poisson} we also set $\mu_{\rm c}=0$, the equilibrium equations (\ref{eq:equi_RM_2D_CC_l0}) further reduce to 
\begin{align}
\text{Div}
\big[
2\mu_{\rm e}\,\text{sym}  (\text{D}\widetilde{u}^{\sharp} - \widetilde{P}^{\sharp} )
\big]
&=
\widetilde{f}
\, ,
\label{eq:equi_RM_2D_CC_mc0_l0}&
2\mu_{\rm e}\,\text{sym}  (\text{D}\widetilde{u}^{\sharp} - \widetilde{P}^{\sharp} )
- 2 \mu_{\rm m}\,\text{sym}\,\widetilde{P}^{\sharp}
- \mu_{\rm M}\, L_{\rm c}^2 \, \widetilde{a} \, \text{Curl} \, \text{Curl}_{\text{2D}} \, \widetilde{P}^{\sharp}
=
\widetilde{M}
\, .
\end{align}
This represents the most simple set of equations for the plane strain relaxed micromorphic model.
In components we have
\begin{align}
\mu_{\rm e} \left(-2 P_{11,1}-P_{12,2}-P_{21,2}+u_{1,22}+2 u_{1,11}+u_{2,12}\right)
&=
f_1 \, ,
\notag
\\
\mu_{\rm e} \left(-P_{12,1}-P_{21,1}-2 P_{22,2}+u_{1,12}+2 u_{2,22}+u_{2,11}\right)
&=
f_2 \, ,
\notag
\\
\mu_{\rm M}\,L_{\rm c}^2 \widetilde{a} \left(P_{11,22}-P_{12,12}\right)-2 (\mu_{\rm e}+\mu_{\rm m}) P_{11}+2 \mu_{\rm e} u_{1,1}
&=
M_{11} \, ,
\\
- \mu_{\rm M}\,L_{\rm c}^2 \widetilde{a} \left(P_{11,12}-P_{12,11}\right)-(\mu_{\rm e}+\mu_{\rm m}) P_{12}-(\mu_{\rm e}+\mu_{\rm m}) P_{21}+\mu_{\rm e} \left(u_{1,2}+u_{2,1}\right)
&=
M_{12} \, ,
\notag
\\
\mu_{\rm M}\,L_{\rm c}^2 \widetilde{a} \left(P_{21,22}-P_{22,12}\right)-(\mu_{\rm e}+\mu_{\rm m}) P_{12}-(\mu_{\rm e}+\mu_{\rm m}) P_{21}+\mu_{\rm e} \left(u_{1,2}+u_{2,1}\right)
&=
M_{21} \, ,
\notag
\\
- \mu_{\rm M}\,L_{\rm c}^2 \widetilde{a} \left(P_{21,12}-P_{22,11}\right)-2 (\mu_{\rm e}+\mu_{\rm m}) P_{22}+2 \mu_{\rm e} u_{2,2}
&=
M_{22} \, .
\notag
\end{align}
\subsection{Subclasses of the relaxed micromorphic model as singular limits}
%
%
%

\subsubsection{The isotropic micro-stretch model in dislocation form as a particular case of the relaxed micromorphic model}
The micro-stretch model in dislocation format \cite{neff2014unifying,scalia2000extension,de1997torsion,neff2009mean,kirchner2007mechanics} can be obtained from the relaxed micromorphic model by letting formally $\mu_{\rm  micro}\to\infty$, while $\kappa_{\rm  micro}<\infty$.
For bounded energy, the micro-distortion tensor $P$ must be devoid  from the deviatoric component $\text{dev} \, \text{sym} \, P = 0 \Leftrightarrow P = A + \omega \mathbbm{1}$, $A \in \mathfrak{so}(3)$, $\omega \in \mathbb{R}$.
The expression of the strain energy for this model in dislocation format can then be written as \cite{neff2014unifying} (using Curl as the curvature measure)
\begin{align}
W \left(\text{D}u, A,\omega,\text{Curl}\left(A + \omega \mathbbm{1}\right)\right) 
=
&
\, \mu_{\rm M} \left\lVert \text{dev} \, \text{sym} \, \text{D}u \right\rVert^{2}
+ \frac{\kappa_{\text{e}}}{2} \text{tr}^2 \left(\text{D}u - \omega \mathbbm{1} \right) 
+ \mu_{c} \left\lVert \text{skew} \left(\text{D}u - A \right) \right\rVert^{2}
+ \frac{9}{2} \, \kappa_{\rm m} \, \omega^2
\label{eq:energy_MST}
\\*
&
+ \frac{\mu_{\rm M}\, L_{\rm c}^2}{2} \,
\left(
a_1 \, \left\lVert \text{dev sym} \, \text{Curl} \, A \right\rVert^2
+ a_2 \, \left\lVert \text{skew} \,  \text{Curl} \left(A + \omega \mathbbm{1}\right) \right\rVert^2
+ \frac{a_3}{3} \, \text{tr}^2 \left(\text{Curl} \, A \right)
\right) \, ,
\notag
\end{align}
since $\text{Curl} \left(\omega \mathbbm{1}\right) \in \mathfrak{so}(3)$.
The equilibrium equations, in the absence of body forces,   are obtained by variation of $(u, A, \omega)$ respectively and read
\begin{align}
\text{Div}\overbrace{\left[
2\mu_{\rm M}\,\text{dev}\,\text{sym} \, \text{D}u
+ \kappa_{\text{e}} \text{tr} \left(\text{D}u - \omega \mathbbm{1}\right) \mathbbm{1}
+ 2\mu_{c}\,\text{skew} \left(\text{D}u - A\right) \right]}^{\mathlarger{\widetilde{\sigma}}\coloneqq}
&= f \, ,
\notag
\\*
2\mu_{c}\,\text{skew} \left(\text{D}u - A\right)-\mu_{\rm M}\, L_{\rm c}^2 \, \text{skew} \, \text{Curl}\left(
a_1 \, \text{dev} \, \text{sym} \, \text{Curl} \, A \, 
+ a_2 \, \text{skew} \, \text{Curl} \left(A + \, \omega \mathbbm{1}\right) \,
+ \frac{a_3}{3} \, \text{tr} \left(\text{Curl} \, A \right)\mathbbm{1} \, 
\right) &=\text{skew}\,M \, ,
\label{eq:equi_MST}
\\*
\text{tr}
\bigg[
\kappa_{\text{e}} \text{tr} \left(\text{D}u - \omega \mathbbm{1}\right) \mathbbm{1}
- \kappa_{\rm micro} \text{tr} \left( \omega \mathbbm{1}\right) \mathbbm{1}
-\mu_{\rm M}\, L_{\rm c}^2 \,  a_2 \, \text{Curl}\,
\text{skew} \, \text{Curl} \left(\omega \mathbbm{1} + A\right) 
\bigg]
&= \text{tr}(M) \,.
\notag
\end{align}
Under the plane-strain hypothesis only the in-plane components of the kinematic fields are different from zero and they only depend on $(x_1,x_2)$.
The structure of the kinematic fields ($\widetilde{u}$,$\widetilde{A}$,$\omega$) are
\begin{align}
\widetilde{u}
&
=
\left(
\begin{array}{ccc}
u_{1} \\ 
u_{2} \\
0
\end{array}
\right)
\, ,
\qquad\quad
\widetilde{A}
=
\left(
\begin{array}{cc|c}
0 & A_{12} & 0 \\ 
-A_{22} & 0 & 0 \\ 
\hline
0 & 0 & 0
\end{array}
\right)
\, ,
\qquad\quad
\omega \widetilde{\mathbbm{1}_2}
=
\omega
\left(
\begin{array}{cc|c}
1 & 0 & 0 \\ 
0 & 1 & 0 \\ 
\hline
0 & 0 & 0
\end{array}
\right)
\, ,
\notag
\\*
\text{Curl}(\widetilde{A}+\omega\widetilde{\mathbbm{1}_2})
&=
\left(
\begin{array}{ccc}
0 & 0 & A_{12,1}-\omega_{,2} \\
0 & 0 & A_{12,2}+\omega_{,1} \\
0 & 0 & 0 \\
\end{array}
\right)
\, ,
\label{eq:gen_kine_stretch}
\\*
\text{Curl} \, \text{Curl}(\widetilde{A}+\omega\widetilde{\mathbbm{1}_2})
&=
\left(
\begin{array}{cc|c}
A_{12,12}-\omega_{,22} & \omega_{,12}-A_{12,11} & 0 \\
A_{12,11}+\omega_{,12} & -A_{12,12}-\omega_{,11} & 0 \\
\hline
0 & 0 & 0 \\
\end{array}
\right)
\,.
\notag
\end{align}
Under the plane-strain assumption, the equilibrium equations in components read now
\begin{align}
-2 \mu_{\rm c} A_{12,2}+(\kappa_{\rm e}+\mu_{\rm e}) u_{1,11}+\kappa_{\rm e} u_{2,12}-2 \kappa_{\rm e} \omega_{,1}+(\mu_{\rm c}+\mu_{\rm e}) u_{1,22}-\mu_{\rm c} u_{2,12}
&=
f_1\, ,
\notag
\\
2 \mu_{\rm c} A_{12,1}+(\kappa_{\rm e}-\mu_{\rm c}) u_{1,12}+(\kappa_{\rm e}+\mu_{\rm e}) u_{2,22}-2 \kappa_{\rm e} \omega_{,2}+(\mu_{\rm c}+\mu_{\rm e}) u_{2,11}
&=
f_2 \, ,
\label{eq:eq_eqa_MS_4}
\\
\frac{1}{2} \mu_{\rm M}\,L_{\rm c}^2 \widetilde{a} \left(A_{12,22}+A_{12,11}\right)+\mu_{\rm c} \left(-2 A_{12}+u_{1,2}-u_{2,1}\right)
&=
\frac{M_{12}-M_{21}}{2} \, ,
\notag
\\
\frac{1}{2} \mu_{\rm M}\,L_{\rm c}^2 \widetilde{a} \left(\omega_{,22}+\omega_{,11}\right)-2 (\kappa_{\rm e}+\kappa_{\rm m}) \omega +\kappa_{\rm e} \left(u_{1,1}+u_{2,2}\right)
&=
\frac{M_{11}+M_{22}}{2} \, .
\notag
\end{align}
%
%
%
%
\subsubsection{The isotropic Cosserat model in dislocation form as a particular case of the relaxed micromorphic model}
If we take the limit for $\lambda_{\rm micro},\mu_{\rm micro} \to \infty$ ($\mathbb{C}_{\rm micro} \to \infty$), the isotropic relaxed micromorphic model is particularised to the linear Cosserat model \cite{neff2014unifying,ghiba2023cosserat}. 
The expression of the strain energy for the isotropic Cosserat continuum can be equivalently written in dislocation format as (using Curl as the curvature measure)
\begin{align}
W \left(\text{D}u, A,\text{Curl}\,A\right) = &
\, \mu_{\rm M} \left\lVert \text{sym} \, \text{D}u \right\rVert^{2}
+ \mu_{\rm c} \left\lVert \text{skew} \left(\text{D}u - A \right)\right\rVert^{2}
+ \frac{\lambda_{\rm M}}{2} \text{tr}^2 \left(\text{D}u \right)
\label{eq:energy_Cos}
\\*
&+ \frac{\mu_{\rm M}\, L_{\rm c}^2}{2}
\left(
a_1 \, \left \lVert \text{dev} \, \text{sym} \, \text{Curl} \, A\right \rVert^2 \, 
+ a_2 \, \left \lVert \text{skew} \, \text{Curl} \, A\right \rVert^2 \, 
+ \frac{a_3}{3} \, \text{tr}^2 \left(\text{Curl} \, A \right)
\right)  \, .
\notag
\end{align}
The Cosserat model features the classical displacement filed $u \in \mathbb{R}^{3}$ and the infinitesimal micro-rotation tensor $A \in \mathfrak{so}(3)$, i.e. $A$ is a skew-symmetric second order tensor.
The system of equilibrium equations reads
\begin{align}
\text{Div}
\big[
\overbrace{
2\mu_{\rm e}\,\text{sym}\,\text{D}u
+ 2\mu_{\rm c}\,\text{skew} \left(\text{D}u - A \right)
+ \lambda_{\rm e}\,\text{tr} (\text{D}u) \mathbbm{1}
}^{
\mathlarger{\sigma}\coloneqq
}
\big]
&=
f
\, ,
\notag
\\
\label{eq:equi_Coss}
2\mu_{\rm c}\,\text{skew} \left(\text{D}u - A \right)-\text{skew}\,\text{Curl}
\Big(
\underbrace{
\mu_{\rm M}\, L_{\rm c}^2
\left( a_1 \, \text{dev sym} \, \text{Curl} \, A +
a_2 \, \text{skew} \, \text{Curl} \, A +
\frac{a_3}{3} \, \text{tr} \left(\text{Curl} \, A\right)\mathbbm{1}\right)
}_{
\mathlarger{m}\coloneqq
}
\Big)
&=
\text{skew} \, M
\, .
\end{align}
Here, $\mu_{\rm c}> 0$ is called the Cosserat couple modulus.
The skew-operator in equation (\ref{eq:equi_Coss})$_2$ appears because of the reduced kinematics and $\text{skew} \, M$ is the skew-symmetric part of the body volume moment tensor.
Note that there is \textit{no} equation like $\text{Div}\,\sigma_{\rm micro}=\text{Div}\,\text{skew}\,M$ here and taking $\mu_{\rm c}>0$ is mandatory for coupling both equations in (\ref{eq:equi_Coss}).

Under the plane-strain hypothesis only the in-plane components are different from zero and they only depend on $(x_1,x_2)$.
The structure of the kinematic fields are reported below in  (\ref{eq:gen_kine_Cos})
\begin{align}
u
&
=
\left[
u_{1} \, , \,
u_{2} \, , \,
0
\right]^{\rm T}
\, ,
\qquad\qquad
\text{D}u
=
\left(
\begin{array}{ccc}
u_{1,1} & u_{1,2} & 0 \\
u_{2,1} & u_{2,2} & 0 \\
0 & 0 & 0 \\
\end{array}
\right)
\, ,
\qquad\qquad
\widetilde{A}
=
\left(
\begin{array}{ccc}
0 & A_{12} & 0 \\ 
-A_{12} & 0 & 0 \\ 
0 & 0 & 0
\end{array}
\right)
\, ,
\label{eq:gen_kine_Cos}
\\*[5pt]
\text{Curl} \,\widetilde{A}
&
= 
\left(
\begin{array}{ccc}
0 & 0 & A_{12,1} \\
0 & 0 & A_{21,2} \\
0 & 0 & 0 \\
\end{array}
\right)
\, ,
\qquad
\text{skew} \,\text{Curl} \,\text{Curl} \,\widetilde{A}
= 
\left(
\begin{array}{ccc}
0 &
-(A_{12,11}+A_{12,22}) & 0 \\
A_{12,11}+A_{12,22} &
0 & 0 \\
0 & 0 & 0 \\
\end{array}
\right)
\, .
\notag
\end{align}
Moreover, since 
\begin{align}
\text{tr}(\text{Curl} \,\widetilde{A})=0 \, ,
\quad\text{and}\quad
\left\lVert \text{dev sym} \, \text{Curl} \, \widetilde{A}\right\rVert^2=\left\lVert \text{sym} \, \text{Curl} \, \widetilde{A}\right\rVert^2=\left\lVert \text{skew} \,  \text{Curl} \, \widetilde{A}\right\rVert^2=\frac{1}{2}\left\lVert\text{Curl} \, \widetilde{A}\right\rVert^2 \, ,
\end{align}
under the plane-strain hypothesis, the model will just depend on one cumulative parameter $\widetilde{a}\coloneqq\frac{(a_1+a_2)}{2}$, and the equilibrium equations (\ref{eq:equi_Coss}) reduce to (see the $\sharp$-notation in  (\ref{eq:gen_kine}))
\begin{align}
\text{Div}
\big[
\overbrace{
2\mu_{\rm e}\,\text{sym}\,\text{D}\widetilde{u}^{\sharp}
+ 2\mu_{\rm c}\,\text{skew} \left(\text{D}\widetilde{u}^{\sharp} - \widetilde{A}^{\sharp} \right)
+ \lambda_{\rm e}\,\text{tr} (\text{D}\widetilde{u}^{\sharp}) \mathbbm{1}
}^{
\mathlarger{\sigma}\coloneqq
}
\big]
&=
\widetilde{f}
\, ,
\notag
\\*
2\mu_{\rm c}\,\text{skew} \left(\text{D}\widetilde{u}^{\sharp} - \widetilde{A}^{\sharp} \right)
-
\mu_{\rm M}\, L_{\rm c}^2
\widetilde{a}\,
\text{skew}\,\text{Curl}
\underbrace{ \, \text{Curl}_{\text{2D}} \, \widetilde{A}^{\sharp}
}_{
\mathlarger{m}\coloneqq
}
&=
\text{skew} \, \widetilde{M}
\, .
\end{align}
Note the additional appearance of the skew-operator due to the reduced kinematics of the Cosserat model.
Moreover, the Cosserat model is only operative for positive Cosserat couple modulus $\mu_{\rm c}>0$, in contrast to the relaxed micromorphic model. Finally, the equilibrium equations in component form read
\begin{align}\begin{split}
-2 \mu_{\rm c} A_{12,2}+(\lambda_{\rm e}-\mu_{\rm c}+\mu_{\rm e}) u_{2,12}+(\lambda_{\rm e}+2 \mu_{\rm e}) u_{1,11}+(\mu_{\rm c}+\mu_{\rm e}) u_{1,22}
&=f_1,\\*
2 \mu_{\rm c} A_{12,1}+(\lambda_{\rm e}-\mu_{\rm c}+\mu_{\rm e}) u_{1,12}+(\lambda_{\rm e}+2 \mu_{\rm e}) u_{2,22}+(\mu_{\rm c}+\mu_{\rm e}) u_{2,11}
&=f_2,\\*
\frac{1}{2} \mu_{\rm M}\, L_{\rm c}^2\, \widetilde{a} (A_{12,22}+A_{12,11})+\mu_{\rm c} (-2 A_{12}+u_{1,2}-u_{2,1})
&=\frac{M_{12}-M_{21}}{2}.\end{split}
\end{align}
\subsubsection{Classical isotropic linear elasticity in plane strain}
The plane-strain system of standard classical linear elasticity ($L_{\rm c} \to 0$) reads
\begin{align}
\text{Div}
\big[
\overbrace{
2\mu_{\rm eM}\,\text{sym}\,\text{D}\widetilde{u}^{\sharp}
 + \lambda_{\rm M}\,\text{tr} (\text{D}\widetilde{u}^{\sharp}) \mathbbm{1}
}^{
\mathlarger{\sigma}\coloneqq
}
\big]
&=
\widetilde{f},\\
\intertext{and the component form is}
(\lambda_{\rm M}+\mu_{\rm M}) u_{2,12}+(\lambda_{\rm M}+2 \mu_{\rm M}) u_{1,11}+\mu_{\rm M} u_{1,22}
&=f_1,\notag\\*
(\lambda_{\rm M}+\mu_{\rm M}) u_{1,12}+(\lambda_{\rm M}+2 \mu_{\rm M}) u_{2,22}+\mu_{\rm M} u_{2,11}
&=f_2,\notag
\end{align}
The Fourier system in this case assumes the well-known form
\begin{equation}
\begin{aligned}
-\left((\lambda_{\rm M}+2\mu_{\rm M})\xi_1^2+\mu_{\rm M} \xi_2^2\right)\widehat{u}_1-(\lambda_{\rm M}+\mu_{\rm M})\xi_1\xi_2 \widehat{u}_2=\widehat{f}_1,\\
-(\lambda_{\rm M}+\mu_{\rm M})\xi_1\xi_2 \widehat{u}_1-\left((\lambda_{\rm M}+2\mu_{\rm M})\xi_2^2+\mu_{\rm M}\xi_1^2\right)\widehat{u}_2=\widehat{f}_2,
\end{aligned}
\end{equation}
and the Fourier determinant becomes
\begin{equation}
\det \mathbb{A}_{\rm lin.elast}(\xi)=\mu_{\rm M}(\lambda_{\rm M}+2\mu_{\rm M})\upxi^4.
\end{equation}

\subsection{Properties of the second kind modified Bessel functions}
\label{Bessel}
Here we show some well known relations regarding the second kind modified Bessel functions $K_n[z]$ that have been used in the derivation of the Green's functions in \eqref{Green} and \eqref{CMu} of the relaxed micromorphic medium. Also we derive some useful limits that were employed for passing from the general relaxed micromorphic model to other generalized continua.

The modified Bessel functions $K_n[r]$ are solutions of 
the ODE
\begin{equation}
z^2 u''(z)+z u'(z)-(z^2+n^2)u(z)=0.
\end{equation}
Some useful recurrence relations for the second kind modified Bessel functions $K_n[r]$ are \cite{gradshteyn2014table}:
\begin{equation}
K_{n+1}[z]=K_{n-1}[z]+\frac{2n}{z}K_{n}[z], \qquad  K_n[z]=K_{-n}[z], \quad n \ge 0
\end{equation}
If $z=(x_1^2+x_2^2)^{1/2}>0$, we derive the first and second derivatives of $K_n[z]$ w.r.t $x_i$ as
\begin{equation}
\begin{cases}
&\partial_{x_i} K_n[z]=-\frac{x_i}{2z}\left( K_{n+1}[z]+K_{n-1}[z]\right)\\
&\partial_{x_i} \partial_{x_j} K_n[z]=\frac{x_i x_j}{4z^2}(K_{n+2}[z]+2K_n[z]+K_{n-2}[z])-\frac{1}{2z}\left(\delta_{ij}-\frac{x_i x_j}{z^2} \right)\left(K_{n+1}[z]+K_{n-1}[z]\right)
\end{cases}, \qquad n \ge 0.
\end{equation}
where $\delta_{ij}$ is the Kronecker delta. These equations have been employed for the derivation of the Green's functions of the relaxed micromorphic plane strain theory.

For small argument $z \to 0$ we have the asymptotic relation \cite{gradshteyn2014table}:
\begin{equation}
\label{asymptsmall}
K_n[z] \sim
\begin{cases}
-\ln\tfrac{z}{2}-b, \enskip &\text{for} \enskip n=0\,,\\
\tfrac{\Gamma[n]}{2} \left(\tfrac{2}{z}\right)^n \enskip 
&\text{for} \enskip n>0\,,
\end{cases}    
\end{equation}
where $b$ is the Euler constant and $\Gamma[\cdot]$ is the Gamma function.

For large argument $z \to \infty$ we have the asymptotic relation \cite{gradshteyn2014table}:
\begin{equation}
\label{asymptlarge}
K_n[z] \sim \sqrt{\frac{\pi}{2z}}\,e^{-z} \quad \text{for} \enskip n \ge 0 \, ,
\end{equation}
which show that all $K_n$ functions become quickly zero at infinity with exponential rate.

We now prove some limits that appear in the main text. 
\begin{equation}
\lim_{z \to 0} \left(\frac{2}{z^2}-K_2[z]\right)=\frac{1}{2}\, , \quad \lim_{z \to 0} \left(\frac{1}{z}-K_1[z]\right)=0\, ,  \quad \lim_{z \to 0} z K_1[a z]=a^{-1}\, ,  \quad \lim_{z \to 0} z^2 K_0[z]=0 \, ,  \quad \lim_{z \to 0} K_0[z]=-\ln z.
\end{equation}
Now the first three limits are easily derived by expanding $K_2[z]$ and $K_1[z]$ in series as $z \to 0$. We have: $K_2[z]=2/z^2-1/2+O(z^2)$ and  $K_1[z]=1/z+O(z)$. The last limit is a direct consequence of\eqref{asymptsmall} and the fact that $\lim_{z \to 0} z^n \ln z=0$, $n>0$. The above results cover the limit cases \eqref{limitpure}, \eqref{limitpure2}, \eqref{limclCM1}  where $\ell_2 \to \infty$ or $\mu_{\rm c}=0$.

Accordingly, we have
\begin{equation}
\lim_{z \to \infty} z^2 \,K_0[z]=0 \, , \qquad \lim_{z \to \infty} z\, K_1[z]=0, \qquad \lim_{z \to \infty} \left(\frac{2}{z^2}-K_2[z]\right)=0\, ,
\end{equation}
which are direct consequence of \eqref{asymptlarge}. The above results cover the limit cases \eqref{lmclF}, \eqref{limclCM} where $\ell_j \to 0$.

\end{alphasection}
}

\end{document}